\documentclass{amsart}
\usepackage{a4wide}
\usepackage{comment}
\usepackage{microtype}
\usepackage{stackengine}
\usepackage{mathtools}
\stackMath
\usepackage{amscd}

\usepackage{amsmath,amssymb,amsfonts, accents, stmaryrd}
\usepackage{xstring}
\usepackage{xcolor}

\usepackage{eucal}
\usepackage{url}
\usepackage{latexsym}
\usepackage{amsthm}
\usepackage[latin1]{inputenc}
\usepackage[shortlabels]{enumitem}

\usepackage[colorlinks=true, linkcolor={blue!50!black},
citecolor=blue, pdfhighlight=/O, 
ocgcolorlinks=true]{hyperref}
\hypersetup{bookmarksdepth=2}

\usepackage{graphicx}
\usepackage{color}

\usepackage{tikz-cd}
\tikzset{Isom/.style={above,every to/.append style={edge node={node [sloped, allow upside down, auto=false]{$\sim$}}}}}

\usepackage[normalem]{ulem}

\usepackage[capitalise]{cleveref}
\crefformat{equation}{(#2#1#3)}

\newtheoremstyle{introthms}
	{}{}{\itshape}{}{\bfseries }{}{ }
	{\thmname{#1} \thmnumber{#2}. \thmnote{\bfseries{(#3)}}}

\theoremstyle{introthms}
\newtheorem{introthm}{Theorem}

\newtheorem{introcor}[introthm]{Corollary}

\swapnumbers

\numberwithin{equation}{section}
\theoremstyle{plain}
\newtheorem{theorem}[equation]{Theorem}

\newtheorem{lemma}[equation]{Lemma}

\newtheorem{proposition}[equation]{Proposition}

\newtheorem{corollary}[equation]{Corollary}

\theoremstyle{definition}

\newtheorem{definition}[equation]{Definition}

\newtheorem{construction}[equation]{Construction}
\newtheorem{example}[equation]{Example}

\newtheorem{observation}[equation]{Observation}
\newtheorem{remark}[equation]{Remark}

\newtheorem{notation}[equation]{Notation}

\newtheorem*{uremark}{Remark}

\setenumerate{label=(\arabic*),itemsep=3pt,topsep=3pt, leftmargin=*}

\newcommand{\ie}{i.e.\@}

\newcommand{\lortho}{curved orthofibration}
\newcommand{\lorthos}{curved orthofibrations}

\newcommand{\icat}{$\infty$-category}
\newcommand{\icats}{$\infty$-categories}
\newcommand{\igpds}{$\infty$-groupoids}
\newcommand{\igpd}{$\infty$-groupoid}
\newcommand{\itcat}{$(\infty,2)$-category}
\newcommand{\itcats}{$(\infty,2)$-categories}

\newcommand{\C}{\cat{C}}
\newcommand{\D}{\cat{D}}
\newcommand{\E}{\cat{E}}

\newcommand{\CC}[0]{\cat{C}}
\newcommand{\DD}[0]{\cat{D}}

\newcommand{\Del}[0]{\mb{\Delta}}

\newcommand{\Cat}{\cat{Cat}}
\newcommand{\Grayop}{\mathrm{OpGray}}
\newcommand{\Ortholax}{\mathrm{CrvOrtho}}

\newcommand{\localortho}{\Ortholax}

\DeclareMathOperator{\LCocart}{LCocart}
\DeclareMathOperator{\RCocart}{RCocart}

\DeclareMathOperator{\LCart}{LCart}
\newcommand{\sGpd}{\mathrm{s}\Gpd}
\newcommand{\CAT}{\textsc{Cat}}
\newcommand{\tCat}{\textbf{\textup{Cat}}}

\newcommand{\sSet}[0]{\cat{sSet}}
\newcommand{\Gpd}{\cat{Gpd}}
\newcommand{\sS}[0]{\Gpd}
\newcommand{\GPD}{\textsc{Gpd}}

\newcommand{\Cocart}{\cocart}
\newcommand{\Cart}{\cart}

\DeclareMathOperator{\ortho}{Ortho}
\DeclareMathOperator{\bifib}{Bifib}
\DeclareMathOperator{\cocart}{Cocart}
\DeclareMathOperator{\cart}{Cart}

\renewcommand{\S}{\mathrm{S}}

\newcommand{\Cartlax}{\cart^{\mathrm{opl}}}

\newcommand{\Leftfib}{\mathrm{LFib}}

\newcommand{\markSet}{\sSet^{+}}
\newcommand{\markCat}{\Cat(\markSet)}
\newcommand{\scSet}{\sSet^{\scale}}

\newcommand{\Gray}{\mathrm{Gray}}

\newcommand{\ssS}{\mathrm{s}\Gpd}

\newcommand{\AdTrip}{\cat{AdTrip}}

\newcommand{\Strco}{\mathrm{Str}^{\mathrm{cc}}}
\newcommand{\Strbi}{\mathrm{Str}^{\mathrm{bi}}}
\newcommand{\Strcart}{\mathrm{Str}^{\mathrm{ct}}}
\newcommand{\Strort}{\mathrm{Str}^{\mathrm{oc}}}

\newcommand{\Unco}{\mathrm{Un}^{\mathrm{cc}}}
\newcommand{\Uncart}{\mathrm{Un}^{\mathrm{ct}}}
\newcommand{\Unort}{\mathrm{Un}^{\mathrm{oc}}}

\newcommand{\Dualco}{\mathrm{D}^\mathrm{cc}}
\newcommand{\Dualcart}{\mathrm{D}^\mathrm{ct}}
\newcommand{\SDualco}{\mathrm{S}\Dualco}
\newcommand{\SDualcart}{\mathrm{S}\Dualcart}
\newcommand{\Funcat}{\mathcal{F}} 
\newcommand{\asscat}{\mathrm{ac}}

\DeclareMathOperator{\cst}{cst}

\newcommand{\Fun}[0]{\cat{Fun}}
\newcommand{\core}{\iota}

\newcommand{\nerve}{\mathrm N}
\newcommand{\Q}{\mathrm Q}

\newcommand{\ADual}{{\Span^\perp}}
\newcommand{\AADual}{{\Span^{\perp 2}}}

\newcommand{\Span}{\mathrm{Span}}

\DeclareMathOperator{\Tw}{Tw}

\DeclareMathOperator{\Ar}{Ar}
\newcommand{\TTw}{\Tw^{(2)}}
\newcommand{\TwL}{\Tw^{\ell}}
\newcommand{\TwR}{\Tw^{r}}

\newcommand{\op}[0]{\mm{op}}
\newcommand{\eg}{\mm{eg}}
\newcommand{\fb}{\mm{fw}}
\newcommand{\glob}{\int}
\newcommand{\rev}{\mathrm{rev}}  
\newcommand{\ing}{\mm{in}}

\newcommand{\lax}{\mathrm{lax}}
\newcommand{\oplax}{\mathrm{opl}}

\newcommand{\scale}{\mathrm{sc}}

\newcommand{\lcl}{\emph{l}}    
\newcommand{\rcl}{\emph{r}}
\newcommand{\eff}{\mskip-2mu\textit{eff}}

\newcommand{\rto}[1]{\stackrel{#1}{\rt}}
\newcommand{\lto}[2][\;]{\xleftarrow[#1]{#2}}
\newcommand{\rt}[0]{\longrightarrow}
\newcommand{\lt}[0]{\longleftarrow}

\newcommand{\mm}[1]{\mathrm{#1}}
\newcommand{\mb}[1]{\mathbf{#1}}

\newcommand{\mf}[1]{\mathfrak{#1}}

\newcommand{\ul}[1]{\underline{#1}}
\newcommand{\ol}[1]{\overline{#1}}

\newcommand{\cat}[1]{
\StrLen{#1}[\mystrlen]
\ifnum\mystrlen=1 {#1}
\else \mathrm{#1}
\fi}

\newcommand{\tcat}[1]{\mathbf{#1}}

\newcommand{\ho}{\mathrm{h}}
\newcommand{\Hom}[0]{\mm{Hom}}
\newcommand{\Map}[0]{\mm{Map}}

\newcommand{\colim}{\operatornamewithlimits{\mathrm{colim}}}

\newcommand{\cocolon}{\nobreak \mskip6mu plus1mu \mathpunct{}\nonscript\mkern-\thinmuskip {:}\mskip2mu \relax}

\newcommand{\id}[0]{\mm{id}}
\newcommand{\pr}{\mathrm{pr}}
\newcommand{\const}{\cst}
\newcommand{\ev}{\mathrm{ev}}

\newcommand{\Gtimes}{\boxtimes}
\newcommand{\join}{\mathbin{\operatornamewithlimits{\star}}}
\newcommand{\oplslice}{/\!/^\oplax \,}
\newcommand{\laxslice}{\sslash}

\newcommand{\Aut}{\mathrm{Aut}}

\usetikzlibrary{calc}
\tikzset{curve/.style={settings={#1},to path={(\tikztostart)
			.. controls ($(\tikztostart)!\pv{pos}!(\tikztotarget)!\pv{height}!270:(\tikztotarget)$)
			and ($(\tikztostart)!1-\pv{pos}!(\tikztotarget)!\pv{height}!270:(\tikztotarget)$)
			.. (\tikztotarget)\tikztonodes}},
	settings/.code={\tikzset{quiver/.cd,#1}
		\def\pv##1{\pgfkeysvalueof{/tikz/quiver/##1}}},
	quiver/.cd,pos/.initial=0.35,height/.initial=0}

\DeclareFontFamily{U}{cbgreek}{}
\DeclareFontShape{U}{cbgreek}{m}{n}{
        <-6>    grmn0500
        <6-7>   grmn0600
        <7-8>   grmn0700
        <8-9>   grmn0800
        <9-10>  grmn0900
        <10-12> grmn1000
        <12-17> grmn1200
        <17->   grmn1728
      }{}
\DeclareFontShape{U}{cbgreek}{bx}{n}{
        <-6>    grxn0500
        <6-7>   grxn0600
        <7-8>   grxn0700
        <8-9>   grxn0800
        <9-10>  grxn0900
        <10-12> grxn1000
        <12-17> grxn1200
        <17->   grxn1728
      }{}

\DeclareRobustCommand{\Qoppa}{%
  \text{\usefont{U}{cbgreek}{\normalorbold}{n}\symbol{21}}%
}
\makeatletter
\newcommand{\normalorbold}{%
  \ifnum\pdf@strcmp{\math@version}{bold}=\z@ bx\else m\fi
}
\makeatother

\newcommand{\QF}{\Qoppa}

\title{Two-variable fibrations, factorisation systems and $\infty$-categories of spans}

\author[R.~Haugseng]{Rune Haugseng}
\address{Institutt for matematiske fag, NTNU Trondheim, Norway}
\email{rune.haugseng@ntnu.no}
\author[F.~Hebestreit]{Fabian Hebestreit}
\address{Department of Mathematics, University of Aberdeen, Scotland}
\email{fabian.hebestreit@abdn.ac.uk}
\author[S.~Linskens]{Sil Linskens}\address{Mathematisches Institut, RFWU Bonn, Germany}
\email{linskens@math.uni-bonn.de}
\author[J.~Nuiten]{Joost Nuiten}
\address{IMT, Universit\'e de Toulouse III, France}
\email{joost.nuiten@math.univ-toulouse.fr}

\begin{document}

 \begin{abstract}
We prove a universal property for $\infty$-categories of spans in the
generality of Barwick's adequate triples, explicitly describe the cocartesian fibration corresponding to the span functor, and show that the latter restricts to a self-equivalence on the class of orthogonal adequate triples, which we introduce for this purpose. 

As applications of the machinery we develop we give a quick proof of
Barwick's unfurling theorem, show that an orthogonal factorisation
system arises from a cartesian fibration if and only if it forms an
adequate triple (generalising work of Lanari), extend the description
of dual (co)cartesian fibrations by Barwick, Glasman and Nardin to
two-variable fibrations, explicitly describe parametrised adjoints
(extending work of Torii), identify the orthofibration classifying the
mapping category functor of an $(\infty,2)$-category (building on work
of Abell\'an Garc{\'\i}a and Stern), formally identify the
unstraightenings of the identity functor on the
$\infty$-category of $\infty$-categories with the (op)lax under-categories of a point, and deduce a certain naturality property of the Yoneda embedding (answering a question of Clausen).
\end{abstract}

\setcounter{tocdepth}{1}
\maketitle
\tableofcontents

\section{Introduction}\label{sec:intro}
Categories of spans typically appear in algebra and geometry as a convenient way to encode co-and contravariant functorialities (transfer and restriction) and base change isomorphisms between these. For example, span categories of finite $G$-sets classically appear in representation theory as the indexing categories of Mackey functors \cite{Dress}. In the $\infty$-categorical setting, spans play an even more significant role, because they often provide the only feasible method to organise base change isomorphisms in a homotopy-coherent way. As such, they naturally appear in equivariant homotopy theory (via the theory of spectral Mackey functors \cite{Barwick, Nardin}), in motivic homotopy theory \cite{BachH} and in higher algebra \cite{Cranch, Harpaz, HaugsengElmanto}.

In each of these cases, one is interested in diagrams indexed by a certain \icat{} of spans that is constructed informally as follows: from the data of an \icat{} $X$, together with two classes of maps called \emph{ingressive} (or forwards, denoted $\rightarrowtail$) and \emph{egressive} (or backwards, denoted $\twoheadrightarrow$), one constructs an \icat{} $\Span(X)$ with the same objects as $X$ and with morphisms from $x$ to $y$ given by spans
$$
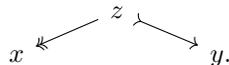
\begin{tikzcd}[row sep=0.4pc]
& z\arrow[ld, two heads]\arrow[rd, rightarrowtail] \\
x & & y.
\end{tikzcd}$$
Composition proceeds by pullback of spans, which of course requires that pullbacks of ingressive and egressive maps exist and remain ingressive and egressive. 
From a geometric point of view, one can think of $\Span(X)$ as an \icat{} of (combinatorial) bordisms in $X^{\op}$, the two legs of a span giving the two boundary inclusions. As such, span categories also arise in various situations as objects of interest in themselves. Notably, they appear in algebraic $K$-theory \cite{BarwickQ,BGT} and its hermitian refinements \cite{Schlichting,CDH2} via Quillen's $\mathcal{Q}$-construction \cite{QuillenHigherK}; see in particular \cite{RaptisSteimle2, CDH2, HS} for works on $K$-theory where such a geometric perspective on spans is brought to fruition.

The above informal description of $\Span(X)$ has been substantiated by Barwick, who constructs an explicit functor sending each \emph{adequate triple} (of an \icat{} $X$ and two classes of maps) to a quasi-category $\Span(X)$ \cite{Barwick}. The first purpose of this text is to study some further abstract properties of this construction of span \icats. For example, we give a description of span \icats{} by a universal property:
\begin{introthm}\label{thma}
There is an adjoint pair $\begin{tikzcd}[cramped, column sep=1.5pc]\TwR \colon \Cat
  \arrow[r, yshift=0.5ex] & \arrow[l, yshift=-0.5ex] \AdTrip \cocolon
  \Span \end{tikzcd}$, where $\TwR(A)$ is the twisted arrow \icat{} of
$A$ with the ingressive and egressive morphisms being those which induce equivalences in the target and source component,
respectively.
\end{introthm}
Here we follow the same notational conventions as in \cite{PartI} and
denote the \icat{} of small \icats{} by $\Cat$ (and similarly that of
\igpds{} by $\Gpd$, etc.).
Let us point out that our proof of this result relies on a Segal space construction of span \icats{}, instead of the point-set approach of Barwick; a slightly different proof was sketched by Raskin in \cite{Raskin}. \Cref{thma} also has a $2$-categorical upgrade, which identifies the diagram \icats{} $\Fun(A,\Span(X))$ with a span \icat{} of the \icat{} of diagrams $\TwR(A) \rightarrow X$. This statement has already been applied in hermitian K-theory, where it forms the basis for the parametrised algebraic surgery of \cite{HS}. 

In addition, we give a description of the `universal span category', i.e.\ the cocartesian fibration classified by the functor $\Span\colon \AdTrip\rt \Cat$:
\begin{introthm}\label{thmb}
The \emph{opposite} of the lax under-category $\ast\laxslice\AdTrip$ admits the structure of an adequate triple in which a map
$$\begin{tikzcd}
& \ast\arrow[rd, "y", ""{name=s, swap}]\arrow[ld, "x"{swap}]\\
X\arrow[rr, "f"{swap}]\arrow[Rightarrow, from=s, shorten <=1ex, shorten >=5ex, "\mu"{pos=0.2, swap}] & & Y
\end{tikzcd}$$
is egressive if $\mu\colon y \to f(x)$ is egressive, and
  ingressive if $\mu$ is ingressive and $f$ is an equivalence.
The natural map $\Span\big((\ast\laxslice\AdTrip)^{\op}\big)\rt \AdTrip$ extracting the morphisms  underlying egressive maps (and inverting those underlying ingressive ones) is then equivalent to the cocartesian fibration classified by $\Span\colon \AdTrip\rt \Cat$.
\end{introthm}

As a consequence of \Cref{thmb}, we obtain a new proof of the main result of \cite{BGN}, which asserts that the cocartesian unstraightening of a functor $F\colon B\rightarrow \Cat$ can be constructed from its cartesian unstraightening via a span construction.

Note that the \icat{} $\Span(X)$ has the following features: the egressive arrows from $X$ are reverted, the ingressive arrows from $X$ are kept and each map factors uniquely as a `reverted-egressive', followed by an ingressive map. This suggests that taking span \icats{} acts as an involution on adequate triples arising from (orthogonal) factorisation systems. The second part of the paper studies such \emph{orthogonal} adequate triples, whose egressive and ingressive arrows form a factorisation system and in which all commuting squares of the form
\begin{equation}\label{diag:ambigressive}\begin{tikzcd}
x'\arrow[r, rightarrowtail]\arrow[d, two heads] & x\arrow[d, two heads]\\
y'\arrow[r, rightarrowtail] & y
\end{tikzcd}\end{equation}
are cartesian. Our most significant result about these is the following: 
\begin{introthm}\label{thmc}
Taking span \icats{} gives rise to a $C_2$-action on the full subcategory of orthogonal adequate triples
\[\ADual \colon \AdTrip^\perp \longrightarrow \AdTrip^\perp.\]
\end{introthm}
The simplest example of an orthogonal adequate triple is given by a product $A\times B$, with ingressive maps coming from $A$ and egressive maps coming from $B$. More generally, every cartesian fibration $p\colon X\rt A$ defines the structure of an orthogonal adequate triple on $X$, for which a map is ingressive if it is $p$-cartesian and egressive if it is fibrewise, i.e.\ its image under $p$ is invertible. In fact, every orthogonal adequate triple arises in this way:
\begin{introthm}\label{thmd}
An orthogonal factorisation system arises from a cartesian fibration as the cartesian and fibrewise maps if and only if it forms an orthogonal adequate triple. Moreover, this construction restricts to an equivalence between the $\infty$-category of orthogonal adequate triples and the $\infty$-category of cartesian fibrations whose fibres have contractible realisation.
\end{introthm}

This result extends previous work of Lanari \cite{LanariCart}, who established such a correspondence for pointed cartesian fibrations (i.e.\ those whose base and total \icat{} have a terminal object which is preserved by the fibration) and what he calls cartesian factorisation systems.
Using the relation between orthogonal adequate triples and cartesian
fibrations, the span duality from \Cref{thmc} also gives an
auto-equivalence of the \icat{} of cartesian fibrations. This reduces
to the duality functor of Barwick, Glasman and Nardin from \cite{BGN}.

Similarly, one can use \Cref{thmc} to obtain duality functors for fibrations over a product $A\times B$. For example, one can consider maps of orthogonal triples $p=(p_1, p_2)\colon X\rt A\times B$ such that:
\begin{enumerate}
\item the ingressive maps of $X$ are exactly the $p$-cartesian lifts of arrows in $A$
\item $p_2$ is a cocartesian fibration and for any square \eqref{diag:ambigressive} in $X$ where $x\to y$ is $p_1$-cocartesian, the base change $x'\rt y'$ is $p_1$-cocartesian as well.
\end{enumerate}
Such maps are called \emph{orthofibrations} in \cite{PartI}, or \emph{two-sided fibrations} in \cite[Section 7.1]{RiehlVerity}, and can be
straightened to functors $A^{\op}\times B\rt \Cat$, which also
correspond to cartesian fibrations over $A \times B^\op$ via
unstraightening. Using that $\ADual$ sends $A\times B$ (with the adequate triple structure mentioned above) to $A\times B^{\op}$, the dualisation of Theorem \ref{thmc} also restricts to an equivalence
\begin{equation}\label{eq:span ortho dual}
\ADual\colon \ortho(A, B)\simeq \Cart(A\times B^{\op})
\end{equation}
between orthofibrations and cartesian fibrations, taking opposite categories at the level of fibres. There are similar dualities for the various other types of two-variable fibrations considered in \cite{PartI}; For example, Theorem \ref{thmc} also restricts to a duality between \lorthos{} over $A\times B$ and op-Gray fibrations over $A\times B^{\op}$ (see \cite{PartI} or Section \ref{sec:product} for definitions), which under straightening \cite{LurieGoo} correspond to $2$-functors out of the Gray tensor product $B\boxtimes A^{\op}$. 

Relying on To\"en's equivalence $\Aut(\Cat) \simeq C_2$ \cite{Toen}, we extend another result from \cite{BGN} and show that any composite of dualisations and (un)straightenings as above is uniquely determined by its action on fibres:

\begin{introthm}\label{thme}
The three functors $\Cat^\op \times \Cat^\op \rightarrow \textsc{Cat}$ given by
\[(A,B) \quad \longmapsto \quad  \Fun(A \times B,\Gpd),\quad \ \Fun(A \times B, \Cat) \quad \text{and} \quad \Fun(A \Gtimes B,\mathbf{Cat})\]
have automorphism groups $*, C_2$, and $\ast$, respectively. The non-trivial automorphism in the middle case is given by post-composition with $(-)^\op \colon \Cat \rightarrow \Cat$. (Here $\mathbf{Cat}$ denotes the $(\infty,2)$-category of \icats{}, as
in \cite{PartI}.)
\end{introthm}

For example, this implies that the dualisation equivalence \eqref{eq:span ortho dual} arising from \Cref{thmc} coincides with the dualisation constructed in \cite{PartI} using straightening and unstraightening (though this can also be deduced from \Cref{thmb} directly). 
As an application of this explicit dualisation for two-variable fibrations, we consider the fibrations classified by the enhanced mapping functor
\[\Map_{\mathbf{X}}\colon X^\op\times X\rightarrow \Cat\]
associated to an $(\infty, 2)$-category $\mb{X}$ with underlying \icat{} $X$. On the one hand, this functor is classified by a cartesian fibration
\[(s,t)\colon\TwR(\mathbf{X})\rt X\times X^\op\]
from the \emph{oplax twisted arrow \icat{}} of $\mb{X}$, explicitly constructed in \cite{GS}. We show (\Cref{thm:CartDualAr}) that dually (and upon taking opposite categories), this functor also classifies the orthofibration
\[(s, t)\colon \Ar^{\oplax}(\mathbf{X})\rt X\times X\]
from the \icat{} underlying the \emph{oplax arrow \icat{}}, defined using the (in principle unrelated) Gray tensor product of Gagna, Harpaz and Lanari from \cite{GHL-Gray}. This is a typical example of an identification of dual fibrations which seems difficult to see by passing through a form of straightening and unstraightening (as in \cite{PartI}).

In a different direction, \Cref{thme} shows that straightening a cartesian fibration over $A\times B$ to a functor $A^{\op}\times B^{\op}\rt \Cat$ is naturally equivalent to straightening it first over $A$ and then over $B$. This should not be surprising, but deducing it from the definitions does not seem entirely obvious (nevertheless, we provide a second proof as \ref{cor:strintwovar}). As a consequence, we answer a recent question of Clausen:
\begin{introcor}\label{corf}
The Yoneda embedding $\cat{C} \hookrightarrow \mathcal P(\cat{C})$ canonically extends to a natural transformation of functors $\Cat \rt \CAT$  from the inclusion to the composite
\[\Cat \xrightarrow{\Fun(-,\Gpd)} (\CAT^\mathrm R)^\op \simeq \CAT^\mathrm{L} \subseteq \CAT.\]
\end{introcor}
Again, this result is certainly expected, and we were surprised to learn that it is apparently not contained in the literature. Let us point out that, essentially by definition, the Yoneda embedding does define a natural transformation to the functor $\mathcal{P}\colon \Cat \rt \CAT$ defined using that $\mathcal P(\cat{C})$ is the free cocompletion of $\cat{C}$. However, it is not a priori clear that this functor agrees with the one appearing in the corollary in a fashion compatible with the Yoneda embedding (but this follows from the result above as well).

\begin{uremark}
Along with \cite{PartI} this article is part of a recombination of our earlier preprints \cite{Haugseng} and \cite{HLN}. In \cite{PartI}, we described the dualisation of two-variable fibrations, as well as its applications to monoidal and parametrised adjunctions, along the lines of  \cite{Haugseng}; These applications were also discussed in terms of the dualisation functor $\ADual$ in \cite{HLN} (in particular Proposition A and Theorem C in there). The results from \cite{HLN} concerning span \icats{} and the comparison of various straightening functors for two-variable fibrations (in particular Theorem B and Corollary D there) are contained in the present paper.
\end{uremark}

\subsection*{Organisation}
In Section \ref{sec:span} we redevelop Barwick's theory of span \icats{} for adequate triples in terms of complete Segal spaces, and then discuss some further properties of the span construction in \Cref{sec:functors}. In particular, here we prove Theorem \ref{thma} as \Cref{modelfuncob} and deduce Theorem \ref{thmb} as a corollary of \Cref{thm:span_funct_cocart}. 

We then introduce orthogonal adequate triples in Section \ref{sec:orthtrip} and establish Theorem \ref{thmc} as \Cref{D_is_equiv}. The relation between orthogonal adequate triples and fibrations is discussed in Section \ref{sec:fibs}, where \Cref{thmd} is proven as \Cref{prop:orthogonal=cart fib}. In Section \ref{sec:product} we show how \Cref{thmc} induces dualities between various types of two-variable fibrations considered in \cite{PartI}, and use this to identify oplax arrow and twisted arrow \icats{} as duals in Section \ref{sec:OplaxAr}. We deduce \Cref{corf} in the final Section \ref{sec:Yoneda} as \Cref{yonedafunctorial}.

Finally, the uniqueness of all these dualities, i.e.\ Theorem \ref{thme}, is derived in Appendix \ref{sec:straightening} as \Cref{discrete autos}.

\subsection*{Conventions}

In order to declutter notation we will write $\Gpd, \Cat$ and $\Cat_2$ for the \icats{} of $\infty$-groupoids (or spaces), $\infty$-categories and $(\infty, 2)$-categories, respectively.

The letter $\iota$ will denote the core of an
$\infty$-category, \ie{} the $\infty$-groupoid spanned by its
equivalences. By a subcategory of an $\infty$-category we
mean a functor such that the induced morphisms on mapping \igpds{} and
cores are inclusions of path components. A subcategory is \emph{full} if the functor furthermore induces
equivalences on mapping \igpds{}, while it is \emph{wide} if the functor induces an
equivalence on cores. Similarly, a sub-$2$-category of an
$(\infty, 2)$-category is a functor inducing subcategory inclusions
on mapping \icats{} and a monomorphism on underlying
\igpds{}; we say such a sub-$2$-category is \emph{1-full} if it
locally full, \ie{} is
given by full subcategory inclusions on mapping \icats{}.

Throughout, we shall use small caps such as $\CAT$ to indicate large
variants of $\infty$-categories and boldface such as
$\tCat$ to denote $(\infty,2)$-categories. We have
also reserved sub- and superscripts on category names to refer to
changes on morphisms, e.g.\ $\cart(A) \subseteq \Cartlax(A)$.

We will write $\Ar(C)$ for the arrow \icat{} $\Fun([1],C)$ of an \icat{} $C$, and $\TwL(C)$ and $\TwR(C)$ for the two versions of the twisted arrow category, geared so that the combined source-target map defines a left fibraton in the former, and a right fibration in the latter case, see  \Cref{ex:twar}. Finally, we write $\Lambda_0[2]$ for the span category $1\leftarrow 0\rightarrow 2$ and $\Lambda_2[2]$ for the cospan.

\subsection*{Acknowledgments}
It is a pleasure to thank Dustin Clausen, Emanuele Dotto, Yonatan Harpaz, Gijsbert Scheltus Karel Sebastiaan Heuts, Corina Keller, Achim Krause, Markus Land, Denis Nardin, Thomas Nikolaus, Stefan Schwede, Wolfgang Steimle and Lior Yanovski for several very useful discussions. We further thank Shaul Barkan and Jan Steinebrunner for pointing out an omission in the verification of orthogonality for the factorisation system on span categories in the original version of this paper. \\

During the preparation of this text FH and SL were members of the Hausdorff Center for Mathematics at the University of Bonn funded by the German Research Foundation (DFG) and FH was furthermore a member of the cluster ``Mathematics M\"unster: Dynamics-Geometry-Structure'' at the University of M\"unster under grant nos. EXC 2047 and EXC 2044, respectively. FH would also like to thank the Mittag-Leffler Institute for its hospitality during the research program ``Higher algebraic structures in algebra, topology and geometry'', supported by the Swedish Research Council (VR) under grant no.\ 2016-06596. FH and JN were further supported by the European Research Council (ERC) through the grants `Moduli spaces, Manifolds and Arithmetic', grant no. 682922, and `Derived Symplectic Geometry and Applications', grant no. 768679, respectively. SL was supported by the DFG Schwerpunktprogramm 1786 'Homotopy Theory and Algebraic Geometry' (project ID SCHW 860/1-1)\\

\section{Adequate triples and \icats{} of spans}\label{sec:span}

In this section we shall review the theory of adequate triples and their associated span \icats{}, as developed by Barwick in \cite{Barwick} under the name effective Burnside \icats{}. 
We will use the opportunity to present an alternate viewpoint on the material by translating the assertions along the equivalence between $\infty$-categories and complete Segal \igpds{} (i.e.\ complete Segal spaces). This will allow for simpler proofs with far less explicit combinatorics, and will form the basis for our analysis of functors to span \icats{} in the next section.

Since composition in span \icats{} proceeds by pullback, one uses the following data as input for span \icats{}: 
\begin{definition}
  An \emph{adequate triple} $(X, X_{\mathrm{in}}, X^{\mathrm{eg}})$ consists of an $\infty$-category $X$ together with two wide subcategories $X_\mathrm{in}$ and $X^\mathrm{eg}$, whose morphisms are called \emph{ingressive} (or \emph{forwards}, denoted $\rightarrowtail$) and \emph{egressive} (or \emph{backwards}, denoted $\twoheadrightarrow$), respectively, such that
\begin{enumerate}
    \item for any ingressive morphism $f \colon y\rightarrowtail x$ and any egressive morphism $g \colon x'\twoheadrightarrow x$, there exists a pullback
\[\begin{tikzcd}
y' \arrow[d, "g'"{left}] \arrow[r, "f'"] & x' \arrow[d, "g", twoheadrightarrow] \\
y \arrow[r, "f", rightarrowtail]       & x,                    
\end{tikzcd}\]
\item and in any such pullback $f'$ is again ingressive and $g'$ egressive.
\end{enumerate}
Squares whose horizontal arrows are ingressive and whose vertical arrows are egressive are called \emph{ambigressive}, and \emph{ambigressive cartesian} if they are furthermore pullback diagrams. A functor 
\[
F\colon (X,X_\ing,X^\eg)\rightarrow (Y,Y_\ing,Y^\eg)
\] of adequate triples is given by a functor $F\colon X\rightarrow Y$
which preserves ambigressive pullbacks. Therefore we may define $\AdTrip$, the \icat{}
of adequate triples, as the subcategory of $\Fun(\Lambda_2[2],\Cat)$
spanned by the adequate triples and those natural transformations
whose evaluation at $2 \in \Lambda_2[2]$ preserves ambigressive
pullback squares.
\end{definition}

Note that being a natural transformation boils down to preserving ingressive and egressive maps, since these form subcategories. We shall often drop them from notation to avoid cluttering.

\begin{notation}
Given a morphism $p:Y \rightarrow X$ of adequate triples we write
$p_\ing\colon Y_\ing \rightarrow X_\ing$ and $p^\eg\colon Y^\eg
\rightarrow X^\eg$ for the restriction of $p$ to ingressives and
egressives.
\end{notation}

\begin{example}\label{ex:adequate triple}
Let us record the following immediate examples:
\begin{enumerate}
\item For any wide subcategory $T$ of an \icat{} $S$, we have the adequate triples $(S, T, \core{S})$ and $(S, \core{S}, T)$ where the ingressives are maps in $T$ and the egressives are equivalences, and vice versa. In particular, for any \icat{} $S$ we always have the adequate triples $(S,S,\core{S})$ and $(S,\core{S},S)$.
\item If $S$ admits all pullbacks, we have the adequate triple $(S,S,S)$ where all maps are both egressive and ingressive. This gives rise to a fully faithful embedding $\Cat^{\mathrm{pb}}\hookrightarrow \AdTrip$ of the subcategory of $\Cat$ spanned by those \icats{} which admit pullbacks, and those functors which preserve pullbacks.

\item If $(X, X_\ing, X^{\eg})$ is an adequate triple, then $X^\rev=(X, X^{\eg}, X_{\ing})$ is an adequate triple.
\end{enumerate}
\end{example}
\begin{lemma}
The \icat{} of adequate triples admits all limits and one simply computes $\lim\limits \big(X, X_{\ing}, X^{\eg}\big)\simeq \big(\lim\limits X, \lim\limits X_{\ing}, \lim\limits X^{\eg}\big)$ for every $X \colon I \rightarrow \AdTrip$. Likewise, $\AdTrip$ has filtered colimits given by $\colim\limits \big(X, X_{\ing}, X^{\eg}\big)\simeq \big(\colim\limits X, \colim\limits X_{\ing}, \colim\limits X^{\eg}\big)$.
\end{lemma}
\begin{proof}
Since the limit of a natural wide subcategory inclusion remains a wide
subcategory inclusion, $\big(\lim\limits X, \lim\limits X_{\ing},
\lim\limits X^{\eg}\big)$ is indeed a triple of a category together
with two wide subcategories (and likewise for filtered
colimits). Notice that an ambigressive square in $\lim X$ is cartesian
if and only if its image in each $X_i$ is cartesian. Furthermore, a
cospan $x\rightarrowtail x'\twoheadleftarrow x''$ in $\lim X$ admits
an ambigressive pullback if its images in each $X_i$ do (in this case
those pullbacks assemble into an object of the limit). This implies
that $\lim X$ is an adequate triple and that a map $Z\rt \lim X$
preserves ambigressive pullback squares if and only each composite
$Z\rt \lim X\rt X_i$ does. This implies that the (pointwise) limit in
$\Fun(\Lambda_2[2], \Cat)$ also provides the limit in the subcategory
$\AdTrip\hookrightarrow \Fun(\Lambda_2[2], \Cat)$.

A dual argument applies to the colimit of a filtered diagram. In this case, the image of an ambigressive pullback square in some $X_i$ under the map $X_i\rt \colim X_i$ remains an ambigressive pullback (since at the level of mapping $\infty$-groupoids, filtered colimits commute with pullbacks). Conversely, any cospan $x\rightarrowtail x'\twoheadleftarrow x''$ in $\colim X$ arises as the image of a cospan in some $X_i$, whose pullback in $X_i$ then provides the desired ambigressive pullback in $\colim X$. It follows that $\colim X$ is an adequate triple and that a map of triples $\colim X\to Z$ preserves ambigressive pullback squares if and only if each $X_i\to \colim X\to Z$ does. This implies that $\colim X$ is also the colimit in the subcategory $\AdTrip\hookrightarrow \Fun(\Lambda_2[2], \Cat)$.
\end{proof}
\begin{lemma}\label{lem:adtrip_internal_hom}
Let $X$ and $Y$ be two adequate triples and let us write
$\Fun_{\AdTrip}(X,Y)$ for the full subcategory of $\Fun(X, Y)$ spanned
by the morphisms of adequate triples. Say a natural transformation
$\tau \colon F \Rightarrow G$ is ingressive if it is pointwise
ingressive in $Y$ and for every egressive $f \colon x \twoheadrightarrow y$
in $X$ the square
\[\begin{tikzcd}
F(x) \ar[r,"F(f)", two heads] \ar[d,"\tau_x"{swap}, tail] & F(y) \ar[d,"\tau_y", tail]\\
G(x) \ar[r,"G(f)", two heads] & G(y)
\end{tikzcd}\]
is cartesian, and analogously for an egressive natural transformation. This endows $\Fun_{\AdTrip}(X,Y)$ with the structure of an adequate triple, such that the evaluation map $\Fun_{\AdTrip}(X,Y)\times X\to Y$ exhibits it as the internal mapping object in $\AdTrip$. In particular, $\AdTrip$ is cartesian closed.
\end{lemma}
Consequently, one obtains a natural equivalence of adequate triples
\[\Fun_{\AdTrip}(X \times Y,Z) \simeq \Fun_{\AdTrip}(X, \Fun_{\AdTrip}(Y,Z)).\]
\begin{proof}
The pointwise pullback $P$ of an ingressive and an egressive natural transformation $F \twoheadrightarrow H \leftarrowtail G$ exists by assumption on $Y$. To verify that $P$ is indeed a map of adequate triples note that for an ingressive $x \rightarrowtail x'$ the map $P(x) \rightarrow P(x')$ participates in the commutative cube
\[\begin{tikzcd}
       &P(x') & & G(x') \\
P(x) & & G(x) \\ 
& F(x') & & H(x')\\
F(x) && H(x) 
\arrow[from=2-1, to=1-2]
\arrow[tail,from=2-3, to=1-4]
\arrow[tail,from=4-1, to=3-2]
\arrow[tail,from=4-3, to=3-4]
\arrow[two heads,from=2-1, to=2-3]
\arrow[tail,from=2-1, to=4-1]
\arrow[two heads,from=4-1, to=4-3]
\arrow[tail,from=2-3, to=4-3]
\arrow[two heads,from=1-2, to=1-4]
\arrow[tail,from=1-2, to=3-2]
\arrow[two heads,from=3-2, to=3-4]
\arrow[tail,from=1-4, to=3-4]
\end{tikzcd}\]
whose front and back face are cartesian by definition of $P$ and whose bottom face is cartesian since $F \twoheadrightarrow H$ is egressive. Consequently also the top face is cartesian showing that $P(x) \rightarrow P(x')$ is also ingressive. The argument for the preservation of egressives is the same, and that $P$ preserves ambigressives pullbacks follows directly from the pasting laws. 

It is also obvious that the maps $P \rightarrow G$ and $P \rightarrow F$ are pointwise egressive and ingressive, respectively, and the argument above also verifies the second condition in the definition of an egressive transformation, and the argument is dual for the ingressives. In total this shows that $\Fun_{\AdTrip}(X,Y)$ is an adequate triple.

To verify the universal property, let $Z$ be a third adequate triple. In suffices to verify that under the equivalence $\Map_{\Cat}(Z\times X, Y)\simeq \Map_{\Cat}(Z, \Fun(X, Y))$ a functor $F\colon Z\times X\rt Y$ is a map of adequate triples if and only if the corresponding map $F'\colon Z\rt \Fun(X, Y)$ takes values in $\Fun_{\AdTrip}(X, Y)$ and determines a map of adequate triples. To see this, note that any ambigressive square in $Z\times X$ is uniquely a composite of four types of squares: 
\[\begin{tikzcd} (z,x_0) & (z,x_1) && (z_0,x) & (z_1,x) \\
                        (z,x_2) & (z,x_3) && (z_2,x) & (z_3,x) \\
(z_0,x_0) & (z_1,x_0) && (z_0,x_0) & (z_0,x_1) \\
(z_0,x_2) & (z_1,x_2) && (z_2,x_0) & (z_2,x_1)
\arrow[from=1-1, to=1-2, tail]
\arrow[from=1-1, to=2-1, two heads]
\arrow[from=2-1, to=2-2, tail]
\arrow[from=1-2, to=2-2, two heads]
\arrow[from=1-4, to=1-5, tail]
\arrow[from=1-4, to=2-4, two heads]
\arrow[from=2-4, to=2-5, tail]
\arrow[from=1-5, to=2-5, two heads]
\arrow[from=3-1, to=3-2, tail]
\arrow[from=3-1, to=4-1, two heads]
\arrow[from=4-1, to=4-2, tail]
\arrow[from=3-2, to=4-2, two heads]
\arrow[from=3-4, to=3-5, tail]
\arrow[from=3-4, to=4-4, two heads]
\arrow[from=4-4, to=4-5, tail]
\arrow[from=3-5, to=4-5, two heads]
\end{tikzcd}\]
Now, $F$ preserving ambigressive (cartesian) squares of the first type is equivalent to $F'$ taking values in $\Fun_{\AdTrip}(X, Y)$. Preserving the ambigressive (automatically cartesian) squares as in the lower row then corresponds precisely to $F'\colon Z\rt \Fun_{\AdTrip}(X, Y)$ preserving ingressive and egressive morphisms. Finally, $F$ preserves the cartesian ambigressive squares as in the upper right corner if and only if $F'$ preserves ambigressive pullbacks.
\end{proof}
One can generate more interesting examples of adequate triples using the following criterion:
\begin{proposition}\label{prop:adequate triple}
Let $p\colon Y\rt X$ be a functor and let $(X, X_\mm{in},
X^{\mm{eg}})$ be an adequate triple such that $Y$ has all
$p$-cartesian lifts over $X_{\ing}$. Then $Y$ is part of an adequate triple
$\big(Y, Y^\dagger, p^{-1}(X^{\eg}))$, where a map is ingressive if it
is $p$-cartesian and its image is ingressive in $X$, and egressive if
its image under $p$ is. This structure also upgrades $p$ to a map of
adequate triples.
\end{proposition}

\begin{proof}
Consider the left solid cospan in the diagram
\[\begin{tikzcd}
	{y_3} & {y_2} && {x_1\times_{x_0} x_2} & {x_2} \\
	{y_1} & {y_0} && {x_1} & {x_0}
	\arrow[tail, from=2-1, to=2-2]
	\arrow[two heads, from=1-2, to=2-2]
	\arrow[two heads, from=1-5, to=2-5]
	\arrow[tail, from=2-4, to=2-5]
	\arrow[dashed, two heads, from=1-4, to=2-4]
	\arrow[dashed, tail, from=1-4, to=1-5]
	\arrow[dashed, tail, from=1-1, to=1-2]
	\arrow[dashed, two heads, from=1-1, to=2-1]
\end{tikzcd}\] with  $x_1\rightarrowtail x_0\twoheadleftarrow x_2$ its image in $X$ under $p$.
By asgsumption the right diagram can be completed to an ambigressive pullback in $X$ and we let $y_3\rightarrowtail y_2$ be a $p$-cartesian lift of $x_1\times_{x_0} x_2 \rightarrowtail x_2$. The universal property of $p$-cartesian maps then implies that there exists a unique map $y_3\rightarrow y_1$ which lives over $x_1\times_{x_0} x_2\to x_1$ and makes the left square commute. Note that the left square is ambigressive by definition. To see that it is also a pullback, compute 
\begin{align*}
&\ \Hom_Y(y,y_1) \times_{\Hom_Y(y,y_0)}\Hom_Y(y,y_2) \\
\simeq&\ \Hom_X(p(y),x_1) \times_{\Hom_X(p(y),x_0)} \Hom_Y(y,y_0)\times_{\Hom_Y(y,y_0)} \Hom_Y(y,y_2)  \\
\simeq&\ \Hom_X(p(y),x_1) \times_{\Hom_X(p(y),x_0)} \Hom_Y(p(y),x_2)\times_{\Hom_X(p(y),x_2)} \Hom_Y(y,y_2) \\
\simeq&\ \Hom_X(p(y),x_1 \times_{x_0} x_2) \times_{\Hom_X(p(y),x_2)} \Hom_Y(y,y_3) \\
\simeq&\ \Hom_Y(y,y_2)
\end{align*}
for $y \in Y$, where the first equivalence uses the universal property of the cartesian edge $y_3 \rightarrowtail y_2$ and the last one that of $y_1 \rightarrowtail y_0$. Finally note that the uniqueness of pull-backs implies that every ambigressive pullback is of this form and thus are clearly preserved by $p$, which is thus a map of adequate triples.
\end{proof} 

The case of \Cref{prop:adequate triple} where $p$ is a cartesian fibration admits the following generalisation:
\begin{proposition}\label{prop:cartadequateforspan}
Let $(X, X_\ing, X^{\eg})$ be an adequate triple, $F\colon X^{\op}\rightarrow \AdTrip$ a functor and $p\colon Y\rt X$ its cartesian unstraightening. Then $Y$ has the structure of an adequate triple, in which a map $y_1\to y_0$ is egressive if it factors as an egressive morphism in the fibre $Y_{p(y_1)}\simeq F(p(y_1))$, followed by a $p$-cartesian morphism with egressive image in $X$ (and likewise for the ingressives). Furthermore, this structure makes $p$ a map of adequate triples.
\end{proposition}

An equivalent description of the egressives is given by unstraightening the functor $X^\op \rightarrow \AdTrip \xrightarrow{(-)^\eg} \Cat$ and then restricting the arising cartesian fibration to $X^\eg$, and likewise for the ingressives; this follows immediately from the fact that every map in the domain of a cartesian fibration factors uniquely as a fibrewise followed by a cartesian arrow.

\begin{proof}
Consider a cospan $f\colon y_1\rightarrowtail y_0 \twoheadleftarrow y_2\cocolon g$ where the left arrow is ingressive and the right arrow is egressive in $Y$. To show that this admits an ambigressive pullback, note that $f$ decomposes into two egressive maps in $Y$, one contained in a fibre of $p$ and the other $p$-cartesian (and likewise for $y_2\twoheadrightarrow y_0$). By the pasting lemma for pullbacks, it therefore suffice to show the existence of an ambigressive pullback in the arising special cases.
	
First, suppose that both $f$ and $g$ are contained in a single fibre $Y_x\simeq F(x)$. Since the fibre forms an adequate triple, there exists an ambigressive pullback of $f$ and $g$ within $Y_x$. This square remains a pullback square in $Y$ by \cite[Proposition 4.3.1.10]{HTT} and stays ambigressive by definition.

Next, suppose that $g$ is a $p$-cartesian lift of an egressive arrow. Then Proposition \ref{prop:adequate triple} implies the existence of a pullback 
\[\begin{tikzcd}
y_3 & y_2 \\
y_1 & y_0
\arrow[tail, from=2-1, to=2-2,"f"]
\arrow[two heads, from=1-2, to=2-2,"g"]
\arrow[dashed, from=1-1, to=1-2]
\arrow[two heads, dashed, from=1-1, to=2-1]
\end{tikzcd}\]
that maps to an ambigressive pullback in $X$. It also implies that the left vertical arrow is again $p$-cartesian and that the upper horizontal one is $p$-cartesian or fibrewise if $f$ is (the former by switching the roles of ingressives and egressives in \ref{prop:adequate triple}. In both cases it follows that the square is ambigressive.

The final case, i.e. if $f$ is assumed a $p$-cartesian lift of an ingressive arrow follows in the same fashion by again reversing the roles of the egressives and ingressives in Proposition \ref{prop:adequate triple}.
\end{proof}

\begin{example}\label{ex:twar}
Let $(X, X_{\ing}, X^{\eg})$ be an adequate triple and let $\TwR(X)$ be its twisted arrow \icat{}, with the convention that the source and target define a right fibration. In particular, $s\colon \TwR(X)\rt X$ is a cartesian fibration and Proposition \ref{prop:adequate triple} endows $\TwR(X)$ with the structure of an adequate triple, in which a map $\alpha\rt \beta$ given by
$$\begin{tikzcd}
x \arrow[d, "\alpha"{swap}] \arrow[r, "f"] & y \arrow[d, "\beta"]\\
w & z \arrow[l, "g"] 
\end{tikzcd}$$ is egressive if $f$ is egressive, and ingressive if $f$ is ingressive and $g$ is an equivalence.
\end{example}

\begin{example}\label{ex:twar simplex}
Let $A$ be an \icat{}, considered as an adequate triple $(A,A,\iota{A})$ with all morphisms ingressive and only the equivalences egressive. Example \ref{ex:twar} endows $\TwR(A)$ with the structure of an adequate triple, where a morphism is ingressive if its image under $t\colon \TwR(A)\rt A^{\op}$ is an equivalence, and egressive if its image under $s\colon \TwR(A)\rt A$ is an equivalence. This determines a functor $\TwR \colon \Cat\rt \AdTrip$.

For the construction of span $\infty$-categories this structure will be of particular interest when $A$ is a simplex. For example, $\TwR([2])$ is the poset
\[
\begin{tikzcd}
[row sep=2.4ex, column sep=0.8ex]
 & & (0 \leq 2) \ar[ld, twoheadrightarrow]\ar[rd, rightarrowtail]& & \\
 & (0 \leq 1) \ar[ld, twoheadrightarrow]\ar[rd, rightarrowtail]& & (1 \leq 2) \ar[ld, twoheadrightarrow]\ar[rd, rightarrowtail]& \\
(0 \leq 0) && (1 \leq 1) && (2 \leq 2)
\end{tikzcd}
\]
where the left pointing arrows are egressive and the right pointing ones are ingressive. More generally, the ambigressive pullbacks in $\TwR([n])$ are precisely given by the diagrams
\[\begin{tikzcd}[column sep=0.8pc, row sep=0.8pc]
& (i \leq l) \arrow[dr, rightarrowtail]\arrow[dl, twoheadrightarrow] & \\
(i \leq j) \arrow[dr, rightarrowtail]& &\arrow[dl, twoheadrightarrow] (k \leq l)\\
& (k \leq j). & 
\end{tikzcd}\]
for $i \leq k \leq j \leq l$. 
\end{example}

We next set out to construct Barwick's functor
\[\Span \colon \AdTrip \longrightarrow \Cat\]
sending each adequate triple to its \icat{} of spans, described
informally in the introduction. 
We will define $\Span$ more precisely as a functor into complete Segal \igpds{}. Recall that such functors can be obtained from the following general procedure:
\begin{construction}
Given a functor $C:D\rightarrow A$ between $\infty$-categories, we obtain a `singular complex' functor $\S_C \colon A \rightarrow \mathcal P(D)$ with $\mathcal P(D) = \Fun(D^\op,\Gpd)$ by currying the composition
\[D^\op \times A \xrightarrow{C^\op\times \id} A^\op\times A \xrightarrow{\Map_{A}} \Gpd,\]
i.e.\ $\S_C(X)_d = \Map_A(C(d),X)$. If $A$ is cocomplete, then $\S_C$ is right adjoint to the colimit extension $|{-}|_C:\mathcal P(D) \rightarrow A$, see \cite[Proposition 5.2.6.3]{HTT}.

For the cosimplicial object $\Del \rightarrow \Cat, n \mapsto [n]$ we obtain in this fashion an adjunction 
\[\begin{tikzcd}\asscat \colon \ssS\arrow[r, yshift=0.5ex] & \Cat \cocolon \nerve.\arrow[l, yshift=-0.5ex]\end{tikzcd},\]
where $\ssS = \mathcal P(D)$ is the $\infty$-category of simplicial $\infty$-groupoids. By results of Joyal, Lurie, Rezk and Tierney \cite{JT,LurieGoo,Rezk},
the \emph{nerve} functor $\nerve$ is fully faithful with essential image the full
subcategory of complete Segal objects inside $\ssS$, i.e. those
simplicial \igpds{} $T$ that satisfy the Segal condition and for which
the diagram
\[\begin{tikzcd}[column sep=3pc]
T_0 \arrow[rr,"\Delta"] \arrow[d,"{s}"{left}] && T_0 \times T_0  \arrow[d, "{(s,s)}"] \\
T_3 \arrow[rr, "{(d_{\{0,2\}},d_{\{1,3\}})}"] && T_1 \times T_1
\end{tikzcd}\]
is cartesian; see \cite[Section 10]{Rezk2} for this characterisation of completeness. A general cosimplicial object $C$ in $A$ therefore gives
rise to a functor $A \rightarrow \Cat$ via the composition of $\S_C$
with $\asscat$. When $\S_C$ takes values in (complete) Segal \igpds{} we
have a good understanding of this functor, since there is then
  no need to localise; this happens precisely when the cosimplicial
  object $C$ satisfies the dual version of the Segal and completeness
  conditions.
\end{construction}

\begin{example}
Consider the functor $\S_{B}\colon \Cat\rt \sGpd$ associated to the cosimplicial object
\[B\colon \Del\rt \Cat, \quad [n]\longmapsto[n]\join [n]^\op.\] Then $\S_B(\CC) \simeq \nerve\TwR(\CC)$: This follows from the observation that $\TwR$ arises from a right Quillen functor between the Joyal model structures whose left adjoint sends $[n]$ to $[n]\star [n]^\op$ (see \cite[Section 5.2.1]{HA} or \cite[Proposition 4.13]{HNP}). In particular $\S_B$ takes values in complete Segal \igpds{}.
\end{example}

Applying this construction to the functor $\TwR\colon \Del\rt \AdTrip$
sending each simplex $[n]$ to $\TwR([n])$ with the structure from
Example \ref{ex:twar simplex}, we obtain: 

\begin{definition}\label{def:span}
We define $\Span \colon \AdTrip \rightarrow \Cat$ as the composition 
\[\AdTrip\xrightarrow{\S_{\TwR}} \ssS \xrightarrow{\asscat} \Cat.\]
\end{definition}

\begin{theorem}[Barwick]\label{thm:Spanissegal}
The essential image of $\S_{\TwR}$ is contained in the complete Segal \igpds{}. In other words, there is a natural equivalence
$$
\Map_{\Cat}\big([n], \Span(X)\big)\simeq \Map_{\AdTrip}(\TwR([n]), X).
$$ 
\end{theorem}

We will give a new proof of this result (avoiding Barwick's recourse to a point-set model). Recall also that the category associated to a complete Segal \igpd{} $T$ can be described rather explicitly: For example, the core of $\asscat(T)$ is simply given by $T_0$, and its morphism \igpds{} are given by the fibres of $(d_1,d_0) \colon T_1 \rightarrow T_0^2$. In the case at hand this gives us the desired description of
$\Span(X, X_{\mm{in}}, X^{\mm{eg}})$: objects are objects of $X$,
morphisms from $x$ to $y$ are spans
$$\begin{tikzcd}[column sep=2pc]
x & z\arrow[l, twoheadrightarrow]\arrow[r, rightarrowtail] & y
\end{tikzcd}$$
with the left arrow egressive and the right arrow ingressive, and composition proceeds by pullback of spans. Completeness furthermore implies that $\core\Span(X,X_{\ing},X^\eg) \simeq \core(X)$ via the degenerate spans consisting of identities. Before we dive into the proof, let us record two immediate properties of the construction, the second of which in particular yields inclusions of subcategories 
$$\begin{tikzcd}
\big(X^{\eg}\big)^{\op}\arrow[r] & \Span(X) & X_{\ing}\arrow[l],
\end{tikzcd}$$
when applied to 
\[(X,\iota X,X^\eg) \longrightarrow X \longleftarrow (X,X_\ing,\iota X).\]

\begin{lemma}\label{span_op}
For an adequate triple $(X,X_\mathrm{in}, X^\mathrm{eg})$ with reverse $X^\mathrm{rev}=(X,X^\mathrm{eg},X_\mathrm{in})$ as in Example \ref{ex:adequate triple}, there is a natural canonical equivalence 
\[\Span(X)^\op \simeq \Span(X^\mathrm{rev}).\]
\end{lemma}

\begin{proof}
For each $[n]$, consider $\TwR([n])$ and $\TwR([n]^{\op})$ with the structure of an adequate triple as in Example \ref{ex:twar simplex}. There is a natural equivalence of cosimplical adequate triples $\TwR([-]^{\op})\simeq \TwR([-])^{\rev}$, sending an object $(i\geq j)$ in $\TwR([n]^{\op})$ to the object $(j\leq i)$ in $\TwR([n])$. This induces a natural equivalence of simplicial objects
\[
\Map_{\AdTrip}\big(\TwR([-]^{\op}),X\big)\simeq \Map_{\AdTrip}\big(\TwR([-])^{\rev},X\big)\simeq \Map_{\AdTrip}\big(\TwR([-]),X^{\rev}\big)
\]
so that the simplicial $\infty$-groupoid defining $\Span(X)$ is the reverse of that defining $\Span(X^\rev)$. The result follows since generally $\asscat(T)^\op \simeq \asscat(T^\rev)$ as this is true on simplices.
\end{proof}

\begin{proposition}\label{prop:spanswithcore}
Let $A \subseteq B$ be a  subcategory. Then the triples $(B,A,\iota B)$ and $(B,\iota B,A)$ are adequate and
\[\Span(B,A,\iota B) \simeq A \quad \text{and} \quad \Span(B,\iota B,A) \simeq A^\op.\]
\end{proposition}

\begin{proof}
We shall prove the first claim; the second then follows from \ref{span_op}. The triple $(B,A,\iota B)$ is evidently adequate (cf.\ Example \ref{ex:adequate triple}). Note that the \igpd{} of functors from an adequate triple $(X,X_\ing,X^\eg)$ into $(B,A,\iota B)$ is equivalent to the \igpd{} of functors $X\rightarrow B$ which invert the edges in $X^\eg$ and take those of $X_\ing$ into $A$. Applying this to $X=\TwR([n])$ shows that the remaining claim is equivalent to the assertion that the source map
\[s\colon \TwR([n]) \longrightarrow [n]\]
is a localisation (necessarily at those maps whose source component is an equivalence). This is in fact true for all $\infty$-categories $C$ in place of $[n]$, for example since $\TwR(\C) \rightarrow \C$ is a cocartesian fibration with contractible fibres (which are always localisations, see \cref{lem:cartlocalise} below). The localisation property is, however, particularly easy to see when $C$ admits a terminal object: In that case the source evaluation $\TwR(\C) \rightarrow \C$ is even a Bousfield localisation, with fully faithful left adjoint sending $c \mapsto (c \to \ast)$.
\end{proof}

We now turn to the proof of Theorem \ref{thm:Spanissegal}. Instead of following Barwick's strategy of explicitly filling simplices in a point-set implementation of the above construction, our proof will be a simple adaptation of the argument given in \cite[Section 2.1]{CDH2} for the case of stable \icats{} (with all maps ingressive and egressive). In fact, we will prove a slightly stronger statement, which will ultimately allow us to deduce \Cref{thma}. It uses the following definition:

\begin{definition}\label{def:QA}
Let $A$ be an \icat{} and $X$ an adequate triple. We write $\Q_A(X)$ for the \icat{} $\Fun_{\AdTrip}(\TwR(A), X)$, where $\TwR(A)$ is the adequate triple from Example \ref{ex:twar simplex} and denote by $\Q\colon \Cat^{\op}\times \AdTrip\rt \AdTrip$ the resulting functor.
\end{definition}

In particular, each adequate triple $X$ gives rise to a natural simplicial diagram $\Q_{\bullet}(X)$ in $\Cat$.
\begin{lemma}\label{Qsegal}
Let $X$ be an adequate triple. Then the simplicial \icat{} $\Q_\bullet(X)$ satisfies the Segal and completeness conditions, that is the Segal maps
\[\Q_n(X) \longrightarrow \Q_1(X) \times_{\Q_0(X)} \Q_1(X) \dots \times_{\Q_0(X)} \Q_1(X)\]
are equivalences and
\[\begin{tikzcd}[column sep=3pc]
\Q_0(X) \arrow[rr,"{\Delta}"] \arrow[d,"s"{left}] && \Q_0(X)\times \Q_0(X) \arrow[d, "{(s,s)}"] \\
\Q_3(X) \arrow[rr, "{(d_{\{0,2\}},d_{\{1,3\}})}"] && \Q_1(X) \times \Q_1(X)
\end{tikzcd}\]
is cartesian.
\end{lemma}
Note that this immediately implies Theorem \ref{thm:Spanissegal}, as $\core \colon \Cat \rightarrow \Gpd$ preserves limits.

\begin{proof}
Let $\mathcal J_n \subseteq \TwR([n])$ denote the subposet consisting
of those $(i \leq j)$ with $j \leq i+1$, i.e. the zig-zag along the
bottom. Note that $\mathcal J_n$ decomposes as an iterated pushout
\[\mathcal J_n \simeq \TwR([1]) \cup_{\TwR([0])} \TwR([1])\cup \dots \cup_{\TwR([0])} \TwR([1])\]
along the Segal maps; in fact, the nerve $\nerve(\mathcal J_n)$ is already the iterated pushout of the nerves $\nerve(\TwR([1]))$ in simplicial \igpds{} (see \cite[Proposition 5.13]{HaugsengSpans} for a similar argument). The iterated pullback appearing in the Segal condition is therefore equivalent to the full subcategory $J_n(X)$ of $\Fun(\mathcal J_n,X)$ spanned by those functors taking left-pointing edges in $\mathcal J_n$, i.e. those of the form $(i \leq i+1) \rightarrow (i \leq i)$, to ingressives and right pointing arrows, namely $(i \leq i+1) \rightarrow (i+1 \leq i+1)$, to egressives. Furthermore, this translates the Segal map to the map $\Q_n(X) \rightarrow J_n(X)$ induced by the restriction
\[\Fun(\TwR([n]),X) \longrightarrow \Fun(\mathcal J_n,X).\]
But from the pointwise formula for Kan extensions, one finds that a diagram $F \colon \TwR([n]) \rightarrow X$ lies in $\Q_n(X)$ if and only if it is right Kan extended from its restriction to $\mathcal J_n$, which has to lie in $J_n$. The claim now follows from \cite[Proposition 4.3.2.15]{HTT}, since right Kan extension from a full subcategory is fully faithful.

Similarly, for completeness we first note that the map $P \rightarrow \Q_3(X)$ from the pullback in question to the lower left corner is fully faithful since the degeneracy $\Q_0(X)^2 \rightarrow \Q_1(X)^2$ is (as $|\TwR([1])| \simeq *$). We claim that its essential image consists exactly of those diagrams whose edges are all equivalences; since also $|\TwR([3])| \simeq *$ these are precisely the constant ones, which gives the result.
So consider a diagram
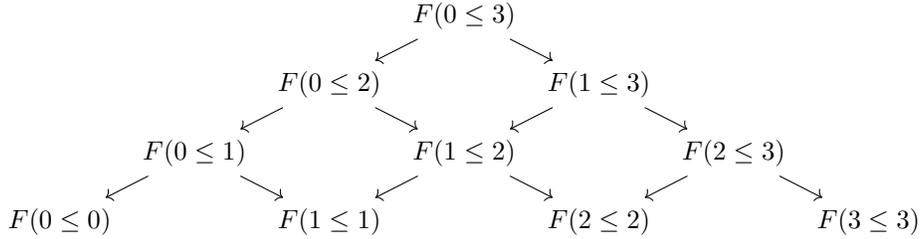
\[
\begin{tikzcd}
[row sep=2ex, column sep=0.8ex]
 & & & F(0 \leq 3) \ar[tail,rd]\ar[two heads,ld]\\
 & & F(0 \leq 2) \ar[two heads,ld]\ar[tail,rd]& &F(1 \leq 3)\ar[tail,rd]\ar[two heads,ld]  \\
 & F(0 \leq 1) \ar[two heads,ld]\ar[tail,rd]& & F(1 \leq 2) \ar[two heads,ld]\ar[tail,rd]& & F(2 \leq 3)\ar[tail,rd]\ar[two heads,ld] \\
F(0 \leq 0) && F(1 \leq 1) && F(2 \leq 2) && F(3 \leq 3)
\end{tikzcd}
\]
all of whose squares are (ambigressive) cartesian and such that the four compositions
\[F(0 \leq 2) \longrightarrow F(0 \leq 0), \quad F(0 \leq 2) \longrightarrow F(2 \leq 2) \]
\[F(1 \leq 3) \longrightarrow F(1 \leq 1), \quad F(1 \leq 3) \longrightarrow F(3 \leq 3) \]
are equivalences. Then it first follows that, as pullbacks of equivalences, also $F(0 \leq 3) \rightarrow F(0 \leq 1)$ and $F(0 \leq 3) \rightarrow F(2 \leq 3)$ are equivalences and then by two-out-of-six the entire outer slopes are. But then by commutativity of the larger rectangles also $F(0 \leq 1) \rightarrow F(1 \leq 1)$ and $F(2 \leq 3) \rightarrow F(2 \leq 2)$ are equivalences, and then as pullbacks thereof also $F(0 \leq 2) \rightarrow F(1 \leq 2)$ and $F(1 \leq 3) \rightarrow F(1 \leq 2)$. Finally, this implies that also $F(1 \leq 2) \rightarrow F(1 \leq 1)$ and $F(1 \leq 2) \rightarrow F(2 \leq 2)$ are equivalences by two-out-of-three.
\end{proof}

Next we shall describe functors into span \icats{}. More precisely, we will show that the equivalence between functors $[n]\rt \Span(X)$ and maps of adequate triples $\TwR([n])\rt X$ from Theorem \ref{thm:Spanissegal} extends to all \icats{}:

\begin{theorem}\label{modelfuncob}
The functors $\begin{tikzcd}\TwR \colon \Cat \arrow[r, yshift=0.5ex] & \arrow[l, yshift=-0.5ex] \AdTrip \cocolon \Span \end{tikzcd}$ form an adjoint pair. In other words, for each \icat{} $A$ and each adequate triple $X$, there is a natural equivalence
\[\Map_{\AdTrip}(\TwR(A),X) \simeq \Map_{\Cat}(A,\Span(X)),\]
where $\TwR(A)$ is as in Example \ref{ex:twar simplex}.
\end{theorem}
\begin{remark}
Using a different argument this result was first sketched by Raskin in \cite[Chapter 20]{Raskin}. The resulting description of functors into span categories is for example also used in \cite{HS} to perform parametrised surgery on the cobordism \icats{} $\mathrm{Cob}(\C,\Qoppa)$ from \cite{CDH2}, which are hermitian refinements of $\Span(\C)$.
\end{remark}
Recall from Definition \ref{def:QA} that we abbreviate
\[\Q_A(X) = \Fun_{\AdTrip}(\TwR(A), X).\]
The theorem will follow readily from the following:

\begin{proposition}\label{lem:Q_preserves_limits}
For each adequate triple $X$, the functor $\Q_{-}(X) \colon \Cat^\op \rightarrow \Cat$ preserves limits. 
\end{proposition}

We will employ the (opposite of the) following criterion:

\begin{lemma}\label{lem:colimoutofcat}
Let $f \colon \D \rightarrow \E$ be a functor from a small to a
cocomplete category, such that the right adjoint in the induced
adjunction
\[\begin{tikzcd}{|{-}|_f} \colon \mathcal P(\D) \arrow[r, yshift=0.5ex] & \arrow[l, yshift=-0.5ex] \E \cocolon {\S_f} \end{tikzcd}\]
is fully faithful. Then for a functor $F \colon \E \rightarrow \C$ to another cocomplete category the following are equivalent: 
\begin{enumerate}
\item $F$ preserves colimits,
\item the natural map $|{-}|_{F\circ f} \Longrightarrow F \circ |{-}|_f$ is an equivalence, and
\item\label{colimoutofcat_3} \begin{enumerate}
\item[(i)] the natural map $f_!(F\circ f) \Longrightarrow F$ is an equivalence, i.e. $F$ is left Kan extended from its restriction along $f$, and
\item[(ii)] $|{-}|_{F\circ f} \colon \mathcal P(\D) \rightarrow \C$ inverts all maps that $|{-}|_f$ inverts.
\end{enumerate}
\end{enumerate}
\end{lemma}

Applying this to $\Del \subset \Cat$ and taking opposites in particular shows that a functor $F \colon \Cat^\op \rightarrow \C$ preserves limits if and only if it is right Kan extended from its restriction along $\Del^\op \subset \Cat^\op$ and the restriction $F \colon \Del^\op \rightarrow \C$ is a complete Segal object in $\C$: The maps
\[\Delta^{\{0,1\}} \cup_{\Delta^{\{1\}}} \dots \cup_{\Delta^{\{n-1\}}} \Delta^{\{n-1,n\}} \longrightarrow \Delta^n \quad \text{and} \quad \Delta^3/\Delta^{\{0,2\}}, \Delta^{\{1,3\}} \longrightarrow \Delta^0\] 
encoding the Segal and completeness conditions are (by design) categorical equivalences, i.e.\ inverted by $\asscat\colon \ssS\rightarrow \Cat$, proving that if $|{-}|_{F^\op} \colon \ssS \rightarrow C^\op$ inverts these then $F \colon \Del^\op \rightarrow C$ is a complete Segal object, and the Yoneda lemma reduces the converse to the case $C = \Gpd$, where it follows from the fact that $|{-}|_{(\nerve C)^\op} \simeq \Hom_{\Cat}(\asscat({-}),C)$ as (colimit preserving) functors $\ssS \rightarrow \Gpd^\op$.

\begin{proof}
Let us point out that any colimiting cocone $G\colon I^{\rhd}\to E$ arises as the image of a colimiting cocone $G'\colon I^{\rhd}\to \mathcal{P}(D)$ under $|{-}|_f$: indeed, one can take $G'$ to be the colimit of the diagram $S_f\circ G_{|I}\colon I\rt \mathcal{P}(D)$. Consequently, a functor out of $E$ preserves colimits if and only if its composition with $|{-}|_f$ does.
Using this, (1) $\Leftrightarrow$ (2) is an immediate consequence of \cite[Lemma 5.1.5.5]{HTT}.

For (2) $\Rightarrow$ (3), (ii) follows immediately. To see (i), note that $F$ agrees with the left Kan extension of $|{-}|_{F\circ f} = F \circ |{-}|_f$ along $|{-}|_f$, since $|{-}|_f$ is a localisation \cite[Proposition 5.2.7.12]{HTT}. In turn, $|{-}|_{F\circ f}$ is the left Kan extension of $F\circ f$ along the Yoneda embedding $h\colon D\rt \mathcal{P}(D)$ \cite[Lemma 5.1.5.5]{HTT}. The claim then follows from transitivity of Kan extensions.

Finally, for (3) $\Rightarrow$ (2) note that (ii) implies that
$|{-}|_{F\circ f}$ descends along $|{-}|_f$ to some functor $G \colon \E
\rightarrow \C$, which is then automatically its left Kan extension
along $|{-}|_f$. Part (i) and transitivity of Kan extensions then imply
that $G \simeq F$, or in other words that $|{-}|_{F\circ f} \simeq F \circ
|{-}|_f$ as desired.
\end{proof}

\begin{proof}[Proof of \cref{lem:Q_preserves_limits}]
Since $\Q(X)$ is a complete Segal object in $\Cat$ by \cref{Qsegal},
we have already shown Item (ii) of Lemma
\ref{lem:colimoutofcat}\ref{colimoutofcat_3}.
So we need only show that $\Q(X) \colon \Cat^\op \rightarrow \Cat$ is right Kan extended from its restriction to $\Del^\op$. By the pointwise formula for Kan extensions this means that the tautological map 
\[\Q_{\DD}(X) \longrightarrow \lim_{n \in (\Del/\DD)^\op} \Q_n(X) \]
is an equivalence for all $\DD \in \Cat$. The diagram $\Del/\DD \rightarrow \Del \hookrightarrow \Cat$, giving rise to the right hand side has colimit $\D$ and is a typical example of a functor $F \colon I \rightarrow \Cat$ such that the natural map 
\begin{equation}\label{goodcolim}\tag{$\ast$}
\colim_{i \in I}\nerve(F_i) \longrightarrow \nerve(\colim_{i \in I}F_i)
\end{equation}
is an equivalence (with both sides evaluating to $\nerve \DD$). We will directly show that $\Q_{(-)}(X)$ preserves all limits over diagrams with this property.

To this end note that there are natural equivalences of simplicial \igpds{}
\[\nerve(\TwR(\CC)) \simeq \big[i \longmapsto \Map_{\Cat}([i] \star [i]^\op,\CC)\big] \simeq \big[i \longmapsto \Map_{\ssS}(\Delta^i \star (\Delta^i)^\op,\nerve(\CC))\big].\]
In particular, since $\Delta^i \star (\Delta^i)^\op \simeq \Delta^{2i+1}$ is completely compact in $\ssS$, we find an equivalence
\[\begin{tikzcd}\colim_{i \in I}\TwR(F_i) \arrow[r, "\sim"] & \TwR\big(\colim_{i \in I}F_i\big)\end{tikzcd}\]
because of \eqref{goodcolim}; let us warn the reader that the above map need not be an equivalence without any assumption on $F$, as the diagram $[1] \xleftarrow{0} [0] \xrightarrow{1} [1]$ with pushout $[2]$ shows. In the case at hand it follows that the map
\[\Fun(\TwR(\colim_{i \in I} F_i), X) \longrightarrow \lim_{i \in I}\Fun(\TwR(F_i),X)\]
is an equivalence.

We have to show that this restricts to an equivalence between the full subcategory $\Q_{\colim_i F_i}(X)$ on the left (spanned by morphisms of adequate triples) and $\lim_{i \in I}\Q_{F_i}(X)$ on the right. To this end, notice that any ambigressive pullback square in $\TwR(A)$, say 
\[
\begin{tikzcd}
[row sep=2ex, column sep=0.8ex]
  & (a_0 \rightarrow a_3) \ar[tail,rd]\ar[two heads,ld]\\
  (a_0 \rightarrow a_2) \ar[tail,rd]& &(a_1 \rightarrow a_3)\ar[two heads,ld]  \\
 & (a_1 \rightarrow a_2),
\end{tikzcd}
\]
 is contained (up to equivalence) in $\TwR([3])$ for a map $[3]\rt A$, namely that given by the string $a_0 \rightarrow a_1 \rightarrow a_2 \rightarrow a_3$. Consequently, a functor $\TwR(A) \rightarrow X$ is a morphism of adequate triples if and only if each restriction
\[\TwR([3]) \longrightarrow \TwR(A) \longrightarrow X\]
is such. To apply this, note that one has
\[\Map_{\Cat}\big([3],\colim_{i \in I}F_i\big) \simeq \nerve_3(\colim_{i \in I} F_i) \simeq \colim_{i \in I}\nerve_3(F_i) \simeq \colim_{i \in I}\Map_{\Cat}([3],F_i)\]
whenever $F$ satisfies \eqref{goodcolim}, so that a functor $\TwR(\colim_i F_i) \rightarrow X$ lies in $\Q_{\colim_i F_i}(X)$ if and only if for every $j \in I$ and map $[3] \rightarrow F_j$, the composite
\[\TwR([3]) \longrightarrow \TwR(F_j) \longrightarrow \TwR(\colim_{i } F_i) \longrightarrow X\]
lies in $\Q_3(X)$. But this is the case if and only if the functor restricts to an element of $\Q_{F_j}(X)$ for every $j \in I$ and thus by definition if and only if it defines an element of $\lim_{i}\Q_{F_i}(X)$, as desired.
\end{proof}

\begin{proof}[Proof of Theorem \ref{modelfuncob}]
By Proposition \ref{lem:Q_preserves_limits} $\iota\Q(X) \simeq \Map_{\AdTrip}(\TwR(-),X) \colon \Cat^\op \rightarrow \Gpd$ preserves limits as does $\Map_{\Cat}(-,\Span(X))$, so by Lemma \ref{lem:colimoutofcat} they are both right Kan extended from $\Del^\op \subset \Cat^\op$.  Thus they agree if their restrictions to $\Del^\op$ agree, which is the case by Theorem \ref{thm:Spanissegal}.
\end{proof}

In fact, one can upgrade Theorem \ref{modelfuncob} to describe the full functor \icat{} $\Fun(K,\Span(X))$:

\begin{corollary}\label{FunofSpan}
There is a canonical equivalence
	\[\Fun(A,\Span(X)) \simeq \Span(\Q_A(X)),\]
	natural in $A \in \Cat$ and $X \in \AdTrip$. 
\end{corollary}

We warn the reader that this result is misstated as $\Fun(A,\Span(X)) \simeq \Q_A(X)$ in \cite[Section 20.9]{Raskin}, which is only true on groupoid cores.

\begin{proof}
The left adjoint $\TwR\colon \Cat\rightarrow \AdTrip$ preserves the cartesian product. By adjunction, the right adjoint $\Span \colon \AdTrip \rightarrow \Cat$ then sends the internal mapping object $\Q_A(X)=\Fun_{\AdTrip}(\TwR(A), X)$ from \Cref{lem:adtrip_internal_hom} to the internal mapping object $\Fun(A , \Span(X))$ in $\Cat$.
\end{proof}

\begin{example}
Since the functors $s \colon \TwR(A) \rightarrow A$ and $t \colon \TwR(A) \rightarrow A^\op$ are localisations, it follows easily from the formula of \cref{FunofSpan} that the functors
\[\Fun(A,X^\ing) \longrightarrow \Fun(A,\Span(X)) \quad
 \text{and} \quad \Fun(A^\op,X^\eg)^\op \longrightarrow
 \Fun(A,\Span(X))\]
are inclusions of (non-full) subcategories for all $A \in \Cat$. Alternatively, this statement can also be reduced to the case $A=[0]$ (which we recorded before \cref{span_op}) by means of \cref{inclusionfull}.

In a different direction note that even if $X$ admits pullbacks and carries the trivial structure of an adequate triple with all maps ingressive and egressive, the same is typically not true for $\Q_A(X) = \Fun_{\AdTrip}(\TwR(A),X)$; this phenomenon is the basis for the notions of forwards and backwards surgery in \cite{HS}. 
\end{example}

\section{Cartesian fibrations and \icats{} of spans}\label{sec:functors}
The goal of this section is to describe the cocartesian fibration associated to the functor $\Span \colon \AdTrip \rightarrow \Cat$, which is itself given by a certain span \icat{} of adequate triples. As a direct application we will obtain a new and fairly direct proof of the main result of \cite{BGN}, which identifies the dual cartesian fibration of a cocartesian fibration in terms of a span construction. We shall more generally analyse the interaction of span categories and (co)cartesian fibrations in Section \ref{sec:fibs} below. 

To understand (co)cartesian edges in span \icats{}, we start by reproducing a criterion due to Barwick. In fact, we provide a proof different from \cite{Barwick} and inspired by the proof of the additivity theorem for Grothendieck--Witt theory in \cite[Section 2.6]{CDH2}, see Example \ref{remarksaboutbarwickI} for the details of this connection.

\begin{theorem}[Barwick]\label{barwick_cocart}
Let $p\colon Y\rightarrow X$ be a morphism in $\AdTrip$ and let $f$ be an edge
in $X_\ing$ such that the following conditions hold:
\begin{enumerate}
    \item Every pullback of $f$ along an egressive edge has a lift in
      $Y_\ing$ which is simultaneously $p$-cocartesian and
      $p_\ing$-cocartesian, for any choice of lift of its
      source.
    \item Consider any commutative square $\sigma$ in $Y$ 
    \begin{center}
       \begin{tikzcd}
w \arrow[tail, r, "g'"] \arrow[two heads, d, "g"'] & v \arrow[d] \\
u \arrow[r, rightarrowtail, "\psi"]                                 & y                   
    \end{tikzcd}
    \end{center}

such that $p(\sigma)$ is an ambigressive pullback in $X$, the edge $g'$ is ingressive, the morphism $g$ is egressive, and the morphism $\psi$ is a $p$-cocartesian and ingressive lift of $f$. Then $g'$ is $p$-cocartesian if and only if the square is an ambigressive pullback.
\end{enumerate}
Then an edge $x\rightarrow y$ of $\Span(Y)$, represented by a span
\begin{center}
    \begin{tikzcd}[row sep=1pc, column sep=1.7pc]
  & u \arrow[ld, twoheadrightarrow, "\phi"'] \arrow[rd, rightarrowtail, "\psi"] &   \\
x &                                                           & y
\end{tikzcd}
\end{center} such that $\psi$ covers $f$, is $\Span(p)$-cocartesian if $\phi$ is $p^\eg$-cartesian and $\psi$ is $p$-cocartesian as indicated.
\end{theorem}

\begin{proof}
Unwinding definitions, we have to show that for any span $\begin{tikzcd}x &\ar[two heads,l] w \ar[tail,r] & z\end{tikzcd}$ in $Y$ the solid diagram (ignoring the numbers for a moment)
\begin{center}
\begin{tikzcd}
   &                                           & w \arrow[two heads, lldd, bend right] \arrow[tail, rrdd, bend left] \arrow[ld, "(1)", dashed] \arrow[rd, "(2)"', dashed] &                                                                          &    \\
   & u \arrow[ld, twoheadrightarrow, "\phi"] \arrow[rd, rightarrowtail, "\psi"'] &                                                                                                                                                                       & \bullet \arrow[ld, "(3)",  dashed] \arrow[rd, "(4)"',  dashed] &    \\
x&                                           & y                                                                                                                                                                    &                                                                          & z
\end{tikzcd}
\end{center} 
admits an essentially unique dashed filling lying over a given entirely solid diagram in $\Span(X)$, such that all left pointing arrows are egressive and all right pointing arrows ingressive, and the top square is cartesian. We then first observe that the second condition on $f$ (the image of $\psi$ in $X$) implies that for such square in which the arrows $w \rightarrow u$ and $w \rightarrow \bullet$ are egressive and ingressive, respectively the assertion that the square is cartesian and $\bullet \rightarrow y$ is egressive is equivalent to the assertion that the map $w \rightarrow \bullet$ is $p$-cocartesian, so we may instead show that there is a unique filler with this property.

We do so by filling the diagram step by step, as indicated by the numbers in the above diagram, essentially uniquely each time: 

\begin{enumerate}
    \item There exists a unique egressive filler because $\phi$ is $p^\eg$-cartesian.
    \item The first assumption on $f$ provides a $p$-cocartesian left, that is automatically unique and also $p_\ing$-cocartesian and ingressive.
    \item There exists a unique filler making the middle square commute because $w \rightarrow \bullet$ is $p$-cocartesian.
    \item There is a unique ingressive filler because $w\rightarrow \bullet$ is $p_\ing$-cocartesian.
\end{enumerate}
\end{proof}
In the special case where $f$ is an equivalence, the two assumptions are automatic, so we obtain:

\begin{corollary}\label{barwick_cocart_equiv}
If $p\colon Y\rightarrow X$ is a morphism in $\AdTrip$ then a span in $Y$ of the form
\begin{center}
    \begin{tikzcd}[row sep=1pc, column sep=1.7pc]
  & u \arrow[ld, twoheadrightarrow, "\phi"'] \arrow[rd, "\sim"] &   \\
x &                                                           & y,
\end{tikzcd}
\end{center}
represents a $\Span(p)$-cocartesian morphism whenever $\phi$ is $p^\eg$-cartesian.
\end{corollary}

\begin{remark}\label{remarksaboutbarwickI}
We feel obliged to point out two oversights in the statement of the above result in \cite[12.2 Theorem]{Barwick}:
\begin{enumerate}[(i)]
\item Barwick requires that the edge $\phi$ be $p$-cartesian and not $p^\eg$-cartesian. Our proof above hopefully makes it transparent why this is not enough. For an explicit counterexample consider $X = Y = [1]^2$ with $p$ the identity, where we equip the source with the triple structure with the horizontal maps (and the identities) ingressive and everything but the horizontal maps egressive and the target with the same ingressives but all maps egressive.
\item As he is working at the point-set level, Barwick has to show that $\Span(p)$ is an inner fibration. To this end he assumes that $p$ is an inner fibration, but in fact requires the stronger assumption that $p$ is an isofibration, as can be seen by lifting a $2$-simplex in $\Span(Y)$ of the form
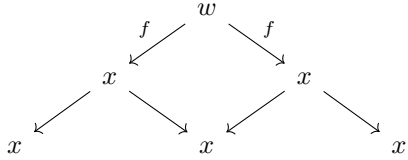
\begin{center}
\begin{tikzcd}[row sep=1.2pc, column sep=1.8pc]
   &                                           & w  \arrow[ld, "f"'] \arrow[rd, "f"] &                                                                          &    \\
   & x \arrow[ld] \arrow[rd] &                                                                                                                                                                       & x \arrow[ld] \arrow[rd] &    \\
x&                                           & x                                                                                                                                                                    &                                                                          & x
\end{tikzcd}
\end{center} 
where $f$ is some equivalence and all other maps identities.

In order to carry out iterated span constructions as in \cite{BGN} at the point-set level one therefore needs to check that if $p$ is an isofibration, then so is $\Span(p)$. This is indeed the case as a direct consequence of equivalences in span \icats{} being precisely the spans of equivalences (which itself follows from the completeness assertion for $\core\Q(X)$ in Lemma \ref{Qsegal}, see also \cite[Proposition 3.4]{Barwick}).\end{enumerate}
\end{remark}

\begin{example}

In particular in light of the previous remark, let us recount the original purpose of Theorem \ref{barwick_cocart} in \cite{Barwick}, which is Barwick's unfurling construction. The input is an adequate triple $(X, X_\ing, X^\eg)$ and a functor $F \colon X^\op \rightarrow \Cat$ such that 
\begin{enumerate}
\item for every ingressive $f \colon x \rightarrow y$ the induced functor $f^* \colon F(y) \rightarrow F(x)$ admits a left adjoint $f_!$, and 
\item the image under $F$ of every ambigressive pullback in $X$ is left-adjointable, i.e.\ given the left square in 
\[\begin{tikzcd}
x \ar[r, tail, "f"] \ar[d, two heads, "g"] & y \ar[d,two heads, "h"] & & F(x) \ar[r,"f_!"]& F(y)  \\
w \ar[r, tail, "i"] & z & & F(w) \ar[u,"g^*"] \ar[r,"i_!"] & F(z) \ar[u,"h^*"']
\end{tikzcd}            \]
the right one commutes via the Beck-Chevalley transformation $f_!g^* \Rightarrow h^*i_!$.
\end{enumerate}
The output is an extension of $F \colon (X^\eg)^\op \rightarrow \Cat$ to a functor $\mathrm{Unf}(F) \colon \Span(X,X_\ing,X^\eg) \rightarrow \Cat$, the \emph{unfurling} of $F$, whose additional functoriality is given by the left adjoints to the images of the ingressive edges in $X$. 

It is constructed by considering the cartesian unstraightening $p \colon Y\rightarrow X$ of $F$, which carries the structure of an adequate triple with ingressives $p^{-1}(X_\ing)$ and egressives $Y^\dagger$, the subcategory formed by the $p$-cartesian lifts of degres in $X^\eg$, by \Cref{prop:adequate triple}. The claim is that the induced map
\[\Span(p) \colon \Span(Y, p^{-1}(X_\ing),p^\dagger) \longrightarrow \Span(X,X_\ing,X^\eg)\]
is a cocartesian fibration and that its cocartesian straightening is the desired extension \[\mathrm{Unf}(F) \colon \Span(X,X_\ing,X^\eg) \rightarrow \Cat.\]

The first claim follows from an application of Theorem \ref{barwick_cocart}: Given a span 
\[\begin{tikzcd} x & \ar["g"',l, two heads] y \ar["f",r,tail] & z\end{tikzcd}\] 
in $X$ and an element $c \in F(x) \subseteq Y$, we can build a $\Span(p)$-cocartesian lift as follows. Per construction we can first choose a $p$-cartesian edge $\phi \colon d \rightarrow c$ covering $g$. By the first assumption on $F$ and the fibrational characterisation of adjunctions, the edge $f$ then further admits a locally $p$-cocartesian lift $\psi \colon d \rightarrow e$, which is in fact $p$-cocartesian by \cite[Corollary 5.2.2.4]{HTT}. Note that this span in $Y$ indeed defines a morphism $c \rightarrow e$ in $\Span(Y, p^{-1}(X_\ing),p^\dagger)$. To see that it is $\Span(p)$-cocartesian we check the conditions of \Cref{barwick_cocart}. For the first condition, note that all egressive pullbacks of $f$ admit $p$-cocartesian lifts by the same argument that gives the existence of $\phi$, since such pullbacks are again ingressive. That these lifts are also $p_\ing$-cocartesian follows directly from the definition of the ingressives in $Y$. Next, a square in $Y$ as in the second condition automatically has the form of the solid square
\[\begin{tikzcd}
g^*c \ar[dd,"\mathrm{cart}", two heads] \ar[tail, r,"\mathrm{cocart}"] & f_!g^* c \ar[rd, dotted, "\sim"'] \ar[r, tail] & d \ar[d] \\
& & h^*i_!c \ar[d,"\mathrm{cart}", two heads] \\
c  \ar[rr,"\mathrm{cocart}", tail]&  & i_!c 
\end{tikzcd}\]
if it covers the ambigressive pullback
\[\begin{tikzcd}
x \ar[r, tail, "f"] \ar[d, two heads, "g"] & y \ar[d,two heads, "h"]\\
w \ar[r, tail, "i"] & z 
\end{tikzcd}            \]
in $X$. The dotted map is the Beck-Chevalley transformation, which is an equivalence by the second assumption on $F$, and the upper right triangle commutes by a simple diagram chase. The task is to show that the large square is an ambigressive pullback if and only if the upper right horizontal map is an equivalence. But if the diagram is ambigressive the map $d \rightarrow i_!c$ is by definition $p$-cartesian, which indeed forces the upper right hand corner to consist of equivalences. And conversely, if the corner consists of equivalences then the right hand composition is cartesian, making the square ambigressive and also forcing it to be cartesian as it covers a cartesian base (see the proof of Proposition \ref{prop:adequate triple} above for the relevant calculation).

To apply (the corrected version of) Barwick's criterion we finally have to check that $\phi$ is $p^\eg$-cartesian. But since it is $p$-cartesian we are saved by the left-cancellation property of cartesian edges.

Finally, it remains to check that the composite
\[(X^\eg)^\op \simeq \Span(X,\core X, X^\eg) \longrightarrow \Span(X,X_\ing,X^\eg) \xrightarrow{\mathrm{Unf}(F)} \Cat\]
agrees with the restriction of $F$. Unwinding definitions this composite is given as the cocartesian straightening of $\Span(Y,p^{-1}(\core X),Y^\dagger) \rightarrow \Span(X,\core X, X^\eg)$ and the statement that this is indeed the restriction of $F$ is the main result of \cite{BGN}, which we shall reprove momentarily in Corollary \ref{BGNidentification} below.
\end{example}

\begin{example}
Suppose that $\C$ is stable and consider the target projection $p=t \colon \Ar(\C) \rightarrow \C$, regarded as a map of adequate triples with all maps both ingressive and egressive. One readily checks that every edge $f$ in $\C$ satisfies the assumptions of Theorem \ref{barwick_cocart}. It follows that $\Span(\Ar(\C)) \rightarrow \Span(\C)$ is a cocartesian fibration for every stable $\C$, and thus by Lemma \ref{span_op} also a cartesian fibration. Taking classifying \igpds{} (i.e.~classifying spaces) preserves pullback squares whose right edge is both a cocartesian and cartesian fibration, as can be shown from Quillen's Theorem B (see \cite{Steimle,Steinebrunner} or the proof of \cite[Theorem 2.5.1]{CDH2} for an account in the present language). Therefore we conclude that 
\[|\Span(\C)| \xrightarrow{c \mapsto (c \rightarrow 0)}  |\Span(\Ar(\C))| \xrightarrow{t} |\Span(\C)|\]
is a fibre sequence. Since $\mathcal K(\C) = \Omega|\Span(\C)|$ is one possible definition of the algebraic K-space of $\C$, looping this fibre sequence once recovers Waldhausen's additivity theorem
\[\mathcal K(\Ar(\C)) \simeq \mathcal K(\C) \times \mathcal K(\C)\] 
since $t$ is clearly split and $\mathcal K(\C)$ an $\mathbb E_\infty$-group using the direct sum in $C$. In this way, Proposition \ref{barwick_cocart} connects to the discussion of additivity in \cite[Section 2]{CDH2}.
\end{example}

As a further application of \Cref{barwick_cocart} we can identify the cocartesian fibration corresponding to the functor $\Span \colon \AdTrip \rightarrow \Cat$. More generally, let $X$ be an \icat{} and consider a diagram of adequate triples
\[F\colon X^{\op}\rt \AdTrip.\]
Writing $p\colon Y\rt X$ for the cartesian fibration classified by the underlying functor of $F$, we can use Proposition \ref{prop:cartadequateforspan} to endow $Y$ with the structure of an adequate triple:

\begin{notation}\label{not:fibration of adtrip}
By forgetting different pieces of the data, the functor $F\colon X^{\op}\rt \AdTrip$ gives three diagrams of \icats{} $F,F_{\ing}$ and $F^{\eg}$, which are classified by three cartesian fibrations 
\[p\colon Y= \Uncart(F) \rt X, \quad p_{\ing}\colon Y_{\ing} = \Uncart(F_\ing) \rt X \quad \text{and} \quad p^{\eg}\colon Y^{\eg} = \Uncart(F^\eg) \rt X,\] respectively, where we generally denote by
\[\Uncart \colon \Fun(X^\op,\Cat) \longrightarrow \Cart(X) \quad \text{and} \quad \Unco \colon \Fun(X,\Cat) \longrightarrow \Cocart(X)\]
the cartesian and cocartesian unstraightening functors respectively. In our specific situation we further write $Y_{\fb}, Y_\ing^{\fb},Y^\eg_{\fb}$ for the pullback of $Y,Y_{\ing},Y^{\eg}$ along $\iota X\rightarrow X$.
\Cref{prop:cartadequateforspan} then shows that $p\colon Y\rt X$ is part of a map of adequate triples
$$\begin{tikzcd}
Y^{\glob}=\big(Y, Y_{\ing}^{\fb}, Y^{\eg}\big)\arrow[r] & (X, \core X, X)
\end{tikzcd}$$
which gives rise to a map
 \begin{equation}\label{eq:span fibration}\Span(p)\colon \Span(Y^{\glob})\rightarrow \Span(X,\iota X, X)\simeq X^{\op}.
 \end{equation}
\end{notation}

We will show that this is the cocartesian fibration classifying $\Span \circ F$. First we note:

\begin{lemma}\label{lem:cocart_Span_unstr}
	The functor $\Span(p)$ is a cocartesian fibration.
\end{lemma}
\begin{proof}
	Corollary \ref{barwick_cocart_equiv} immediately implies that $\Span(p)$ is a cocartesian fibration, with cocartesian edges given by spans of the form:
	\begin{tikzcd}[cramped, column sep=1.2pc] y' & y\arrow[l,two heads]\arrow[r, equals] & y,\end{tikzcd} where $y\rightarrow y'$ is a $p^\eg$-cartesian edge in $Y$.	
\end{proof}

We will now describe the functor classified by this cocartesian fibration:

\begin{theorem}\label{thm:span_funct_cocart}
	Let $X$ be an \icat{}. Then there exists a natural equivalence  
	\[\alpha\colon \Unco(\Span\circ(-))\Longrightarrow \Span((\Uncart(-))^{\glob})\] of functors $\Fun(X^\op,\AdTrip)\rightarrow \Cocart(X^{\op})$. In other words, for each $F \colon X^\op \rt \AdTrip$, the induced cocartesian fibration \eqref{eq:span fibration} is classified by the functor $\Span\circ F\colon X^{\op}\rt \Cat$.
\end{theorem}

We will use the following description of (co)cartesian unstraightening from \cite[Theorem 7.4]{GHN}:

\begin{theorem} \label{thm:lax_colim_unstr}
	For a functor $F\colon B\rt \Cat$, there is a natural equivalence
	\[ 
	\Unco(F) \simeq \colim\Big( \TwR(B) \xrightarrow{(s,t)} B \times B^{\op} 
	\xrightarrow{F \times B_{-/}} \Cat\Big).
	\]
	Dually, for a functor $F\colon B^{\op}\rt \Cat$, there is a natural equivalence
	\[ \Uncart(F) \simeq 
	\colim \Big(\TwR(B^{\op}) \xrightarrow{(s, t)} B^{\op} \times B
	\xrightarrow{F \times B{/-}} \Cat\Big).
	\]
\end{theorem}

We will apply this in the setting of \ref{not:fibration of adtrip} to compute the cocartesian fibration classified by $\Span(F)\colon X^{\op}\rt \Cat$. We find that $\Unco(\Span(F))$ is equivalent to the colimit of the diagram
\[\begin{tikzcd}[column sep=3pc]
\TwR(X^{\op}) \arrow[r, "{(s,t)}"] & X^{\op} \times X
\arrow[rr, "{\Span(F) \times (X^{\op})_{-/}}"] & &  \Cat.
\end{tikzcd}\]
Now note that there are natural equivalences ${-/}X^{\op}\simeq (X{/-})^{\op} \simeq \Span(X{/-},\iota(X{/-}),X{/-})$ and recall that $\Span$ commutes with products of adequate triples to conclude that this colimit is the same as that of the diagram
\[
\TwR(X^{\op}) \xrightarrow{(s,t)} X^{\op} \times X
\xrightarrow{F \times (X{/-},\iota(X{/-}),X{/-})} \AdTrip \xrightarrow{\Span} \Cat.
\]
We shall prove \Cref{thm:span_funct_cocart} by showing that the colimit of the first two functors in this composition is precisely $Y^{\glob}=(Y,Y_{\ing}^{\fb},Y^{\eg})$ and that this colimit is preserved by $\Span$.

We can directly apply the equivalence of Theorem \ref{thm:lax_colim_unstr} to conclude that
\begin{equation}\label{eq:unstr as colims}
Y\simeq \colim F(-)\times X{/-},\quad Y^{\eg}\simeq \colim F^\eg(-) \times X{/-},  \quad \text{ and }\quad Y_{\ing} \simeq \colim F_{\ing}(-)\times X{/-}
\end{equation}
where all colimits are taken over $\TwR(X^{\op})$. To exhibit $Y_{\ing}^{\fb} = Y_{\ing}\times_X \iota X$ as a colimit, we need the following lemma:
\begin{lemma}\label{lem:core as colim}\label{lem:fibrewise ing as colim}
Let $F\colon B^{\op}\rt \Cat$ be a functor. Then the functor $\Uncart(F)\times_B \core B\rt \Uncart(F)$ is naturally equivalent to the functor
\begin{equation}\label{eq:restrict to core}
\colim_{\TwR(B^{\op})} F\times \iota(B{/-})\rt \colim_{\TwR(B^{\op})} F\times B{/-}.
\end{equation}
\end{lemma}
\begin{proof}
Consider the full subcategory of $\Fun\big(\Fun(B^{\op}, \Cat), \Cat\big){/\Uncart}$ spanned by those natural transformations $\Phi\rt \Uncart$ with the following two properties:
\begin{enumerate}
\item\label{it:core on pt} At the terminal diagram, $\Phi(\ast)\rt \Uncart(\ast)\simeq B$ exhibits the inclusion of the core of $B$.
\item\label{it:base change} For each $F\colon B^{\op}\rt \Cat$, the induced square
\begin{equation}\label{eq:naturalitysquare}
\begin{tikzcd}
\Phi(F) \arrow[r] \arrow[d] & \Uncart(F) \arrow[d] \\
\Phi(\ast) \arrow[r, hook]                                & B.                   
\end{tikzcd}
\end{equation}
is cartesian.
\end{enumerate}
This subcategory is contractible and the natural transformation $\Uncart(-)\times_B \core B\rt \Uncart(-)$ is (by definition) an object in there. It will therefore suffice to show that the composite natural transformation 
$$
\Phi(F) \coloneqq \colim_{\TwR(B^{\op})}\big(F\times \core(B{/-})\big)\rt \colim_{\TwR(B^{\op})}\big(F\times B{/-}\big) \simeq \Uncart(F)
$$
satisfies conditions \ref{it:core on pt} and \ref{it:base change} as well. 

Let us start by computing $\Phi(F)$. To this end, let us write $\Ar(B)_s=\Ar(B)\times_{B} \core B$, where the structure map $\Ar(B) \rightarrow B$ is evaluation at the source. Then $t\colon \Ar(B)_s\rt B$ is the left fibration classified by the functor $b\mapsto \core(B{/b})$ and $s\colon \Ar(B)_s\rt \core B$ admits a fully faithful left adjoint $\cst\colon \core B\hookrightarrow \Ar(B)_s$ taking degenerate arrows. Let us now consider the following two pullbacks of \icats{}
$$\begin{tikzcd}
D\arrow[r, hookrightarrow]\arrow[d, "\pi"{swap}] & C\arrow[r, "q"]\arrow[r] \arrow[d] & \TwR(B^{\op})\arrow[d, "t"]\\
\core B\arrow[r, hookrightarrow, "\cst"] & \Ar(B)_s\arrow[r, "t"] & B
\end{tikzcd}$$
as well as the following diagram of \icats{}
$$\begin{tikzcd}
F'\colon C\arrow[r,"q"] & \TwR(B^{\op})\arrow[r, "s"] & B^{\op}\arrow[r, "F"] & \Cat.
\end{tikzcd}$$
Since $q\colon C\to \TwR(B^{\op})$ is the left fibration classifying $\core(B{/-})\colon \TwR(B^{\op})\to B\to \sS$ and $F'$ is constant along the fibres of $q$, the left Kan extension $q_!F'$ is naturally equivalent to the diagram $F(-)\times \core(B{/-})$ of which we want to compute the colimit \cite[Proposition 4.3.3.10]{HTT}.
By transitivity of Kan extensions, it therefore suffices to compute $\colim_C F'$.

Note that the vertical arrows in the above diagram are all cartesian fibrations, so that the fully faithful inclusion $D\hookrightarrow C$ admits a (localising) right adjoint  \cite[Corollary 5.2.7.11]{HTT}
$$\begin{tikzcd}
p=(\id, s)\colon C=\TwR(B^{\op})\times_B \Ar(B)_s\arrow[r] & \TwR(B^{\op})\times_B \iota B=D.
\end{tikzcd}$$
By construction, the diagram $F'$ only depended on the first factor $\TwR(B^{\op})$; Consequently, it factors as $F'= F''p$, for some $F''\colon \DD\rightarrow \Cat$. Because $p$ is a localisation it is in particular cofinal, so we conclude that $\colim_C F'\simeq \colim_D F''$.

Let us first compute the left Kan extension $\pi_!F''$ of $F''$ along the cocartesian fibration $\pi\colon D\to \core B$. We can compute this Kan extension as a colimit over the fibres, which are equivalently the fibres of $t\colon \TwR(B^{\op})\to B$. Unravelling the definitions then shows that
$$
(\pi_!F'')(b)=\colim\big(B^{\op}/b\rt B^{\op}\rto{F} \Cat\big) \simeq F(b)
$$
since $B^{\op}/b$ has a terminal object. We conclude that there is a natural equivalence
$$
\Phi(F)=\colim_{\TwR(B^{\op})} F\times \iota(B{/-})\simeq \colim_{D} F'' \simeq \colim_{\core B} \pi_!F'' \simeq \colim_{\core B} F\simeq \Uncart(F)\times_B \core B.
$$
In particular, we see that $\Phi(\ast)$ is an \igpd{}, so that \ref{it:core on pt} $\Phi(\ast)\to B$ is a core inclusion, being a colimit of core inclusions. For \ref{it:base change}, note that the vertical maps in \eqref{eq:naturalitysquare} are cartesian fibrations, so that it suffices to show that the maps on vertical fibres are equivalences. By naturality of the transformation $\Phi\rt \Uncart$ in the base $B$, this map on fibres map is simply the map when $B=*$. In this case the map \eqref{eq:restrict to core} is clearly an equivalence.
\end{proof}
\begin{proposition}\label{prop:colimit triple}
The diagram 
\[\Psi\colon \TwR(X^{\op}) \xrightarrow{(s,t)} X^{\op} \times X
\xrightarrow{F \times (X{/-},\iota(X{/-}),X{/-})} \AdTrip\]
has colimit given by the adequate triple $Y^{\glob}=(Y, Y_\ing^{\fb}, Y_\eg)$.
\end{proposition}
\begin{proof}
The equivalences \eqref{eq:unstr as colims} and \Cref{lem:fibrewise ing as colim} show that the colimit of $\Psi$ in $\Fun(\Lambda_2[2], \Cat)$ is given by the triple $Y^{\glob}$. To see that this is also a colimit in the subcategory of adequate triples, it suffices to verify that the following are equivalent for a functor $\phi\colon Y^{\glob}=\colim \Psi\rt Z$:
\begin{enumerate}
\item $\phi$ preserves ambigressive pullback squares.
\item the composite $\Psi(f^{\op}\colon x_1\to x_0)\to \colim \Psi\to Z$ does so for each $f^{\op}\in \TwR(X^{\op})$. 
\end{enumerate}
Now note that (on underlying \icats{}) the map $\Psi(f^{\op}\colon
x_1\to x_0)\rt \colim \Psi$ can be identified with the functor
$\lambda_f\colon F(x_1)\times X{/x_0}\rt Y$ sending $(y, g\colon
x_2\to x_0)$ to $g^*f^*(y)$. Because $g^*$ is a map of adequate triples for every map $g$ in $X$, $\lambda_f$ is a map of adequate triples, which shows that (1) implies
(2).

For the converse, note that every ambigressive pullback square in $Y$ can be obtained as a pasting of two types of squares: (a) the pullback of a fibrewise ingressive map along a $p$-cartesian lift of some $g\colon x_2\to x_0$ and (b) the pullback of an ingressive and egressive map in a single fibre $Y_{x_0}\simeq F(x_0)$. It suffices to show that each of these squares is the image of an ambigressive pullback square under some $\lambda_f$. For (a), we can take the image under $\lambda_{\mm{id}_{x_0}}$ of the square in $F(x_0)\times X{/x_0}$ formed by an ingressive arrow in $F(x_1)$ and the arrow $g\to \mm{id}_{x_0}$ in $X/x_0$. For (b), we simply take the image under $\lambda_{\mm{id}_{x_0}}$ of an ambigressive square in $F(x_0)\times \{\mm{id}_{x_0}\}$. 
\end{proof}
Applying $\Span$ to the colimiting cocone of \Cref{prop:colimit triple}, we therefore obtain a natural cocone of the diagram \[
\TwR(X^{\op}) \xrightarrow{(s,t)} X^{\op} \times X
\xrightarrow{F \times (X{/-},\iota(X{/-}),X{/-})} \AdTrip \xrightarrow{\Span} \Cat
\]
whose tip is the $\infty$-category $\Span(Y^{\glob})$. This induces a natural functor 
\[
\alpha_F\colon \Unco(\Span(F))\simeq \colim\limits_{\TwR(X^{\op})}\Big(\Span(F) \times (X^{\op})_{-/} \Big)\rt \Span(Y^{\glob})=\Span\big(\Uncart(F)^{\glob}\big).
\] 

Note that the construction of $\alpha_F$ is functorial in the diagram $F$. Applied to the terminal diagram $F=\ast$, it yields the identity map of $\Span(X, \iota X, X)$. Using the identification $\Span(X, \iota X, X) \simeq X^{\op}$ from \Cref{prop:spanswithcore}, we then obtain the desired natural transformation $\alpha\colon  \Unco(\Span\circ (-))  \Longrightarrow\Span(\Uncart(-)^{\glob})$ between functors $\Fun(X^{\op}, \AdTrip)\rt \Cat{/X^{\op}}$.

\begin{lemma}\label{lem:Phi_preserves_cocart}
For each $F\colon X^{\op}\rt \AdTrip$, the functor $\alpha_F$ preserves cocartesian edges over $X^\op$.
\end{lemma} 
\begin{proof}
Note that $F$ comes with a natural transformation $\core F\rt F$, where each $\core F(x)$ comes equipped with the (only possible) adequate triple structure where all morphisms are both ingressive and egressive. Unstraightening, this identifies $\Uncart(\core F)^{\glob}$ with the subcategory of $\Uncart(F)$ on the cartesian arrows, with all arrows egressive and only the equivalences being ingressive.

The map $\core F\rt F$ induces a natural transformation
$\core\Span(F)\simeq \Span(\core F)\rt \Span(F)$, where the first equivalence uses that $\Span(\core F)\simeq \core(F)$ (\Cref{prop:spanswithcore}). By naturality, we then obtain a commuting diagram
$$\begin{tikzcd}
\Unco\big(\core \Span(F)\big)\arrow[d]\arrow[r, "\alpha_{(\core F)}"] & \Span(\Uncart(\core F)^{\glob})\arrow[d]\\
 \Unco\big(\Span(F)\big)\arrow[r, "\alpha_F"] &\Span(\Uncart(F)^{\glob}).
\end{tikzcd}$$
The result now follows since both vertical arrows can be identified with the inclusion of the wide subcategory of cocartesian arrows: For the left arrow this is evident and for the right arrow, this follows from our description of $\Uncart(\core F)^{\glob}$ and \Cref{barwick_cocart_equiv}.
\end{proof}

\begin{proof}[Proof of Theorem \ref{thm:span_funct_cocart}]
We have seen that $\alpha$ constitutes a natural transformation between
functors into $\Cat/X^{\op}$. By Lemma \ref{lem:Phi_preserves_cocart},
the value of $\alpha$ at a diagram $F$ is a map between cocartesian
fibrations over $X^\op$ which preserves cocartesian edges. It then
suffices to show that $\alpha_F$ is an equivalence when restricted to
every fibre \cite[Corollary 2.4.4.4]{HTT}. Since $\Span$ preserves pullbacks (in particular fibres), we can thus reduce to the case $X=\ast$, where the map is
equivalent to the identity by inspection.
\end{proof}

\begin{remark}\label{rem:Span_cart_fib}
Given a functor $G\colon X^\op\rt \Cat$, one way to write down the cartesian fibration encoding it is to write down the cocartesian fibration which classifies $(-)^\op\circ G$, and then take its opposite. In the case that $G=\Span(F)$, Lemma \ref{span_op} shows that $G^\op = \Span(F^\rev)$. Therefore the previous result allows us to describe its cocartesian straightening. Finally we can again apply Lemma \ref{span_op} to describe the opposite of the result. In total one finds that 
\[
\Span(p)\colon\Span(Y,Y_\ing,Y^\eg_{\fb})\rt \Span(X,X,\core X) \simeq X
\] is a cartesian fibration which classifies the functor $\Span(F)\colon X^\op\rt \Cat$. 
\end{remark}

As a special case we recover the key insight of \cite{BGN}.
\begin{corollary}\label{BGNidentification}
	Let $p\colon Y \rt X$ be a cartesian fibration and define $Y^\vee \rt X^\op$ to be the functor \[ \Span(p)\colon \Span(Y,Y\times_X\iota X,Y^{\dagger}) \rt \Span(X,\iota X, X) \simeq X^{\op} \] 
	where $Y^\dagger$ is the subcategory of $p$-cartesian morphisms (cf.\ \Cref{prop:adequate triple}). This gives rise to an equivalence $\SDualco \colon \Cart(X)\rt \Cocart(X^\op)$ which fits into a commuting triangle
	\[
	\begin{tikzcd}
		\Cart(X) \arrow[rr, "\SDualco"] \arrow[rd, "\Strcart"'] &                  & \Cocart(X^{\op}) \arrow[ld, "\Strco"] \\
		& \Fun(X^{\op},\Cat). &                                         
	\end{tikzcd}
	\]
\end{corollary}

\begin{proof}
Note that once we have exhibited the triangle as commutative, the two-out-of-three property will immediately imply that $\SDualco$ is an equivalence. We will provide a natural equivalence $\SDualco(\Uncart(F))\rt \Unco(F)$ for functors $F\colon X^{\op}\rt \Cat$. To this end, consider the functor $\delta\colon \Cat\rt \AdTrip$ sending $C\mapsto (C, C, \iota C)$ and note that $\Span\circ \delta\simeq \mm{id}$. One then has natural equivalences 
$$
\SDualco(\Uncart(F))\simeq\Span\big(\Uncart(\delta(F))^{\int}\big)\simeq \Unco\big(\Span(\delta(F))\big)\simeq \Unco(F)
$$ where the middle equivalence is \Cref{thm:span_funct_cocart}.
\end{proof}

As another special case, we may consider the universal example. 
\begin{construction}\label{con:lax slice}
Let us define the \emph{lax under-category} $\ast\laxslice\AdTrip$ to be the domain of the cocartesian fibration classifying the functor $(-)^{\op}\colon \AdTrip\rt \Cat$, see the next remark for a justification. One can identify an object of $\ast\laxslice\AdTrip$ with a tuple $(X, x)$ of an adequate triple and an object $x\in X$, and a map $(X, x)\to (Y, y)$ with a map $f\colon X\to Y$ and a map $\mu\colon y\to f(x)$ in $Y$. We will say that $(f, \mu)$ is egressive if $\mu$ is egressive, and ingressive if $\mu$ is ingressive and $f$ is an equivalence.
\end{construction}
\begin{remark}
If we write $*\laxslice\Cat$ for the domain of the cocartesian fibration classified by $(-)^\op\colon \Cat\rt \Cat$, then $\ast\laxslice\AdTrip = \ast\laxslice\Cat \times_{\Cat} \AdTrip$. We will show in Section \ref{sec:OplaxAr} (see specifically Remark \ref{rmk:universal_cart_fib}) that $\ast\laxslice\Cat$ really is equivalent to the lax under-category, defined a priori via a construction in $(\infty,2)$-category theory, justifying the notation. 
\end{remark}

\begin{corollary}
Construction \ref{con:lax slice} defines the structure of an adequate triple on $\big(\ast\laxslice\AdTrip\big)^{\op}$ and the map
$$
\Span(p)\colon \Span\big((\ast\laxslice\AdTrip)^{\op}\big)\rt \Span(\AdTrip^{\op}, \iota\AdTrip^{\op}, \AdTrip^{\op})\simeq \AdTrip
$$
is a cocartesian fibration, classified by the functor $\Span\colon \AdTrip\rt \Cat$.
\end{corollary}

\begin{proof}
Apply Theorem \ref{thm:span_funct_cocart} to the identity functor on $\AdTrip$ and observe that the adequate triple structure from \ref{not:fibration of adtrip} corresponds precisely to the one from \Cref{con:lax slice} under taking opposite categories.
\end{proof}

\begin{remark}
Applying the conclusions of Remark \ref{rem:Span_cart_fib}, we obtain that if we declare a map $(f, \mu)$ in $\ast\laxslice\AdTrip$ to be ingressive if $\mu$ is ingressive and egressive if $\mu$ is egressive and $f$ is an equivalence, then we obtain a \emph{different} adequate triple structure on $(\ast\laxslice\AdTrip)^{\op}$, which gives rise to a cartesian fibration $\Span((\ast\laxslice\AdTrip)^{\op})\rt \AdTrip^{\op}$, classifying the functor $\Span\colon \AdTrip\rt \Cat$.
\end{remark}

\begin{remark}
In \cite{Nardin}, Nardin proves that for an orbital $\infty$-category $T$ the $T$-$\infty$-category of
$T$-commutative monoids in a $T$-$\infty$-category $C$ with finite
$T$-products is equivalent to the $T$-$\infty$-category
$\ul{\Fun}_T^\times(\ul{\mathbf{A}}^{\eff}(T),C)$, where
$\ul{\mathbf{A}}^{\eff}(T)$ is a cocartesian fibration whose fibre over
$t\in T$ is equivalent to $\Span((\mathbf{F}_{T}){/t})$, where
$\mathbf{F}_T$ is the finite coproduct completion of $T$. This
indicates that $T$-commutative monoids are computed by a
generalisation of Mackey functors. By inspecting the definition of
$\ul{\mathbf{A}}^{\eff}(T)$ from \cite[Definition 4.10]{Nardin}, one observes
that $\ul{\mathbf{A}}^{\eff}(T) = \Span(Y^{\glob})$, where $Y$ is the
cartesian unstraightening of the functor underlying \[T^\op\rightarrow
  \AdTrip,\quad t \mapsto ((\mathbf{F}_{T}){/t},(\mathbf{F}_{T}){/t},(\mathbf{F}_{T}){/t}),\] where the functoriality is given by pullback. Therefore Theorem \ref{thm:span_funct_cocart} shows that the functoriality of the $T$-$\infty$-category $\ul{\Fun}_T^\times(\ul{\mathbf{A}}^{\eff}(T),C)$ is induced by the standard functoriality of $\Span(-)$.
\end{remark}

\section{Orthogonal adequate triples and their duals}\label{sec:orthtrip}

Work of Barwick, Glasman and Nardin uses the span
category construction from \Cref{sec:span}
to produce an equivalence between cartesian fibrations and cocartesian
fibrations of \icats{} \cite{BGN}. In the previous section we provided an alternative proof of this fact by exhibiting it as a specific example of a general phenomenon: The cocartesian fibration classifying
$\Span\circ F$ is given by the span construction applied to a particular adequate triple structure on the cartesian unstraightening of $F$.

In this section, we will describe a different perspective on this result, which is a priori entirely independent from the discussion of (co)cartesian fibrations: The construction of span \icats{} upgrades to an automorphism (of order $2$) of a certain full subcategory of adequate triples spanned by those triple that we will
call \emph{orthogonal}. In Section \ref{sec:fibs} we will show how this construction interacts with various notions of (cartesian) fibrations between such orthogonal adequate triples. In particular we will extend Lanari's results on the relation between cartesian fibrations and factorisation systems from \cite{LanariCart}, and thus
recover the main result from \cite{BGN} in a different way.

To define the notion of an orthogonal adequate triple, let us start by briefly recalling Joyal's notion of a factorisation system on an \icat{} $\C$ (also called an \emph{orthogonal} factorisation system), following \cite[Section 5.2.8]{HTT}. It consists of two classes of arrows $\C_\lcl$ (left) and $\C_\rcl$ (right), denoted $\twoheadrightarrow$ and $\rightarrowtail$, such that
\begin{enumerate}[(i)]
  \item both classes of maps are closed under
      retracts,
\item each solid commutative diagram
\[\begin{tikzcd} x \ar[r] \ar[d,"f"', twoheadrightarrow] & y \ar[d,"g", rightarrowtail] \\
                 x' \ar[r]\ar[ru,dashed] & y'\end{tikzcd}\]
admits a unique dashed filler whenever $f \in \C_\lcl$ and $g \in \C_\rcl$, and
\item every map in $\C$ factors as $\begin{tikzcd}[cramped] {\cdot}\arrow[r, twoheadrightarrow, "f"] & {\cdot} \arrow[r, rightarrowtail, "g"] & {\cdot}\end{tikzcd}$ with $f \in \C_\lcl$ and $g \in \C_\rcl$.
\end{enumerate}
Let us point out that the notation for the left and right classes is inspired by the notation for the epi-mono factorisation system, and is \emph{opposite} to the convention typically employed for model categories.

Both $\C_\lcl$ and $\C_\rcl$ then define subcategories of $\C$ that contain all equivalences and each class uniquely determines the other \cite[Proposition 5.2.8.6 \& 5.2.8.11]{HTT}.
\begin{remark}\label{rem:functorial factorisation}
The factorisation in the last item is essentially unique, and this
uniqueness is in fact equivalent to the first two items \cite[Proposition 5.2.8.17]{HTT}. More precisely, writing $\cat{Fact}(C)\subseteq \Fun([2], C)$ for the full subcategory of functors sending $0\leq 1$ to $\C_\lcl$ and and $1\leq 2$ to $\C_\rcl$, one has that restriction to the arrow $0\leq 2$ defines an equivalence
\[\begin{tikzcd}
\circ \colon \cat{Fact}(C)\arrow[r, "\sim"] & \Ar(C).
\end{tikzcd}\]
An (essentially unique) section provides a \emph{functorial} factorisation of morphisms in $C$.
\end{remark}

\begin{definition}\label{def:orthoad}
We call an adequate triple $(X,X_\ing,X^\eg)$ \emph{orthogonal} if every ambigressive square in $X$ is cartesian, and the egressives ($\twoheadrightarrow$) and ingressives ($\rightarrowtail$) are the left and right classes respectively of an orthogonal factorisation system on $X$. We denote the full subcategory of orthogonal adequate triples by $\AdTrip^\perp$.
\end{definition}
\begin{remark}
Since every ambigressive square in an orthogonal adequate triple is automatically cartesian, it follows that $\AdTrip^{\perp}$ is a \emph{full} subcategory of $\Fun(\Lambda_2[2], \Cat)$ (in contrast to $\AdTrip$).
\end{remark}

The simplest examples are $(A,\iota A,A)$ and $(A,A,\iota A)$ for any \icat{} $A$. We will generate more interesting examples by means of the following observation, which follows directly from the mapping properties of cartesian edges:

\begin{proposition}\label{lem:carttriport}
If $X$ is an orthogonal adequate triple and $p\colon Y\rt X$ has all $p$-cartesian lifts over $X_\ing$, then the adequate triple $\big(Y, Y^\dagger, p^{-1}(X^{\eg}))$ from Proposition \ref{prop:adequate triple} is again orthogonal. 
\end{proposition}

In particular, for a cartesian fibration $p \colon Y \rightarrow X$ the triple structure on $Y$ with ingressives the $p$-cartesian edges and egressives the fibrewise maps is orthogonal; in this way our constructions contain those from \cite{BGN}.

\begin{proof}
	The fact that any ambigressive square is cartesian follows
        immediately from the proof of Proposition \ref{prop:adequate
          triple}. To show that the ingressives and egressives form a
        factorisation system on $Y$, let us start by noting that both
        classes are stable under retracts (for the cartesian morphisms, this uses that they are characterised by a lifting property relative to the base X). Next consider a commutative square
	 \[
	 \begin{tikzcd}
		 x \arrow[r] \arrow[d, "\phi"', two heads] & y \arrow[d, "\psi", tail] \\
		 x' \arrow[r] \arrow[ru, dotted]   & y',                         
	 \end{tikzcd}
	 \]
where $\phi\in Y^{\eg}$ and $\psi\in Y_{\ing}$. We need to exhibit a unique dotted arrow making the diagram commute. After applying $p$ there exists a unique dotted edge in $X$ making the square commute. Since $\psi$ is $p$-cartesian we find that there is a unique dotted edge in $Y$ which makes the square commute. Finally, given an edge $f\colon x\rightarrow y$ in $Y$, one can factor $p(f)=gh$ with $g\in X_{\ing}$ and $h\in X^{\eg}$. Taking a $p$-cartesian lift $g'$ of $g$, we can factor $f= g'h'$ with $h'\in p^{-1}(X^\eg)$, so that $Y^{\dagger}$ and $p^{-1}(X^{\eg})$ form a factorisation system.
\end{proof}
\begin{example}\label{ex:product}
A product $A\times B$ of two \icats{} admits the structure of an orthogonal adequate triple, denoted $(A, B)^{\perp}$, with ingressives $A\times \core B$ and egressives $\core A\times B$.
\end{example}

\begin{example}\label{exrecoll}
Let $X$ be a stable \icat{}, together with reflective stable subcategories exhibiting $X$ as a recollement \cite[Definition A.8.1]{HA}
$$\begin{tikzcd}
U\arrow[r, hookrightarrow, "j_*"{swap}, yshift=-0.7ex] & X \arrow[l, "j^*"{swap}, yshift=0.7ex]\arrow[r, "i^*", yshift=0.7ex] \arrow[r, hookleftarrow, "i_*"{swap}, yshift=-0.7ex] & Z;
\end{tikzcd}$$
explicitly, this means that $i^*$ and $j^*$ jointly detect equivalences and $j^*i_*$ sends every object to the terminal object, but there are several other characterisations, see \cite[Appendix A.2]{CDH2}. In particular, by \cite[Lemma A.2.5]{CDH2} $i_*$ admits a further right adjoint $i^!$ such that $i_*i^!x\rt x\rt j_*j^*x$ is a fibre sequence for each $x\in X$.

We now claim that $X$ has the structure of an orthogonal adequate triple, in which a map is ingressive or egressive if its image under $i^!$  or $j^*$ is an equivalence, respectively. By \cite[Corollary 2.6.1]{CDH2}
 $j^*$ is a cartesian fibration with an edge in $X$ being $j^*$-cartesian if and only if it is an $i^!$-equivalence. Thus \cref{lem:carttriport} (applied with target $(U,U,\iota U)$) yields the structure of an orthogonal triple on $X$ as desired.
\end{example}

\begin{example}\label{ex:arrow triple}
For each \icat{} $A$, the source projection $s\colon \Ar(A)\rt A$ is a cartesian fibration. Applying \Cref{lem:carttriport} to the triple $(A,A,\iota A)$, we find that $\Ar(A)$ admits the structure of an orthogonal adequate triple in which a map is egressive (ingressive) if its image under the source map $s\colon \Ar(A)\to A$ (target map $t\colon \Ar(A)\to A$) is an equivalence. For every $A\to B$, the induced functor $\Ar(A)\rt \Ar(B)$ preserves ingressive and egressive maps, so that we obtain a functor
$$\begin{tikzcd}
\Ar\colon \Cat\arrow[r] & \AdTrip^{\perp}.
\end{tikzcd}$$
\end{example}
In fact, \Cref{ex:arrow triple} provides the \emph{free} orthogonal adequate triple generated by $A$:
\begin{proposition}\label{prop:free triple}
The functor $\Ar\colon \Cat\rt \AdTrip^{\perp}$ is the left adjoint to the forgetful functor $\AdTrip^{\perp}\rt \Cat$ taking underlying \icats{}. More precisely, the degeneracy map $\cst\colon A\rt \Ar(A)$ provides the unit transformation exhibiting $\Ar$ as the left adjoint to the forgetful functor.
\end{proposition}
\begin{proof}
It will suffice to define a natural counit map and provide natural homotopies for the triangle identities. Given an adequate triple $X$, let $\epsilon_X\colon \Ar(X)\rt X$ be the unique functor that fits into a commuting diagram
$$\begin{tikzcd}
\cat{Fact}(X)\arrow[d, "\circ"{swap}, "\sim"]\arrow[r, hookrightarrow] & \Fun([2], X)\arrow[d, "\mm{ev}_1"]\\
\Ar(X)\arrow[r, dotted, "\epsilon_X"] & X.
\end{tikzcd}$$
Here the left vertical functor is the equivalence from Remark \ref{rem:functorial factorisation}, a section of which provides a functorial egressive-ingressive factorisation. Since all three solid maps are functorial in the adequate triple $X$, the map $\epsilon_X$ is functorial in $X$. Explicitly, the value of $\epsilon_X$ on a map $\mu\colon f\to g$ in $\Ar(X)$ given by the square
\[
\begin{tikzcd}
\cdot \arrow[r, "\mu_0"] \arrow[d, "f"] & \cdot \arrow[d, "g"] \\
\cdot \arrow[r, "\mu_1"] &	\cdot
\end{tikzcd}
\] 
is given by the middle vertical arrow in the unique diagram
$$\begin{tikzcd}[cramped, row sep=1.1pc]
\cdot\arrow[dd, bend right=30, "f"{swap}]\arrow[r, "\mu_0"]\arrow[d, two heads] & \cdot \arrow[d, two heads]\arrow[dd, bend left=30, "g"]\\
\cdot\arrow[d, rightarrowtail]\arrow[d, rightarrowtail]\arrow[r] & \cdot\arrow[d, rightarrowtail]\\
\cdot\arrow[r, "\mu_1"] & \cdot \end{tikzcd}$$
From this we see that $\epsilon_X$ is a map of adequate triples. Indeed, if $\mu$ is ingressive (i.e.\ $\mu_1$ is an equivalence), then the middle horizontal map is ingressive by the right cancellation property for ingressives \cite[Proposition 5.2.8.6 (3)]{HTT}, and likewise in the egressive case where $\mu_0$ is an equivalence.

It remains to provide the triangle identities. If $X$ is an adequate triple, then the composite
$$\begin{tikzcd}
X\arrow[r, "\cst"] & \Ar(X)\arrow[r, "\epsilon_X"] & X
\end{tikzcd}$$
is evidently naturally equivalent to the identity (using that the functorial factorisation $\Ar(X)\to \cat{Fact}(X)$ sends degenerate arrows to constant diagrams). Furthermore, for an \icat{} $A$, the composite
$$\begin{tikzcd}
\Ar(A)\arrow[r, "\cst"] & \Ar(\Ar(A))=\Fun([1]\times [1], A)\arrow[r, "\epsilon_{\Ar(X)}"] & \Ar(A)
\end{tikzcd}$$
sends an arrow $f$ to $f\rto{\sim} f$ and then applies the ingressive-egressive factorisation to this natural equivalence. This is again naturally equivalent to the identity.
\end{proof}

\begin{proposition}\label{prop:dualadequatetrip}
  For any adequate triple $(X,X_\ing,X^\eg)$, the subcategory
  inclusions $(X^\eg)^\op, X_\ing \to \Span(X)$ from
  \Cref{prop:spanswithcore} give an orthogonal factorisation system on
  $\Span(X)$. If the triple $(X,X_\ing,X^\eg)$ is furthermore
  orthogonal, then $(\Span(X), X_\ing, (X^\eg)^{\op})$ is again an
  orthogonal adequate triple.
\end{proposition}

\begin{proof}
  To prove the first claim we use the criterion from \Cref{rem:functorial factorisation}, i.e.\ we show that the map
\[\begin{tikzcd}
\mm{\circ}\colon \cat{Fact}(\Span(X,X_\ing,X^\eg))\arrow[r] & \Ar(\Span(X,X_\ing,X^\eg)).
\end{tikzcd}\]
is an equivalence. We first observe that this map is an equivalence on cores for all adequate triples: This assertion translates to the statement that the \igpd{} of dotted extensions of the solid diagram
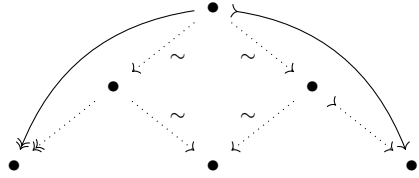
\begin{center}
\begin{tikzcd}
        &                                                             & \bullet \arrow[rd, "\sim"', dotted] \arrow[ld, "\sim", dotted] \arrow[lldd, two heads, bend right] \arrow[rrdd, tail, bend left] &                                                                   &         \\
        & \bullet \arrow[ld, dotted, two heads] \arrow[rd, "\sim", dotted] &                                                                                                                                  & \bullet \arrow[ld, "\sim"', dotted] \arrow[rd, tail, dotted] &         \\
\bullet &                                                             & \bullet                                                                                                                          &                                                                   & \bullet
\end{tikzcd} 
\end{center} 
is contractible, which is obvious. 

To conclude, we now observe that it suffices to show that
\[\begin{tikzcd}\mm{\circ}_* \colon \Hom_{\Cat}([n],\cat{Fact}(\Span(X,X_\ing,X^\eg))) \arrow[r] & \Hom_{\Cat}([n],\Ar(\Span(X,X_\ing,X^\eg)))\end{tikzcd}\]
is an equivalence for all $n \geq 0$ (in fact $n=1$ suffices, but this will not simplify the argument). But by Corollary \ref{FunofSpan} we can rewrite the target as
\[\Hom_{\Cat}([n] \times [1],\Span(X,X_\ing,X^\eg)) \simeq \core\Ar(\Span(\Q_n(X,X_\ing,X^\eg))).\]
Similarly, the source is equivalent to a set of path components in $\Hom_{\Cat}([2],\Span(\Q_n(X,X_\ing,X^\eg)))$, and we claim that it precisely corresponds to $\core\cat{Fact}(\Span(\Q_n(X,X_\ing,X^\eg)))$. Thus the statement on groupoid cores applied to the various $\Q_n(X,X_\ing,X^\eg)$ gives the statement in full. To verify the remaining claim one unwinds definitions to find both sides spanned by those diagrams $\TwR[n] \times \TwR[2] \rightarrow X$ whose restriction along the Segal maps
\[\begin{tikzcd}
\bullet                                              & \bullet \arrow[r, tail] \arrow[l, two heads]                                                                              & \bullet                                              & \bullet \arrow[r, tail] \arrow[l, two heads]                                                                              & \bullet                                                          && \bullet                                              & \bullet \arrow[r, tail] \arrow[l, two heads]                                                                              & \bullet\\
\bullet \arrow[u, tail] \arrow[d, two heads] & \bullet \arrow[d, two heads] \arrow[l, two heads] \arrow[u, tail] \arrow[r, tail] & \bullet \arrow[u, tail] \arrow[d, two heads] & \bullet \arrow[d, two heads] \arrow[l, two heads] \arrow[u, tail] \arrow[r, tail] & \bullet \arrow[d, two heads] \arrow[u, tail] & \cdots & \bullet \arrow[u, tail] \arrow[d, two heads] & \bullet \arrow[d, two heads] \arrow[l, two heads] \arrow[u, tail] \arrow[r, tail] & \bullet \arrow[d, two heads] \arrow[u, tail] \\
\bullet                                              & \bullet \arrow[r, tail] \arrow[l, two heads]                                                                              & \bullet                                              & \bullet \arrow[r, tail] \arrow[l, two heads]                                                                              & \bullet                                                          && \bullet                                              & \bullet \arrow[r, tail] \arrow[l, two heads]                                                                              & \bullet\\
\bullet \arrow[u, two heads] \arrow[d, tail] & \bullet \arrow[d, tail] \arrow[l, two heads] \arrow[u, two heads] \arrow[r, tail] & \bullet \arrow[u, two heads] \arrow[d, tail] & \bullet \arrow[d, tail] \arrow[l, two heads] \arrow[u, two heads] \arrow[r, tail] & \bullet \arrow[d, tail] \arrow[u, two heads] & \cdots & \bullet \arrow[u, two heads] \arrow[d, tail] & \bullet \arrow[d, tail] \arrow[l, two heads] \arrow[u, two heads] \arrow[r, tail] & \bullet \arrow[d, tail] \arrow[u, two heads] \\
\bullet                                              & \bullet \arrow[l, two heads] \arrow[r, tail]                                                              & \bullet                                              & \bullet \arrow[l, two heads] \arrow[r, tail]                                                              & \bullet             && \bullet                                              & \bullet \arrow[l, two heads] \arrow[r, tail]                                                              & \bullet          
\end{tikzcd}\]
has all ambigressively marked squares cartesian and all upwards pointing maps equivalences.

Now assume the adequate triple  $(X,X_\ing,X^\eg)$ is orthogonal.
To prove that $(\Span(X), X_\ing, (X^\eg)^{\op})$ is adequate, we must show that $\Span(X)$ admits
ambigressive pullbacks. By Lemma \ref{span_op}, we may equivalently show that $\Span(X^\rev)$ admits ambigressive pushouts, which turns out to be notationally slightly more convenient. Let us start by observing that by adjunction a commuting square
in $\Span(X^\rev)$ with horizontal arrows in $(X_{\ing})^\op$ and vertical arrows in $X^{\eg}$ corresponds to a diagram of the form
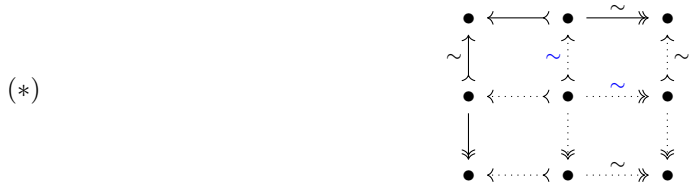
\begin{equation}\tag{$\ast$}
\begin{tikzcd}
\bullet                                              & \bullet \arrow[r, "\sim", two heads] \arrow[l, tail]                                                                              & \bullet                                                               \\
\bullet \arrow[u, "\sim", tail] \arrow[d, two heads] & \bullet \arrow[d, dotted, two heads] \arrow[l, tail, dotted] \arrow[u, "{\color{blue}\sim}", tail, dotted] \arrow[r, "{\color{blue}\sim}", dotted, two heads] & \bullet \arrow[d, dotted, two heads] \arrow[u, "\sim"', tail, dotted] \\
\bullet                                              & \bullet \arrow[l, tail, dotted] \arrow[r, "\sim", dotted, two heads]                                                              & \bullet                                                              
\end{tikzcd}                                                             
\end{equation}
in $X$, including the dotted arrows, whose top right and bottom left squares are pullbacks (so that the arrows with the blue labels are equivalences as well). Such a diagram is uniquely determined by its solid part, i.e.\ the left column and top row. Indeed, we can first fill the top right square with equivalences in the vertical direction. Then we can uniquely factor the composite in the bottom left square as an egressive followed by an ingressive edge (the resulting square is cartesian because $X$ is an orthogonal adequate triple). Finally, we fill the right two squares with equivalences as indicated.

 We can therefore conclude that the triple $(\Span(X), X_\ing, (X^\eg)^\op)$ is both adequate and orthogonal as soon as we can show that any
diagram as above defines a pushout square in $\Span(X^\rev)$. This is, however, rather unpleasant to do directly. We shall instead use the fact that pushout squares in an \icat{} $\C$ are precisely the cocartesian edges of the source map $s \colon \Ar(\C) \rightarrow \C$.
Here we regard the diagram ($\ast$) as a morphism in $\Ar(\Span(X^{\rev}))$ from the left vertical span to the right vertical span, so that the top horizontal span is the image of this morphism under $s$.

By Corollary \ref{FunofSpan} we can identify the functor $s \colon \Ar(\Span(X^\rev)) \rightarrow \Span(X^\rev)$ with 
\[\Span(\Q_1X^\rev) \xrightarrow{\Span(\ev_{(0 \leq 0)})} \Span(X^\rev).\] 
Put into this form, we can apply Barwick's criterion \ref{barwick_cocart} for cocartesian edges, or more precisely Corollary \ref{barwick_cocart_equiv}. It tells us that the diagram ($\ast$) defines a $\Span(\ev_{(0 \leq 0)})$-cocartesian edge if its left pointing half defines an $(\ev_{(0\leq 0)})^\eg$-cartesian edge in $\Q_1X^{\rev}$.

To see that it does, recall from \Cref{lem:adtrip_internal_hom}
that an edge in $\Q_1X^\rev = \Fun_{\AdTrip}(\TwR([1]),X^\rev)$ is
egressive if it is pointwise ingressive and its ambigressive square is
cartesian. Consider thus the (solid) lifting problem (ignoring the red
arrow for a moment)
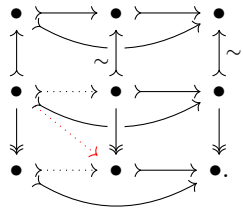
\begin{center}
\begin{tikzcd}
\bullet \arrow[r, tail] \arrow[rr, bend right, tail]                                                                  & \bullet \arrow[r, tail]                                   & \bullet                                    \\
\bullet \ar[rd, dotted, red] \arrow[d, two heads, "f"'] \arrow[u, tail] \arrow[rr, bend right, tail] \arrow[r, dotted, tail] & \bullet \arrow[d, crossing over,two heads] \arrow[u, crossing over, tail, "\sim"{pos=0.3}] \arrow[r, tail] & \bullet \arrow[d, two heads] \arrow[u, tail, "\sim"'] \\
\bullet \arrow[rr, bend right, tail, "h"] \arrow[r, dotted, tail]                                      & \bullet \arrow[r, tail]                                   & \bullet.                                   
\end{tikzcd}
\end{center} 
whose right hand column is the left half of ($\ast$), and whose lower bent rectangle is a pullback. We have to show that it admits a unique dotted filling (whose lower square is then automatically a pullback by pasting).

Considering the top half of the diagram, there is an essentially unique choice for the upper dotted arrow, because the middle upwards pointing map is an equivalence, and the composite of ingressives it is itself ingressive. By composition we also obtain the red arrow.

The bottom row of the diagram together with the middle left dot, then fit into a diagram
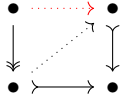
\begin{center}
\begin{tikzcd}
\bullet \arrow[d, two heads, "f"'] \arrow[r, dotted, red]            & \bullet \arrow[d, tail] \\
\bullet \arrow[ru, dotted] \arrow[r, tail, "h"] & \bullet                   
\end{tikzcd}
\end{center} 
Because the egressives and ingressives in $X$ determine an orthogonal factorisation system, this lifting problem has a unique filler. Moreover, because the right class of an orthogonal factorisation system satisfies right cancellation \cite[Proposition 5.2.8.6 (3)]{HTT}, this filler is necessarily ingressive. This provides the necessary lift, and concludes the argument. 
\end{proof}

\begin{definition}
We shall refer to $(\Span(X), X_\ing, (X^\eg)^\op)$ as the \emph{dual} of an orthogonal adequate triple $(X,X_\ing,X^\eg)$, and denote by $\ADual \colon \AdTrip^\perp \rightarrow \AdTrip^\perp$ the resulting functor.
\end{definition}

\begin{remark}
Unravelling the proof above, one arrives at the following description of induced maps out of the ambigressive pushouts in $\Span(X^\rev)$ from the previous proof: Suppose we are given a commutative square in $\Span(X^\rev)$ represented by a diagram
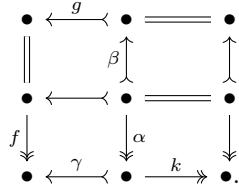
\begin{center}
\begin{tikzcd}
\bullet                                                      & \bullet \arrow[l, "g"', tail] \arrow[r, equal]                                                   & \bullet                                      \\
\bullet \arrow[u, equal] \arrow[d, "f"', two heads] & \bullet \arrow[l, tail] \arrow[d, "\alpha", two heads] \arrow[u, "\beta", tail] \arrow[r, equal] & \bullet \arrow[u, tail] \arrow[d, two heads] \\
\bullet                                                      & \bullet \arrow[l, "\gamma"', tail] \arrow[r, "k", two heads]                                                         & \bullet.                                     
\end{tikzcd}
\end{center}
whose lower left and upper right squares are cartesian. Factor the composite $fg$ 
into an egressive edge $f'$ followed by an ingressive edge $g'$. From the proof above, we learn that the pushout of the two spans starting in the upper left corner is given by the source of $g'$ (or equivalently the target of $f'$). The induced map out of it unwinds to the span
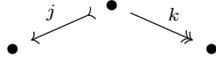
\[\begin{tikzcd}[row sep=0.4pc] & \bullet \ar[ld, "j"', tail] \ar[rd, "k", two heads] & \\ 
\bullet && \bullet
\end{tikzcd}\]
where $j$ is the solution of the lifting problem
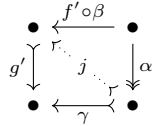
\begin{center}
\begin{tikzcd}
\bullet \arrow[d, "g' "', tail] & \bullet \arrow[d, "\alpha", two heads] \arrow[l, "f'\circ\beta"']   \\
\bullet                            & \bullet \arrow[l, "\gamma", tail] \arrow[lu, "j" description, dotted, tail]
\end{tikzcd}
\end{center}
As a consistency check, note that $j$ is indeed ingressive, as both $\gamma$ and $g'$ are.
\end{remark}

The remainder of this section is dedicated to proving the following result:

\begin{theorem}\label{D_is_equiv}
The functor $\ADual$ is an auto-equivalence of $\AdTrip^\perp$, which is its own inverse.
\end{theorem}

Below in \cref{actionisC2} we will upgrade this equivalence to a $\mathrm{C_{2}}$-action. In the large our strategy of proof is the same as that of Barwick, Glasman and Nardin in \cite{BGN}, who essentially treat the special case corresponding to factorisation systems given by (co)cartesian edges and morphisms lying in a single fibre. 

We have to construct a natural equivalence between $\AADual=\ADual\circ \ADual$ and the identity of $\AdTrip^\perp$. To do this, let us start by describing the composite $\AADual$ more explicitly.
\begin{observation}
Let us write $\core \Q^{\sharp}\colon \ssS\rt \ssS$ for the functor sending a simplicial \igpd{} $S$ to the simplicial \igpd{} given by 
$$
\core \Q_n^{\sharp}(S)=\Map_{\ssS}\big(\nerve\TwR([n]), S\big).
$$ 
By the definition of span \icats{} (see \Cref{Qsegal}), there is a natural transformation of functors $\AdTrip\rt \ssS$
$$\begin{tikzcd}[column sep=1.2pc]
\nerve\big(\Span(X)\big)\arrow[r, hookrightarrow] & \core \Q^{\sharp}(\nerve(X)).
\end{tikzcd}$$
In every simplicial degree, this is given by the inclusion of path components 
\[\Map_{\AdTrip}(\TwR([n]),X) \subseteq \Map_{\Cat}(\TwR([n]),X).\]
Applying this reasoning twice, one sees that for $X$ in $\AdTrip^\perp$, there is a natural map of simplicial \igpds{}
$$\begin{tikzcd}[column sep=1.2pc]
\nerve\big(\AADual(X)\big)\arrow[r, hookrightarrow] & 
\core \Q^{\sharp}\big(\nerve(\ADual(X))\big)
\arrow[r, hookrightarrow] & \core \Q^{\sharp}\big(\core \Q^{\sharp}(\nerve X)\big)
\end{tikzcd}$$
which is an inclusion of path components in each degree.
To unravel the target of the above map, let us denote by 
$$
\TTw(A)=\TwR(\TwR(A))
$$
the twofold iterated twisted arrow \icat{} of an $\infty$-category $A$. We then obtain a natural equivalence
\begin{align*}
\core \Q_n^{\sharp}\big(\core \Q^{\sharp}(\nerve X)\big)&\simeq \Map_{\ssS}\Big(\nerve\TwR([n]), \core \Q^{\sharp}(\nerve X)\Big)\\
&\simeq \lim_{[m]\in \Del/\nerve(\TwR[n])} \Map_{\ssS}\big(\nerve\TwR([m]), \nerve X\big)\\
&\simeq \Map_{\ssS}\big(\nerve\big(\TTw([n])\big), \nerve X\big) \\
&\simeq \Map_{\Cat}(\TTw([n]), X).
\end{align*}
Here the third line uses that $\nerve\TwR(-)$ preserves those colimits of $\infty$-categories that are preserved by the nerve functor (see the proof of Proposition \ref{lem:Q_preserves_limits}).
Summarising, we see that there is a map of simplicial \igpds{}, depending functorially on $X\in \AdTrip^\perp$, which is a degreewise inclusion of path components
\begin{equation}\label{eq:inclusion in double spans}
\begin{tikzcd}
\nerve\big(\ADual\circ\ADual(X)\big)\arrow[r, hookrightarrow] & \Map_{\Cat}\big(\TTw(-), X\big).
\end{tikzcd} 
\end{equation}
\end{observation}
To identify the essential image of \eqref{eq:inclusion in double spans}, let us make the following construction:
\begin{construction}\label{cons:tw2}
Note that $\TTw([n])$ is equivalent to the poset whose objects are tuples $abcd$ with $0\leq a\leq b\leq c\leq d\leq n$, corresponding to a map $(a\leq d)\rt (b\leq c)$ in $\TwR([n])$. The partial order is then given by $abcd\leq a'b'c'd'$ when $a \leq a' \leq b' \leq b \leq c \leq c' \leq d' \leq d$.
We define the following four wide subcategories of $\TTw([n])$:
\begin{enumerate}
    \item\label{it:tw1} $\TTw([n])_1$ is the subcategory spanned by the edges $abcd \rightarrow a'bcd$.
    \item\label{it:tw2} $\TTw([n])_2$ is the subcategory spanned by the edges $ab'cd \rightarrow abcd$.
    \item\label{it:tw3}  $\TTw([n])_3$ is the subcategory spanned by the edges $abcd \rightarrow abc'd$.
    \item\label{it:tw4}  $\TTw([n])_4$ is the subcategory spanned by the edges $abcd' \rightarrow abcd$.
\end{enumerate}
\end{construction}
\begin{proposition}\label{prop:identifying double dual}
Let $X\in \AdTrip^\perp$ and let $[n]\in \Del$. Then the natural
transformation \eqref{eq:inclusion in double spans} identifies the
domain with those path components in $\Map_{\Cat}(\TTw([n]), X)$
consisting of maps $f\colon \TTw([n])\rt X$ that restrict to
$$\begin{tikzcd}[row sep=0pc]
\TTw([n])_1\arrow[r, "f"] & X_{\ing} & \TTw([n])_2\arrow[r,  "f"] & \core X\\
\TTw([n])_3\arrow[r,  "f"] & X^{\eg} & \TTw([n])_4\arrow[r, "f"] & \core X.
\end{tikzcd}$$
\end{proposition}

\begin{proof}

A functor 
$$
f\colon [n]\rt \ADual\circ \ADual(X)$$
corresponds by adjunction to a map $\TwR([n])\rt \ADual(X)$ of adequate triples, i.e.\ a map of $\infty$-categories $f'\colon \TwR([n])\rt \ADual(X)$ such that 
$$
f'(\TwR([n])_\ing)\subseteq \ADual(X)_{\ing}\qquad \text{and}\qquad f'(\TwR([n])^{\eg})\subseteq \ADual(X)^{\eg}.
$$
Here we importantly use that every ambigressive square in the target is automatically a pullback, otherwise the functor $f'$ would have to furthermore preserve ambigressive pullbacks. By \Cref{modelfuncob}, the map of $\infty$-categories underlying $f'$ corresponds itself to a map $\TwR(\TwR([n]))\rt X$ of adequate triples, i.e.\ a map $f''\colon\TwR(\TwR([n]))\rt X$ such that
$$
f''\Big(\TwR\big(\TwR([n])\big)_\ing\Big)\subseteq X_{\ing}\qquad \text{and} \qquad f''\Big(\TwR\big(\TwR([n])\big)^{\eg}\Big)\subseteq X^{\eg}.
$$
Let us now unravel these conditions using the description of $\TTw([n])$ from Construction \ref{cons:tw2}. Note that $\TwR(\TwR([n]))_\ing$ consists of maps of tuples of the form $abcd\leq a'bcd'$, while $\TwR(\TwR([n]))^{\eg}$ consists of maps $abcd\leq ab'c'd$. On the other hand, $\TwR(\TwR([n])_{\ing})$ consists of maps of the form $abcc\leq a'b'cc$ and $\TwR(\TwR([n])^{\eg})$ consists of maps $aacd\leq aac'd'$.

Furthermore, $\ADual(X)_{\ing}$ is the subcategory consisting of spans of the form $\begin{tikzcd}[cramped, column sep=1.5pc]\cdot & \cdot\arrow[r, tail]\arrow[l, "\sim"{above}] & \cdot\end{tikzcd}$, where the right map is in $X_\ing$. Similarly, $\ADual(X)^{\eg}$ is the subcategory whose morphisms are spans of the form $\begin{tikzcd}[cramped, column sep=1.5pc]\cdot & \cdot\arrow[l, two heads]\arrow[r, "\sim"{above}] & \cdot\end{tikzcd}$, where the left map is in $X^\eg$. Combining all this, the above conditions translate as follows:
\begin{enumerate}
\item Every $abcd\leq a'bcd'$ is mapped to $X_\ing$.
\item Every $abcd\leq ab'c'd$ is mapped to $X^\eg$.
\item Every $abcc\leq a'bcc$ is sent to $X_\ing$ and every $abcc\leq ab'cc$ is sent to an equivalence.
\item Every $aacd\leq aac'd$ is sent to $X^\eg$ and every $aacd\leq aacd'$ is sent to an equivalence.
\end{enumerate}
Note that these conditions are certainly implied by the ones from the proposition. Conversely, given these conditions, the ones from the statement follow: for example, the map $abcd\leq abcd'$ fits into a square
$$\begin{tikzcd}
aacd\arrow[r]\arrow[d] & aacd'\arrow[d]\\
abcd \arrow[r] & abcd'
\end{tikzcd}$$
The top horizontal arrow is mapped to an equivalence in $X$ by (4), the vertical arrows are sent to maps over $X^\eg$ by (2). In particular, the bottom arrow maps to $X^\eg$, but also to $X_\ing$ by (1). An edge which is in both the left and right part of an orthogonal factorisation system is an equivalence, and therefore $abcd\leq abcd'$ maps to an equivalence.
\end{proof}
Let us write $\cat{LTw}^{(2)}([n])$ for the localisation of $\Tw^{(2)}([n])$ at the classes of maps \ref{it:tw2} and \ref{it:tw4} and let $\cat{LTw}^{(2)}([n])_1$ and $\cat{LTw}^{(2)}([n])_3$ be its wide subcategories generated by the images of the maps \ref{it:tw1} and \ref{it:tw3} from \Cref{cons:tw2}. Since each map in $\Del$ induces a map $\Tw^{(2)}([m])\rt \Tw^{(2)}([n])$ preserving these four classes of maps, it follows that the triple of $\infty$-categories 
\begin{equation}\label{eq:localised triple}
\big(\cat{LTw}^{(2)}([n]), \cat{LTw}^{(2)}([n])_1, \cat{LTw}^{(2)}([n])_3\big) \quad \in \Fun\big(\Lambda_2[2], \Cat)
\end{equation}
depends functorially on $[n]$. \Cref{prop:identifying double dual} implies that the \igpd{} $\nerve(\AADual(X))_n$ is naturally equivalent to the \igpd{} of maps in $\Fun(\Lambda_2[2], \Cat)$ 
$$
\Big(\cat{LTw}^{(2)}([n]), \cat{LTw}^{(2)}([n])_1, \cat{LTw}^{(2)}([n])_3\Big)\rt \big(X, X_{\ing}, X^{\eg}\big).
$$
We will now identify the localised triples \eqref{eq:localised triple} more explicitly. To do this, consider the composition
$$
[k]\times[1]\rightarrow [k]\star [k] \rightarrow [k]\star [k]^{\op}\star [k]\star [k]^{\op},$$
where the first map is the unique natural transformation between the two inclusions $[k]\rightarrow [k]\star [k]$ of the summands into the join. This induces a natural transformation
$$
z\colon \TTw(-)\rt \Ar(-)
$$
sending a tuple $abcd$ in $\TTw([n])$ to $ac$ in $\Ar([n])$.  Recall from \Cref{ex:arrow triple} that the arrow category $\Ar([n])$ comes equipped with the structure of an adequate triple in which a map is egressive (ingressive) if its image under the source (target) map is an equivalence. The map $z$ then sends the maps \ref{it:tw2} and \ref{it:tw4} to equivalences, the maps \ref{it:tw1} to ingressive arrows and the maps \ref{it:tw3} to egressive arrows. Consequently, we obtain a natural transformation of cosimplicial diagrams in $\Fun(\Lambda_2[2], \Cat)$
\begin{equation}\label{eq:loc tw2}
\Big(\cat{LTw}^{(2)}(-), \cat{LTw}^{(2)}(-)_1, \cat{LTw}^{(2)}(-)_3\Big)\rt \Ar(-).
\end{equation}
Restriction along this diagram of triples then induces a natural map of simplicial \igpds{}
\[
\zeta\colon \nerve(X)\simeq \Map_{\AdTrip^{\perp}}\big(\Ar(-), X\big)\rt \Map_{\Fun(\Lambda_2[2], \Cat)}\big(\cat{LTw}^{(2)}(-), X\big)\simeq\nerve \AADual(X).
\]
for every $X\in \AdTrip^\perp$. Here the first equivalence follows from the fact that $\Ar([n])$ is the free orthogonal adequate triple generated by $[n]$, by \Cref{prop:free triple}.
\begin{lemma}\label{lem:zeta}
For every $X\in \AdTrip^\perp$, the natural transformation $\zeta$ is an equivalence of complete Segal \igpds{}.
\end{lemma}
\begin{proof}
Since the domain and the target of $\zeta$ are complete Segal \igpds{}, $\zeta$ is an equivalence of simplicial \igpds{} as soon as it induces an equivalence in simplicial degrees $0$ and $1$. It therefore suffices to show that the natural transformation \eqref{eq:loc tw2} is an equivalence of triples of \icats{} for $n=0$ or $n=1$.

For $n=0$ this is evident since the domain and codomain of \eqref{eq:loc tw2} are both just a point. For $n=1$, at the level of the underlying \icats{}, we have to verify that the functor $z\colon \TTw([1])\rt \Ar([1])$ exhibits $\Ar([1])$ as the localisation of $\TTw([1])$ at the maps in $\TTw([n])_2$ and $\TTw([n])_4$. The map $z$ is given by the map
$$\begin{tikzcd}
\big(0000\lt 0001\rt 0011\lt 0111\rt 1111\big)\arrow[r] & \big(00\rt 01\rt 11\big)
\end{tikzcd}$$
collapsing the two left-pointing arrows. This is manifestly a localisation. Furthermore, under the map $z$, the noninvertible arrow in $\TwR([1])_1$, i.e.\ $0111\to 1111$, is sent to the ingressive $01\to 11$ and the noninvertible arrow in $\TwR([1])_3$, i.e.\ $0001\to 0011$, is sent to the egressive $00\to 01$. We conclude that the map $\cat{LTw}^{(2)}([1])\rt \Ar([1])$ is an equivalence of (adequate) triples, as desired.
\end{proof} 

\begin{proof}[Proof of Theorem \ref{D_is_equiv}]
Lemma \ref{lem:zeta} shows that after post-composing with the forgetful functor $U\colon \AdTrip\rightarrow \Cat$, there is a natural equivalence
$$\begin{tikzcd}
\zeta\colon U(X)\arrow[r, "\sim"] & U(\AADual(X)).
\end{tikzcd}$$
It remains to verify that this equivalence on the underlying \icats{} also identifies the subcategories of ingressives and egressives. But this is a direct consequence of the naturality of $\zeta$ in $X\in \AdTrip^{\perp}$, since the inclusions $\AADual(X)_{\ing}\hookrightarrow \AADual(X)\hookleftarrow \AADual(X)^{\eg}$ arise as the functors underlying the maps of adequate triples
$$
\AADual(X_\ing,X_\ing,\iota X)\rt \AADual(X)\lt \AADual(X^\eg,\iota X,X^\eg)
$$
and the same for $X_{\ing}\hookrightarrow X\hookleftarrow X^{\eg}$. 
\end{proof}

\section{Cartesian fibrations between orthogonal adequate triples}\label{sec:fibs}
The purpose of this section is to analyse the interaction between
orthogonal adequate triples and cartesian fibrations. In particular,
we show that orthogonal adequate triples $(X,X_\ing,X^\eg)$ uniquely
correspond to cartesian fibrations with contractible fibres, by taking
$X$ to $X \rightarrow X[(X^\eg)^{-1}]$, generalising results of Lanari
\cite{LanariCart}. More generally, we study various kinds of
fibrations between adequate triples that are preserved (or exchanged)
by the dualisation equivalence from \Cref{D_is_equiv}. When the base $X$ arises from a product of two \icats{} (\Cref{ex:product}), we can identify these fibrations with the two-variable fibrations studied extensively in \cite{PartI}; this specialisation will occupy the next section.

The main notion this section will be concerned with is:

\begin{definition}
We will say that a map $p \colon Y \rightarrow X$ of orthogonal adequate triples is an \emph{ingressive cartesian fibration} if $Y$ admits all $p$-cartesian lifts over $X_\ing$ and these precisely make up $Y_\ing$. Given an orthogonal adequate triple $X$, we will write $\Cart_{\ing}(X)\subseteq \AdTrip^{\perp}/X$ for the full subcategory on the ingressive cartesian fibrations. 
\end{definition}

\begin{observation}\label{obs:very explicit} \label{obs:fibres of egr cart}
Let us make the following remarks:
\begin{enumerate} 
\item Let $p\colon Y\rt X$ be an ingressive cartesian fibration. By definition, $Y$ has all $p$-cartesian lifts over $X_{\ing}$ and $Y_{\ing}$ is the wide subcategory on the $p$-cartesian lifts of ingressive morphisms. Since $Y_{\ing}$ and $Y^{\eg}$ form a factorisation system, this determines $Y^{\eg}$ uniquely and we find that $Y$ arises exactly from the construction in \Cref{lem:carttriport}; in particular $Y^\eg = p^{-1}(X^\eg)$.
\item For a map of orthogonal adequate triples $p \colon Y \rightarrow
  X$, an ingressive edge $g$ in $Y$ is $p$-cartesian if and only if it
  is $p_\ing$-cartesian, since lifting an arbitrary map $f \colon x
  \rightarrow x'$ in an orthogonal adequate triple along an ingressive
  $g \colon z \rightarrow x'$ is equivalent to lifting $f' \colon w
  \rightarrow x'$ along $g$ where $f = f'h$ is the factorisation of
  $f$ into an egressive followed by an ingressive.
\item Ingressive cartesian fibrations are stable under pullbacks. In
  particular, if $p\colon Y\rt X$ is an ingressive cartesian fibration
  and $x\in X$, then the fibre over $x$ carries the structure of an
  orthogonal adequate triple $(Y_x, \core Y_x, Y_x)$.
\end{enumerate}
\end{observation}

We start with the following partial converse to \Cref{lem:carttriport}:
\begin{proposition}\label{prop:orth_ad_trip_all_cart}
Let $(X, X_{\ing}, X^{\eg})$ be an orthogonal adequate triple. Then the map 
\[p \colon (X,X_\ing,X^\eg) \longrightarrow (X[(X^{\eg})^{-1}],X[(X^{\eg})^{-1}],\iota(X[(X^{\eg})^{-1}]))\] 
is an ingressive cartesian fibration. In other words, $p \colon X \rightarrow X[(X^\eg)^{-1}]$ is a cartesian fibration, inverts only the egressives in $X$, and the ingressive maps in $X$ are precisely the $p$-cartesian edges.
\end{proposition}

Before we dive into the proof, we spell out some consequences. To this end, let us write 
$\Cart$ for the subcategory of $\Ar(\Cat)$ spanned by objects the cartesian fibrations and by morphisms those squares
	\[
	\begin{tikzcd}
		X \arrow[d, "p"{swap}] \arrow[r, "f"] & Y \arrow[d,"q"] \\
		S \arrow[r]           & T          
	\end{tikzcd}
	\] such that $f$ sends $p$-cartesian edges to $q$-cartesian edges.
By Proposition \ref{prop:orth_ad_trip_all_cart}, the assignment $X\longmapsto (X\rightarrow X[(X^{\eg})^{-1}])$ extends to a functor 
$$
\mathcal{L}^{\eg}\colon \AdTrip^\perp\rt \Cart,
$$
which is fully faithful by the universal property of a localisation. We claim that the assignment of the orthogonal adequate triple $(Y,Y^\dagger, p^{-1}(\iota S))$ to a cartesian fibration $p\colon Y\rightarrow S$ defines a right adjoint to $\mathcal{L}^\eg$: Both 
\[
\Map_{\Cart}\big(X \rightarrow X[(X^\eg)^{-1}],p\big) \quad \text{and} \quad \Map_{\AdTrip}\big(X,(Y,Y^\dagger,p^{-1}(\iota S))\big)
\]
are given by those functors $X \rightarrow Y$ that take ingressives and egressives to $p$-cartesian and fibrewise maps, respectively. 

\begin{proposition}\label{prop:orthogonal=cart fib}
The functor $\mathcal{L}^{\eg}$ restricts to an equivalence between
the category $\AdTrip^\perp$ and the full subcategory of $\Cart$
spanned by those cartesian fibrations whose fibres all have
contractible realisations.
\end{proposition}

Given the discussion above, this follows immediately from the equivalence between the first and third items in the following:

\begin{lemma}\label{lem:cartlocalise}
For a cartesian fibration $p \colon A \rightarrow B$ the following
conditions are equivalent:
\begin{enumerate}
\item $p$ is a localisation,
\item $p$ is cofinal, and
\item the realisations of the fibres of $p$ are contractible.
\end{enumerate}
\end{lemma}

\begin{proof}
It is generally true that localisations are cofinal (and coinitial): If
$p$ is a localisation it follows straight from the definition of Kan extension as an
adjoint to restriction that any functor $F \colon B \rightarrow C$ is right (and left) Kan extended from its restriction along $p$, and thus has the same colimit as $Fp$. This proves $(1) \Rightarrow (2)$. 

To see that $(2) \Leftrightarrow (3)$, we observe that the inclusion 
\[p^{-1}(b) \subseteq b/p, \quad a \longrightarrow (a,\id_b),\]
where $b/p$ denotes the pullback $b/B \times_{B} A$ along $p$,
admits a right adjoint, given by taking $(a,f \colon b \rightarrow p(a))$ to a $p$-cartesian lift of $f$ ending at $a$. Thus $|p^{-1}(b)| \simeq |b/p|$, whence the claim follows from Joyal's cofinality criterion \cite[Theorem 4.1.3.1]{HTT}.

To finally see that $(3) \Rightarrow (1)$, consider for some $\infty$-category $C$ the functor
\[p^* \colon \Fun(B,C) \longrightarrow \Fun^{\mathrm{w}}(A,C),\]
where the superscript on the right denotes those functors inverting all maps that $p$ inverts. We have to show that $p^*$ is an equivalence. We first claim that it has a right adjoint $p_*$. By \cite[4.3.3.7]{HTT} this will follow from
\[b/p \xrightarrow{\mathrm{fgt}} A \xrightarrow{F} C \]
admitting a limit for every $b \in B$ and $F \colon A \rightarrow C$ that inverts the fibrewise maps. But as we just discussed above, the inclusion $p^{-1}(b) \subseteq b/p$ admits a right adjoint and is thus coinitial, so we may instead consider
\[p^{-1}(b) \subseteq A \xrightarrow{F} C, \]
which factors through the localisation $p^{-1}(b) \rightarrow
|p^{-1}(b)|$ since $F$ inverts all fibrewise maps by assumption. Thus
the assumption $|p^{-1}(b)| \simeq *$ implies that $F(a)$ is a limit
of the above diagram for any $a \in p^{-1}(b)$. It also follows
  from the formula for the adjoint as a right Kan extension that the
  unit and counit of the adjunction $(p^*,p_*)$ are equivalences, as desired.
\end{proof}

\begin{example}\label{rmk:pointedcart}
A cartesian fibration $p\colon X\rightarrow S$ is called \emph{pointed} if both $X$ and $S$ have a terminal object and $p$ preserves it. In \cite{LanariCart} Lanari shows that the $\infty$-category of pointed cartesian fibrations is equivalent to a certain subcategory of factorisation systems (he calls them cartesian). 

To see the relationship between this equivalence and ours, note that every fibre $X_s$ of a pointed cartesian fibration admits a final object, given by the source $x$ of a cartesian edge $x \rightarrow \ast$, which lies over the unique map $s \rightarrow \ast$ in $S$. Therefore the fibre of pointed cartesian fibrations have contractible realisations, and we conclude that our equivalence extends that of \cite{LanariCart}.

In fact, in the case of pointed cartesian fibrations Lemma \ref{lem:cartlocalise} can be sharpened: It is easy to check that assigning to some $s \in S$ the source $x$ of a cartesian lift $x \rightarrow \ast$ of $s \rightarrow \ast$ defines a right adjoint to $p$, so that any pointed cartesian fibration is in fact a left Bousfield localisation. As a typical application, one finds that an exact functor between stable $\infty$-categories is a cartesian fibration if and only if it is a right split Verdier projection (in the language of \cite[Appendix A]{CDH2}), and both a cartesian and a cocartesian fibration if and only if it participates (as $j^*$) in a stable recollement as in Example \ref{exrecoll}.
\end{example}

\begin{corollary}\label{actionisC2}
The self-equivalence $\ADual$ upgrades to a $\mathrm{C}_2$-action on $\AdTrip^{\perp}$.
\end{corollary}

\begin{proof}
According to \Cref{BGNidentification} above, the embedding $\mathcal{L}^{\eg} \colon \AdTrip^\perp \rightarrow \Cart$ from \Cref{prop:orthogonal=cart fib} translates $\ADual$ to the functor that takes a cartesian fibration to its fibrewise opposite (defined via straightening, postcomposing with $(-)^\op \colon \Cat \rightarrow \Cat$, and then unstraightening). Since $\Aut(\Cat) \simeq \mathrm{C}_2$ by To\"en's theorem \cite{Toen}, this latter operation defines a $\mathrm{C}_2$-action, hence so does the former. 
\end{proof}

It would be interesting to provide a more direct construction of the coherences for this $\mathrm{C}_2$-action that does not rely on the (un)straightening equivalence, but we have not pursued this. 

As another consequence, we obtain the following description of free cartesian fibrations, a generalisation of which was already proven by different means in \cite{GHN}:
\begin{corollary}\label{cor:free fib}
Consider the functor $\Cart\rt \Cat$ sending a cartesian fibration $p\colon Y\to X$ to its domain $Y$. This functor admits a left adjoint sending $A$ to the cartesian fibration $s\colon \Ar(A)\to A$.
\end{corollary}
\begin{proof}
The forgetful functor factors as $\Cart\rt \AdTrip^{\perp}\rt \Cat$. The first functor has left adjoint $\mathcal{L}^{\eg}$ by \Cref{prop:orthogonal=cart fib} and the second has left adjoint $\Ar$ by \Cref{prop:free triple}. The composite precisely sends $A$ to $s\colon \Ar(A)\rt A$, since this has contractible fibres and the associated orthogonal triple structure on $\Ar(A)$ is precisely that from \Cref{ex:arrow triple}.
\end{proof}

We now come to the proof of \cref{prop:orth_ad_trip_all_cart}.

\begin{proof}[Proof of \cref{prop:orth_ad_trip_all_cart}]
Let us start by recalling the model for $X[(X^{\eg})^{-1}]$ given by the relative Rezk nerve \cite{MazelGee}: for each $[n]$, let us write 
$$\begin{tikzcd}[column sep=1.5pc]
\mm{N}^\mm{rel}(X)(n)\arrow[r, hook] & \Fun([n], X)
\end{tikzcd}$$
for the subcategory of all functors $[n]\rt X$ and natural
transformations which are pointwise egressive. The core of this
\icat{} is simply $\Map_{\Cat}([n], X)$. Taking geometric realisations,
we then obtain a map of simplicial \igpds{}
$$\begin{tikzcd}
\mm{N}(X)=\Map_{\Cat}(-, X)\arrow[r] & \big|\mm{N}^\mm{rel}(X)\big|
\end{tikzcd}$$
such that the induced map on associated $\infty$-categories is the localisation $X\rt X[(X^{\eg})^{-1}]$.

We are now going to give a smaller description of $\big|\mm{N}^\mm{rel}(X)\big|$, which will show in particular that it already satisfies the Segal condition (so that the associated $\infty$-category is just its completion). To this end, consider the full subcategories 
$$
\mm{N}_{\ing}^\mm{rel}(X)(n)\subseteq \mm{N}^\mm{rel}(X)(n)
$$
whose objects are functors $[n]\rt X$ with values in $X_{\ing}$. For each $\alpha\in \mm{N}^\mm{rel}(X)(n)$, there exists an initial $\tilde{\alpha}$ in $\mm{N}_{\ing}^\mm{rel}(X)(n)$ equipped with a map from $\alpha$, given by the unique factorisation (starting at the end) 
$$\begin{tikzcd}
\alpha(0) \arrow[r] \arrow[d, two heads] & \alpha(1)\arrow[r] \arrow[d, two heads] & \dots\arrow[r] & \alpha(n-1)\arrow[r]\arrow[d, two heads] & \alpha(n)\arrow[d, equal]\\
\tilde{\alpha}(0)\arrow[r, tail] & \tilde{\alpha}(1)\arrow[r, tail] & \dots \arrow[r, tail] & \tilde{\alpha}(n-1)\arrow[r, tail] & \tilde{\alpha}(n).
\end{tikzcd}$$
It follows that each inclusion $\mm{N}_{\ing}^\mm{rel}(X)(n)\hookrightarrow \mm{N}^\mm{rel}(X)(n)$ admits a left adjoint.
We thus obtain a diagram of simplicial categories
\begin{equation}\label{diag:localisation}
\begin{tikzcd}
\mm{N}(X_{\ing})\arrow[r]\arrow[d] & \mm{N}_{\ing}^\mm{rel}(X)\arrow[d]\arrow[r] & \big|\mm{N}_{\ing}^\mm{rel}(X)\big|\arrow[d, "\sim"] \\
\mm{N}(X)\arrow[r] & \mm{N}^\mm{rel}(X)\arrow[r] & \big|\mm{N}^\mm{rel}(X)\big| \end{tikzcd}
\end{equation}
where the objects on the left and on the right are simplicial \igpds{} and the right vertical map is an equivalence, since it is given in each simplicial degree by the geometric realisation of a right adjoint functor.

The simplicial \icat{} $\mm{N}_{\ing}^\mm{rel}(X)$ satisfies the Segal conditions: The equivalence $\Fun([n], X)\simeq \Fun([1], X)\times_{X}\dots\times_X \Fun([1], X)$ restricts to an equivalence of \icats{}
\begin{equation}\label{eq:Segal}\begin{tikzcd}
\mm{N}_{\ing}^\mm{rel}(X)([n])\arrow[r] & \mm{N}_{\ing}^\mm{rel}(X)([1])\times_{\mm{N}_{\ing}^\mm{rel}(X)(0)} \dots \times_{\mm{N}_{\ing}^\mm{rel}(X)(0)} \mm{N}_{\ing}^\mm{rel}(X)([1])
\end{tikzcd}\end{equation}
Here the pullbacks are taken along the source and target functors. The target functor
\begin{equation}\label{eq:target kan}\begin{tikzcd}
t\colon \mm{N}_{\ing}^\mm{rel}(X)([1])\arrow[r] & \mm{N}^\mm{rel}_{\ing}(X)(\{1\})=X^{\eg}
\end{tikzcd}\end{equation}
is both a left fibration (since each
$\rightarrowtail \twoheadrightarrow$ fits into a unique ambigressive
square) and a right fibration (since ingressive maps can be pulled
back along egressive maps). This implies that the pullbacks along $t$ in \eqref{eq:Segal}
induce pullbacks after taking realisations (see again \cite{Steimle,Steinebrunner} or the proof of \cite[Theorem 2.5.1]{CDH2} for the more (co)cartesian version of this assertion), so that the equivalent simplicial \igpds{} on the right in \eqref{diag:localisation} also satisfy the Segal conditions.

We will use this description of $X[(X^{\eg})^{-1}]$ in terms of $\mm{N}^{\mm{rel}}_{\ing}(X)$ to identify its over-categories. To this end, recall that for any Segal \igpd{} $S$ and $s\in S(0)$, the simplicial \igpd{} 
\begin{equation}\label{eq:slice}
\big(S/s\big)(n) = S([n+1])\times_{S(\{n+1\})}\{s\}
\end{equation}
is again a Segal \igpd{}, such that the associated \icat{} $\asscat(S/s)\simeq \asscat(S)/s$ is a model for the over-category of $s$ \cite[Lemma 2.4.7]{CDH2}. 
Let us now take an object $x\in X$ and consider the maps of simplicial \igpds{} induced by the top row in \eqref{diag:localisation}
\begin{equation}\label{eq:slice simplicial}\begin{tikzcd}
\mm{N}(X_{\ing})/x\arrow[r] & \big|\mm{N}^{\mm{rel}}_{\ing}(X)/x\big|\arrow[r] & \big|\mm{N}^{\mm{rel}}_{\ing}(X)\big|/x.
\end{tikzcd}\end{equation}
Here the middle term is given by first taking \eqref{eq:slice} at the
level of simplicial \icats{}, and then taking realisations. Unraveling the definitions, the simplicial \icat{} $\mm{N}^{\mm{rel}}_{\ing}(X)/x$ is given in degree $n$ by the \icat{} with objects $\alpha(0)\rightarrowtail \dots \rightarrowtail \alpha(n)\rightarrowtail x$ and morphisms
$$\begin{tikzcd}
\alpha(0) \arrow[r, tail] \arrow[d, two heads] & \alpha(1)\arrow[r, tail] \arrow[d, two heads] & \dots\arrow[r, tail] & \alpha(n)\arrow[r, tail]\arrow[d, two heads] & x\arrow[d, equal]\\
\beta(0)\arrow[r, tail] & \beta(1)\arrow[r, tail] & \dots \arrow[r, tail] & \beta(n)\arrow[r, tail] & x.
\end{tikzcd}$$
Since the ingressives and egressives form a factorisation system, all
vertical maps are equivalences. Consequently, the first map in
\eqref{eq:slice simplicial} is an equivalence (even before taking
geometric realisations). In addition, since the target map
\eqref{eq:target kan} is a Kan fibration, so is the map $\mm{N}_{\ing}^\mm{rel}(X)([n+1])\rt \mm{N}_{\ing}^\mm{rel}(X)(\{n+1\})$ (being a composite of its base changes). It follows that the pullbacks \eqref{eq:slice} are preserved under taking classifying $\infty$-groupoids, so that the second map
in \eqref{eq:slice simplicial} is an equivalence. Since $\big|\mm{N}^{\mm{rel}}_{\ing}(X)\big|$ has associated \icat{} $X[(X^{\eg})^{-1}]$, the associated \icat{} of $\big|\mm{N}^{\mm{rel}}_{\ing}(X)\big|/x$ is $X[(X^{\eg})^{-1}]/p(x)$.

All in all, we have therefore found that for each object $x\in X$, the composite functor
$$\begin{tikzcd}
X_{\ing}/x\arrow[r, "\iota"] & X/x\arrow[r, "p"] & X[(X^{\eg})^{-1}]/p(x)
\end{tikzcd}$$
is an equivalence. Under this equivalence, the functor $p$ is identified with the functor $X/x\rt X_{\ing}/x$ sending each $y\to x$ to the ingressive part $y'\rightarrowtail x$ of its functorial egressive-ingressive factorisation. This functor admits a fully faithful right adjoint, whose essential image is given by the ingressive maps to $x$. By \cite[Lemma 2.16]{AyalaFrancis}, this implies that $p$ is a cartesian fibration and that an arrow in $X$ is $p$-cartesian if and only if it is ingressive.

In particular, we can apply \Cref{lem:carttriport} to obtain an orthogonal triple structure on $X$ whose ingressives are the $p$-cartesian arrows and whose egressives are the maps that are inverted by $p$. Since the first class coincides with $X_{\ing}$, the second class coincides with $X^{\eg}$.
\end{proof}

We now turn to more restrictive types of fibrations between orthogonal triples:

\begin{definition}\label{def:egcartfib}
An ingressive cartesian fibration $p\colon Y\rt X$ is said to be
\begin{enumerate}

\item an \emph{op-Gray fibration} if $p^{\eg}\colon Y^{\eg}\rt X^{\eg}$ is a cartesian fibration, and

\item a \emph{\lortho{}} if $p^{\eg}\colon Y^{\eg}\rt X^{\eg}$ is a cocartesian fibration.
\end{enumerate}
For $X$ an orthogonal adequate triple, we will write $\Ortholax(X), \Grayop(X)$ for the subcategories of $\Cart_\ing(X)$ spanned by the \lorthos{} and op-Gray fibrations, and maps between them preserving (in addition) $p^{\eg}$-(co)cartesian morphisms.
\end{definition}

\begin{remark}
To avoid confusion we would like to explicitly mention that while arrows with the decoration $\twoheadrightarrow$ denoted cartesian edges in \cite{PartI}, with our conventions it is exactly the ingressive edges in $Y$ which are all $p$-cartesian and that the egressive edges are \emph{almost never} all $p$-cartesian.
\end{remark}

The naming schema above is inspired by the definitions of \cite{PartI}. We will show in Section \ref{sec:product} that in the case $X= (A\times B,A\times \iota B,\iota A \times B)$, the definitions above reduce to the corresponding objects of \cite{PartI}. 

Note that $\Cart_\ing(A,A,\iota A)$ is simply the $\infty$-category of cartesian fibrations over $A$, as are the two subcategories $\Ortholax$ and $\Grayop$ in this case. On the other hand, $\Cart_{\ing}(A,\iota A,A)$ is just $\Cat/A$, while $\Grayop(A,\iota A,A) = \Cart(A)$ and $\Ortholax(A,\iota A,A) = \mathrm{Cocart}(A)$.

One can informally think about an op-Gray fibration over $X$ as
encoding a lax functor from $X$ to $\Cat$, which is strong on
$2$-simplices in $X_{\ing}$ or $X^{\eg}$, as well as $2$-simplices of
the form $x_0\twoheadrightarrow x_1 \rightarrowtail x_2$. This is substantiated (using the scaled straightening construction of Lurie \cite{LurieGoo}) by the following observation:
\begin{lemma}\label{lem:opgray}
An op-Gray fibration $p\colon Y\rt X$ is in particular a locally cartesian fibration. Furthermore, for each $2$-simplex $\sigma\colon [2]\rt X$ arising from the composition of an ingressive and an egressive map $x_0\twoheadrightarrow x_1\rightarrowtail x_2$, the restriction $\sigma^*(p)\colon \sigma^*Y\rt [2]$ is a cartesian fibration.
\end{lemma}
\begin{proof}
Note that the second part implies that $p$ is a locally cartesian fibration since every morphism in $X$ fits into a $2$-simplex $\sigma$ as indicated. To see that $\sigma^*(p)$ is a cartesian fibration, one simply notes that it admits enough cartesian lifts of $x_1\rightarrowtail x_2$ (since $p$ does) and locally cartesian lifts over $x_0\twoheadrightarrow x_1$. 
\end{proof}
\begin{lemma}\label{lem:cartortho}
Let $p\colon Y\rt X$ be an op-Gray fibration between orthogonal adequate triples. Then the following are equivalent:
\begin{enumerate}
\item\label{it:cart condition1} The underlying functor $p\colon Y\rt X$ is a cartesian fibration.

\item\label{it:cart condition2} For every ambigressive square
\begin{equation}\label{eq:ambigresive cart}\begin{tikzcd}
y_{00}\arrow[r, rightarrowtail, ]\arrow[d, twoheadrightarrow, "f"'] & y_{01}\arrow[d, twoheadrightarrow, "f'"]\\
y_{10}\arrow[r, rightarrowtail] & y_{11} 
\end{tikzcd}\end{equation}
in which $f'$ is $p^{\eg}$-cartesian, $f$ is $p^{\eg}$-cartesian as well.
\end{enumerate}
\end{lemma}
\begin{proof}
Assuming \ref{it:cart condition1}, the uniqueness of cartesian
  lifts implies that an egressive morphism in $Y$ is $p^{\eg}$-cartesian if and only if it is $p$-cartesian. Except for $f'$, all arrows in the square \eqref{eq:ambigresive cart} are therefore $p$-cartesian, which implies that $f'$ is $p$-cartesian as well.

For the converse, \Cref{lem:opgray} already asserts that $p$ is a locally cartesian fibration. It then suffices to verify that locally $p$-cartesian morphisms are closed under composition. To this end, let $f$ and $g$ be two composable locally $p$-cartesian morphisms in $Y$ and consider the diagram
$$\begin{tikzcd}[cramped]
\cdot\arrow[rd, "f"{swap}] \arrow[r, twoheadrightarrow, "f_1"] & \cdot \arrow[d, rightarrowtail, "f_2"]\arrow[r, twoheadrightarrow] & \cdot \arrow[d, rightarrowtail]\\
& \cdot \arrow[rd, "g"{swap}] \arrow[r, twoheadrightarrow, "g_1"] & \cdot\arrow[d, rightarrowtail, "g_2"] \\
& & \cdot
\end{tikzcd}$$
Here the two triangles factor $f$ and $g$ into an egressive followed by an ingressive morphism; \Cref{lem:opgray} shows that the maps $f_1, f_2, g_1$ and $g_2$ are all locally $p$-cartesian. The top right square is obtained by factoring the down-right composite into an egressive, followed by an ingressive map. Condition \ref{it:cart condition2} then implies that all individual arrows depicted in the above picture are locally $p$-cartesian (and the horizontal ones are $p$-cartesian). Since $p^{\eg}$ was a cartesian fibration, the top composite is then locally $p$-cartesian, and since the right vertical composite is $p$-cartesian, the composite $gf$ is locally $p$-cartesian as well. 
\end{proof}
\begin{definition}\label{def:ortho}
We will say that a map $p\colon Y\rt X$ of orthogonal adequate triples is a \emph{cartesian fibration} if it is an op-Gray fibration satisfying the equivalent conditions of \Cref{lem:cartortho}.
Dually, a \lortho{} is called an \emph{orthofibration} if for each ambigressive square \eqref{eq:ambigresive cart} in which $f$ is $p^{\eg}$-cocartesian, $f'$ is $p^{\eg}$-cocartesian as well.

Let us write $\Cart(X)\subseteq \Grayop(X)$ and $\ortho(X)\subseteq
\Ortholax(X)$ for the full subcategories on cartesian and
orthofibrations; this notation is justified by the following observation:
\end{definition}
\begin{observation}\label{obs:cart is cart}
Note that the datum of a cartesian fibration of orthogonal adequate triples $p\colon Y\rt X$ is equivalent to that of a cartesian fibration between the underlying \icats{}. In addition, a map in $\Cart(X)$ is required to preserve cartesian lifts over ingressive and egressive maps in $X$; since every map in $X$ factors into an egressive map followed by an ingressive one, it follows that the maps in $\Cart(X)$ preserve \emph{all} cartesian arrows. In other words, $\Cart(X)$ is simply equivalent to the \icat{} of cartesian fibrations over the underlying \icat{} of $X$.
\end{observation}
The main result about these various types of fibration is that they behave well under dualisation:
\begin{theorem}\label{thm:dualfibs}
Let $X$ be an orthogonal adequate triple. Then the natural equivalence
 \[\ADual\colon \AdTrip^{\perp}/X\rt \AdTrip^{\perp}/\ADual(X)\]
restricts to natural equivalences
\begin{align*}
\Cart_{\ing}(X)&\simeq \Cart_{\ing}(\ADual(X))\\
\Grayop(X)&\simeq \Ortholax(\ADual(X))\\
\Cart(X)& \simeq \ortho(\ADual(X)).
\end{align*}
Furthermore, for each fibration $p\colon Y\rt X$, the dual fibration $\ADual(Y)\rt \ADual(X)$ has fibres given by the \emph{opposites} of the fibres of $p$.
\end{theorem}
In the simple case $X = (A, A, \iota A)$ this result specialises to the statement that taking span categories for the adequate triple structure given by the cartesian and fibrewise maps provides an equivalence
\[\Span \colon \Cart(A) = \Cart_\ing(X) \rightarrow
  \Cart_\ing(\ADual(X))= \Cart(A)\] 
sending a cartesian fibration to the cartesian fibration with opposite fibres. Postcomposing with the equivalence $\Cart(A)\simeq \Cocart(A^{\op})$ sending $Y\rt A$ to $Y^{\op}\rt A^{\op}$, we obtain the equivalence first established by
Barwick, Glasman and Nardin \cite[Theorem 1.4]{BGN}. Note that by \Cref{span_op} this composite equivalence coincides with the equivalence obtained by first exchanging the ingressive and egressive maps and then applying the equivalence from \Cref{BGNidentification}.

In the case $X= (A,\iota A, A)$, on the other hand, the result simply specialises to the equivalence
\[\Cat/A = \Cart_\ing(X) \rightarrow \Cart_\ing(\ADual(X)) = \Cat/A^\op,\] given by taking opposites. In the case of a general $X$, our equivalence combines these two extremes, as we shall discuss in the case of two-variable fibrations from \cite{PartI} in the next section.
\begin{proof}
Note that naturality of these equivalences follows directly from $\ADual$ being a functor and taking over-categories being natural in the base.

Let us start with the first equivalence. Since $\ADual$ is its own inverse, it suffices to verify that for an ingressive cartesian fibration $p\colon Y\rt X$, the dual map $q = \ADual(p)\colon \ADual(Y)\rt \ADual(X)$ is an ingressive cartesian fibration. Recall that an ingressive map in $\ADual(Y)$ is given by a span $y_{11}\lto{\sim} y_{01} \rightarrowtail y_{00}$ in $Y$. 

To see that $q$ is an ingressive fibration, it suffices to verify that all such spans define $q$-cartesian edges in $\Span(Y)$, since there are then clearly enough cartesian lifts of edges in $X_{\ing}$. Equivalently, this means that the reverse span $y_{00}\leftarrowtail y_{01}\rto{\sim} y_{11}$ defines a $q^{\op}$-cocartesian arrow. We can identify the opposite map $q^{\op}\colon \Span(Y)^{\op}\rt \Span(X)^{\op}$ with $\Span(p^{\rev})\colon \Span(Y^{\rev})\rt \Span(X^{\rev})$, i.e.\ we exchange the roles of the ingressive and egressive maps. The reverse span is then $\Span(p^{\rev})$-cocartesian by \Cref{barwick_cocart_equiv}.

To see the second equivalence, we observe that given any $p\colon Y\rt X$ in $\Cart_{\ing}(X)$, the functor $q^{\eg}\colon \ADual(Y)^{\eg}\rt \ADual(X)^{\eg}$ is equivalent to the \emph{opposite} of $p^{\eg}\colon Y^{\eg}\rt X^{\eg}$. This immediately implies that $\ADual$ exchanges \lorthos{} and op-Gray fibrations, and furthermore preserves maps between them that preserve (co)cartesian ingressive arrows.

For the third equivalence, recall from \Cref{prop:dualadequatetrip} that an ambigressive square in $\ADual(Y)$ corresponds to a diagram of the form
$$
\begin{tikzcd}
\bullet                                              & \bullet \arrow[r, tail] \arrow[l, two heads, "\sim"]                                                                              & \bullet                                                               \\
\bullet \arrow[u, two heads, "\alpha"] \arrow[d, tail, "\sim"'] & \bullet \arrow[d, tail, "\sim"] \arrow[l, two heads, "\sim"] \arrow[u, two heads, "\beta"] \arrow[r, tail] & \bullet \arrow[d, tail, "\sim"] \arrow[u, two heads, "\gamma"'] \\
\bullet                                              & \bullet \arrow[l, two heads, "\sim"] \arrow[r, tail]                                                              & \bullet                                                              
\end{tikzcd}                                                             
$$
in $Y$. 
Suppose $p$ is a cartesian fibration. The left span defines a $q^{\eg}$-cocartesian arrow if and only if $\alpha$ is $p^{\eg}$-cartesian. This immediately implies that $\beta$ is also $p^{\eg}$-cartesian. Finally the top right square is cartesian, and therefore we conclude by \ref{lem:cartortho} that $\gamma$ is also $p^{\eg}$-cartesian, so that the right vertical span is $q^{\eg}$-cartesian. This implies by definition that $q$ is an orthofibration of orthogonal triples. A dual argument (again using \ref{lem:cartortho}) shows that the dual of an orthofibration is a cartesian fibration.

Finally, for the statement about fibres, note that $\ADual$ (being an equivalence) preserves fibres, i.e.\ pullbacks along a map (of orthogonal adequate triples) $\ast\rt X$. \Cref{obs:fibres of egr cart} then implies that the fibre of $\ADual(Y)\rt \ADual(X)$ over a point $x$ is given by $\ADual(Y_x, \core Y_x, Y_x)$, whose underlying category is equivalent to $Y_x^{\op}$ by \Cref{prop:spanswithcore}.
\end{proof}

\section{Dualisation and straightening of two-variable fibrations}\label{sec:product}
The purpose of this section is to use the results of the previous section to dualise and straighten various kinds of fibrations over a product of $\infty$-categories. The first part essentially consists of specialising \Cref{thm:dualfibs} to the case where $X$ is the particular orthogonal triple
$$
(A, B)^\perp \coloneqq (A\times B,A\times \core B, \core A\times B)
$$
associated to a product of two \icats{}. We will find that the various fibrations defined over $(A,B)^\perp$ in the previous section recover a subset of the fibrations which we considered in \cite{PartI}. 

In particular, we will extend the explicit description of dual (co)cartesian fibrations to the situation of curved ortho- and Gray fibrations. As an application we give an explicit description of parametrised adjoints from \cite{PartI}, extending previous work of Torii \cite{Torii}.

Let us start by making the various types of fibrations of orthogonal adequate triples $p\colon Y\rt (A, B)^\perp$ appearing in \Cref{sec:fibs} more explicit. Recall that for all these fibrations, the structure of an adequate triple on $Y$ is uniquely determined by the underlying functor (see \Cref{obs:very explicit})
$$
p=(p_1, p_2)\colon Y\rt A\times B.
$$
By Proposition \ref{lem:carttriport}, such $p$ determines an ingressive cartesian fibration if and only if all (ingressive) arrows in $A\times \core B$ have enough cartesian lifts. By \cite[Corollary 2.2.2]{PartI} or \cite[Proposition 2.4.1.3(3)]{HTT}, this is equivalent to 
$$\begin{tikzcd}
Y\arrow[rr, "p"]\arrow[rd, "p_1"{swap}] & & A\times B\arrow[ld, "\mm{pr}_1"]\\
& A
\end{tikzcd}$$
defining a map of cartesian fibrations over $A$ preserving cartesian arrows. In particular, each $\alpha\colon a\rt a'$ in $A$ gives rise to a \emph{cartesian transport functor} $\alpha^*\colon Y_{a'}\rt Y_a$ between the fibres of $p_1$. Following \cite{PartI}, let us write $\LCart(A, B)$ for the \icat{} of $p\colon Y\rt A\times B$ having $p$-cartesian lifts over $A\times \core B$ and maps preserving such $p$-cartesian lifts. In these terms, we have:

\begin{observation}\label{cor:lcart}
There are natural equivalences of \icats{}
$$\Cart_{\ing}((A, B)^\perp)\simeq \LCart(A, B).$$
\end{observation}
\begin{remark}\label{rem:lcocart}
There is an evident analogue $\LCocart(A, B)$ consisting of $p\colon Y\rt A\times B$ having $p$-cocartesian lifts over $A\times \core B$. This \icat{} does not have a good interpretation in terms of adequate triples; one can typically not equip $Y$ with the natural structure of an adequate triple.
\end{remark}

\begin{lemma}\label{lem:over AxB}
Let $p=(p_1, p_2)\colon Y\rt A\times B$ be a map in $\LCart(A, B)$. Then $p$ defines:
\begin{enumerate}
\item\label{it:opgray} an op-Gray fibration in the sense of \Cref{def:egcartfib} if and only if for each $a\in A$, the map $p_2\colon Y_a\rt Y\rt B$ from the fibre of $p_1$ is a cartesian fibration.

\item\label{it:lortho} a \lortho{} in the sense of \Cref{def:egcartfib} if and only if for each $a\in A$, the map $p_2\colon Y_a\rt Y\rt B$ from the fibre of $p_1$ is a cocartesian fibration.
\end{enumerate}
\end{lemma}

\begin{proof}
Observe that the map $Y^{\eg}\rt ((A, B)^{\perp})^{\eg}$ is simply given by the restriction of $p\colon Y\rt A\times B$ to $\core A\times B$. It follows from \cite[Remark 2.2.9]{PartI} that this restriction to $\core A\times B$ is a (co)cartesian fibration if and only if for each $a\in A$, the restriction $p_a\colon Y_a\rt \{a\}\times B$ of $p$ is a (co)cartesian fibration. Finally, note that $p_a$ can be identified with $p_2$ (restricted to the fibre $Y_a$).
\end{proof}

\begin{definition}\label{def:curv ortho}
Let us denote by $\Ortholax(A, B)$ the subcategory of $\Cat/A\times B$ whose objects are curved orthofibrations and whose morphisms are maps preserving locally $p$-cartesian morphisms over $A\times\core B$ and locally $p$-cocartesian morphisms over $\core A\times B$.
We define $\Grayop(A, B)$ similarly.
\end{definition}

We then immediately find:

\begin{corollary}\label{cor:gray=gray}
There are natural equivalences
$$\begin{tikzcd}[row sep=0.2pc]
\Ortholax((A, B)^\perp)\arrow[r, " \sim"] & \Ortholax(A, B) \\
\Grayop((A, B)^{\perp})\arrow[r, " \sim" ] & \Grayop(A, B).
\end{tikzcd}$$
\end{corollary}
\begin{remark}\label{rem:gray}
As one may expect, the notion of a \emph{Gray fibration} $p\colon Y\rt A\times B$ was introduced in \cite{PartI} as the opposite of an op-Gray fibration, i.e.\ there is an equivalence
$$\begin{tikzcd}
(-)^{\op}\colon \Gray(A, B)\arrow[r, "\sim"]  & \Grayop(A^{\op}, B^{\op}).\end{tikzcd}$$
In particular, $\Gray(A, B)$ is a certain subcategory of $\LCocart(A, B)$. As in \Cref{rem:lcocart}, it does not have a good analogue for general orthogonal adequate triples. 
\end{remark}

\begin{remark}
In \cite[Section 2]{PartI} we recorded various other ways to recognise (op-)Gray fibrations and \lorthos{}. For example, $(p_1, p_2)\colon Y\rt A\times B$ is a \lortho{} if $p_1$ is a cartesian fibration and $p_2$ is a cocartesian fibration \cite[Proposition 2.3.3]{PartI}, i.e.\
\[\localortho(A,B) = \LCart(A,B) \cap \RCocart(A,B).\] 
\end{remark}

\begin{lemma}\label{lem:ortho AxB}
A map $p=(p_1, p_2)\colon Y\rt A\times B$ is an orthofibration in the sense of \Cref{def:ortho} if and only if it is a \lortho{} and for each map $\alpha\colon a\rt a'$, the cartesian transport functor $\alpha^*\colon Y_{a'}\rt Y_{a}$ preserves $p_2$-cocartesian arrows.
\end{lemma}
\begin{proof}
Let us start by noting that an arrow in some fibre $Y_a\subseteq Y$ is $p_2$-cocartesian if and only if it defines a cocartesian arrow for the base change of $p$ to $\core A\times B$; see \cite[Corollary 2.2.7]{PartI}. This base change coincides with $p^{\eg}\colon Y^{\eg}\rt ((A\times B)^{\perp})^{\eg}$, so an arrow in $Y_a$ is $p_2$-cocartesian if and only if it defines a $p^{\eg}$-cocartesian arrow.

Using this, it follows that $\alpha^*$ preserves $p_2$-cocartesian arrows if and only if the following holds: for each square in $Y$
$$\begin{tikzcd}
y_{a, b}\arrow[r]\arrow[d] & y_{a', b}\arrow[d] \arrow[rr, mapsto, shorten=3ex, yshift=-1.5pc] & & {(a, b)}\arrow[r, "\alpha"]\arrow[d, "\beta"{swap}] & {(a', b)}\arrow[d, "\beta"]\\
y_{a, b'}\arrow[r] & y_{a', b'} & & {(a, b')}\arrow[r, "\alpha"] & {(a', b')}
\end{tikzcd}$$
in which the vertical arrows are $p$-cartesian and the bottom horizontal arrow is $p_{\ing}$-cocartesian, the top horizontal arrow is $p_2$-cocartesian. This is precisely the condition of \Cref{def:ortho}, since the ambigressive squares in $Y$ are precisely squares of the above form whose vertical maps are $p$-cartesian. 
\end{proof}
\Cref{lem:ortho AxB} asserts that the notion of an orthofibration as defined in \Cref{def:ortho} agrees (over $X=A\times B$) with the notion of an orthofibration employed in \cite{PartI} (see in particular \cite[Proposition 2.3.11]{PartI}). Writing $\ortho(A, B)\subseteq \Ortholax(A, B)$ for the full subcategory spanned by the orthofibrations, we therefore obtain:
\begin{corollary}\label{cor:ortho=ortho}
There are natural equivalences of \icats{} 
$$
\ortho((A, B)^\perp)\simeq \ortho(A, B).
$$
\end{corollary}
All in all, we have found that over the orthogonal triple $(A, B)^\perp=(A\times B, A\times \core B, \core A\times B)$, the \icats{} of fibrations introduced in \Cref{sec:fibs} coincide with those appearing under the same name in \cite{PartI}. 

We will now unravel the content of the dualisation equivalence from \Cref{thm:dualfibs} over the base $X=(A, B)^\perp$. To this end, let us start with the following observation:
\begin{lemma}\label{prop:spans of A x B}
	Let $A$ and $B$ be two $\infty$-categories. Then there is an equivalence of orthogonal adequate triples
	\[\ADual\big((A,B)^\perp\big) \simeq (A,B^{\op})^\perp.\]
\end{lemma}
\begin{proof}
	Note that $(A, B)^\perp$ decomposes as a product of the triples $(A,A,\core A)$ and $(B,\core B, B)$, both of which are evidently orthogonal and adequate (and both these properties are preserved under taking products). Since the span construction preserves products, being a right adjoint, the claim follows from \Cref{prop:spanswithcore}.
\end{proof}
Consequently, \Cref{thm:dualfibs} provides equivalences between certain types of fibrations over $A\times B$ and $A\times B^\op$, respectively. Recall from \Cref{thm:dualfibs} that these equivalences take opposite categories at the level of fibres. To conform with the conventions of \cite{PartI}, we will therefore compose the equivalences of \Cref{thm:dualfibs} with taking opposites, resulting in: 

\begin{definition}\label{d:sdual two-var}
For $\infty$-categories $A$ and $B$ we define functors
$$\begin{tikzcd}[row sep=0.2pc]
\SDualco\colon \LCart(A, B) \arrow[r,"\ADual"] & \LCart(A,B^\op) \arrow[r, "(-)^\op"] & \LCocart(A^{\op}, B) \\
\SDualcart \colon \LCocart(A^\op, B) \arrow[r,"(-)^\op"] & \LCart(A,B^\op) \arrow[r,"\ADual"] & \LCart(A,B)
\end{tikzcd}$$
using the identifications from \cref{cor:lcart} and \cref{prop:spans of A x B}.
\end{definition}

Using Theorem \ref{thm:dualfibs} and the identifications \cref{cor:gray=gray} from \cref{cor:ortho=ortho} we then immediately find:

\begin{corollary}\label{thm:dualisation functor}\label{Dual_eq} 
The functor $\SDualco \colon \LCart(A, B) \rightarrow \LCocart(A^{\op}, B)$ is an equivalence with inverse $\SDualcart$, and restricts to equivalences
\[\Ortholax(A, B) \overset{\sim}{\longrightarrow} \Gray(A^{\op}, B) \quad \text{and} \quad \ortho(A, B) \overset{\sim}{\longrightarrow} \Cocart(A^{\op}\times B)\]
natural in $A, B \in \Cat$. For $A=*$, and thus in particular on fibres this equivalence restricts to the identity.
\end{corollary}

\begin{remark}
For a map $p \colon Y \rightarrow A \times B$ in $\LCart(A, B)$, the dual fibration $\SDualco(p)\colon \SDualco(Y)\rt A^{\op}\times B$ is given explicitly as follows: The objects of $\SDualco(Y)$ are the same as those of $Y$ and a morphism from $y$ to $y'$ is given by a span
\[\begin{tikzcd}[row sep=0.4pc, column sep=1.2pc]
& z\ar[ld]\ar[rd] & &\ar[r,mapsto,"p",yshift=-4ex]&\ & & (a,b) \ar[rd,"{(\mathrm{id},f)}"]\ar[ld,"{(g,\id)}"']\\
y && y' &&& (a',b) && (a,b')
\end{tikzcd}\]
where the left map $z\rightarrow y$ is $p$-cartesian. Dually, for a map $q \colon X \rightarrow A \times B$ in $\LCocart(A,B)$ we find $\SDualcart(Y) \rightarrow A^\op \times B$ has morphisms from $y$ to $y'$ given by diagrams
\[\begin{tikzcd}[row sep=0.4pc, column sep=1.2pc]
& z& &\ar[r,mapsto,"q",yshift=-4ex]&\ & & (a',b') \\
y \ar[ru]&& y'\ar[lu] &&& (a',b)\ar[ru,"{(\mathrm{id},f)}"] && (a,b')\ar[lu,"{(g,\id)}"']
\end{tikzcd}\]
where the right map $y' \rightarrow z$ is $q$-cocartesian.
\end{remark}

Taking $B=\ast$, \cref{thm:dualisation functor} reproduces the dualisation equivalence between cartesian fibrations over $A$ and cocartesian fibrations over $A^{\op}$ from \cite{BGN}, which their main result shows is equivalent (naturally in $A$ and $B$) to 
\[\Cart(A) \xrightarrow{\Strcart} \Fun(A^\op,\Cat) \xrightarrow{\Unco} \Cocart(A^\op);\]
We gave a more direct proof of this fact as Theorem \ref{BGNidentification} above (it also follows from the unicity results from \Cref{sec:straightening}). In \cite[Section 2.5]{PartI} we extended this procedure to an equivalence
\[\Dualco \colon \LCart(A,B)  \xrightarrow{\Strcart} \Fun(A^\op,\Cat/B) \xrightarrow{\Unco} \LCocart(A^\op,B) \cocolon \Dualcart\]
by straightening in one variable. We deduce:

\begin{corollary}\label{corthmb}
	The functors 
	\[
	\begin{tikzcd}\SDualco \colon \LCart(A,B) \arrow[r, yshift=0.5ex] & \arrow[l, yshift=-0.5ex]\LCocart(A^\op,B) \cocolon \SDualcart \end{tikzcd} 
	\]
	agree with the equivalences $\Dualco$ and $\Dualcart$ constructed in \cite[Section 2.5]{PartI}.
\end{corollary}

\begin{proof}
Essentially by construction, the two diagrams 
\[
\begin{tikzcd} \LCart(A,B) \ar[r,"\SDualco"] & \LCocart(A^\op,B) && \LCart(A,B) \ar[r,"\Dualco"] & \LCocart(A^\op,B)\\
			\cart(A)/\pr_1 \ar[u,"\mathrm{fgt}"] \ar[r,"\SDualco"] & \cocart(A^\op)/\pr_1 \ar[u,"\mathrm{fgt}"'] && \cart(A)/\pr_1 \ar[u,"\mathrm{fgt}"] \ar[r,"\Dualco"] & \cocart(A^\op)/\pr_1 \ar[u,"\mathrm{fgt}"']\end{tikzcd}
\]
commute naturally in $A$ and $B$. Here the bottom copy of $\SDualco$ arises from \Cref{d:sdual two-var} by taking $B=\ast$. Since the vertical maps are equivalences, the agreement of the two functors in the one-variable case implies that in the two-variable case. 
\end{proof}

\begin{example}\label{arrowdual}
	As an example, consider the dual of the bifibration $(s,t)\colon \Ar(X)\rt X\times X$, see \cite[Corollary 2.4.7.11]{HTT}, which in particular defines an object in $\ortho(X,X)$. We claim that the functor from \Cref{d:sdual two-var} sends this to
	\[\big(\TwL(X) \longrightarrow X^\op \times X\big) \simeq \SDualco(\Ar(X)\longrightarrow X^2).\]

Let us mention that an analogous identification
\[\big(\TwL(X) \longrightarrow X^\op \times X\big) \simeq \Dualco(\Ar(X)\longrightarrow X^2)\] 
is a direct consequence of \cite[Corollary A.2.5]{HMS}, the argument of which is based on the universal property of twisted arrow \icats{} from \cite[Corollary 5.2.1.22]{HA}. In the (current) absence of a similar universal property for twisted arrow \icats{} of \itcats{}, such a proof does not generalise to this more general situation. The present example should be considered a warm-up for Section \ref{sec:OplaxAr} below, where we provide this extension using the span model for dual fibrations.

Now for the proof: Unravelling the definitions, we find that
$\SDualco$ is given by the opposite of the associated \icat{} of a certain simplicial $\infty$-groupoid $\Phi$, where $\Phi_n$ is given by the \igpd{} of diagrams $\phi\colon \TwR([n])\times[1] \rt X$ that take edges in $\TwR([n])^\eg\times\{0\}$ and $\TwR([n])_\ing\times\{1\}$ to equivalences. We now claim that the localisation of $\TwR([n])\times[1]$ at these subcategories is naturally equivalent to $[n] \join [n]^\op$ via the map
	\begin{equation}\label{localisationtwarn}\tag{$\ast$}\TwR([n]) \times [1] \longrightarrow [n] \join [n]^\op, \quad ((i \leq j),\epsilon)\longmapsto \begin{cases} i_l & \epsilon = 0 \\
	j_r & \epsilon = 1 \end{cases},
	\end{equation}
	where we have used subscripts to indicate join factors. Restriction along this map then shows that $\asscat (\Phi) \simeq \TwR(X),$ and so $\asscat (\Phi)^\op \simeq \TwL(X)$ as required. One can also check that this equivalence lives over $X^\op\times X$.

	To see the claim we note that both restrictions $\TwR([n]) \times \{\epsilon\} \rightarrow X$ of a diagram $\phi$ as considered above take all squares in $\TwR([n])$ to pushout squares in $X$ (namely ones in which two opposite edges are equivalences). Using the pointwise formula for left Kan extensions, one readily checks that $\phi\colon \TwR([n]) \times [1] \rightarrow X$ is then left Kan extended from the subposet $\mathcal J_n \times [1]$, where $\mathcal J_n$ is the arch along the top of $\TwR([n])$ consisting of all $(i \leq j)$ with $i=0$ or $j=n$. It follows that the inclusion $\mathcal J_n \times [1] \rightarrow \TwR([n]) \times [1]$ induces an equivalence upon localisation. 
	
	Now note that $\mathcal J_n$ consists of two copies of $[n]$
        glued along the initial vertex. The claim then follows from
        the fact that the localisation of $[n] \times [1]$ at $[n]
        \times \{1\}$ is given by $[n+1]$ and likewise for the
        localisation at $[n] \times \{0\}$; this realises the
        localisation of $\mathcal J_n$ as the pushout of $[n+1]$ and
        $[1+n]$ along $[1]$, embedded into the former as the terminal
        segment, and into the latter as the initial one. Finally, observe that the localisation map $\mathcal J_n \times [1] \rightarrow [n+1] \cup_{[1]} [1+n] \cong [n] \join [n]^\op$ just described is indeed the restriction of \eqref{localisationtwarn}.
\end{example}

We can also use the identification of $\Dualco$ and the functor $\SDualco$ from \Cref{d:sdual two-var} to describe the fibrewise adjoints constructed in \cite[Section 3.1]{PartI} more explicitly.  For this we need to recall some notation from \cite[Section 3]{PartI}. Fix a parametrised left adjoint \[\begin{tikzcd}D \ar[rr,"L"]\ar[rd,"p"'] && C \ar[ld,"q"]\\
	& B, & \end{tikzcd},\]
i.e.\ a map between two cartesian fibrations (though not necessarily preserving cartesian edges), such that the restrictions $L_b\colon D_b\rt C_b$ to the fibres admit right adjoints $R_b\colon C_b \rt D_b$. In \cite{PartI} we constructed from this data a diagram
\[\begin{tikzcd}\Dualco(p) \ar[rd] && \ar[ll,"R"'] \Dualco(q) \ar[ld] \\
	& B^\op, & \end{tikzcd},\]
such that $R$ restricts to the functors $R_b$ under the identifications $\Dualco(p)_b \simeq D_p$ and $\Dualco(q)_b \simeq C_b$ arising from the naturality of $\Dualco$, and showed that this gives an equivalence between the $(\infty,2)$-category of parametrised left adjoints and that of parametrised right adjoints (after taking opposites appropriately). One can use the equivalence $\SDualco \simeq \Dualco$ to give an explicit description of the functor $R$; see also \cite[Section 3.1]{Torii} for a point-set variant of this construction of fibrewise adjoints. 

 To formulate the statement, observe that for each $\beta\colon b'\to b$, the left adjoints $L_b$ and $L_b'$ and the cartesian transport functors of $p$ and $q$ are related by a natural transformation $\lambda_\beta\colon L_{b'}\beta^*\to \beta^*L_b$. We then have:
\begin{proposition}
	For a parametrised left adjoint $L$ the associated parametrised right adjoint 
		\[\begin{tikzcd}\SDualco(p) \ar[rd] && \ar[ll,"R"']\SDualco(q) \ar[ld]\\
		&B^\op, & \end{tikzcd}\]
	can be described as follows:
	\begin{itemize}
		\item For an object $y_b$ in the fibre $\SDualco(q)_b\simeq C_b$ over $b\in B$, one has $R(y_b)\simeq R_b(y_b)$.
		\item For a map $\beta\colon b'\to b$ in $B$, the functor $R$ is given on morphisms over $\beta^{\op}$ in $\SDualco(q)$  by
		\[\begin{tikzcd} [column sep=0.9pc, row sep=0.3pc]			& \beta^{*}x_{b}\ar[ld,two heads]\ar[rd,"f"] & &\ar[rr,mapsto,yshift=-2ex, shorten=-1ex] & & {} & & \beta^*R_{b}(x_b) \ar[rd, shorten >=-0.8ex, "g"{pos=0.4}]\ar[ld,two heads, shorten >=-0.8ex]\\
			x_{b} && y_{b'} &&& & R_b(x_b) && R_{b'}(y_{b'}).
		\end{tikzcd}\]
		Here the left-pointing arrows are cartesian lifts of $\beta$ in $C$ and $D$ and the right-pointing arrows are fibrewise over $b'$, with $g$ given by the (fibrewise) adjoint to
		\begin{equation}\label{eq:leg}\begin{tikzcd}
				L_{b'}\beta^*R_b(x_b)\arrow[r, "\lambda_{\beta}"] & \beta^*L_bR_b(x_b)\arrow[r, "\epsilon"] & \beta^*x_b\arrow[r, "f"] & y_{b'}.
		\end{tikzcd}\end{equation}
	\end{itemize}
\end{proposition}
\begin{proof}
	A direct construction of $R$ can be carried out exactly as in \cite[Theorem 3.1.11]{PartI}, so we will be brief. Applying cocartesian unstraightening to the functor $L$ yields a \lortho{} $\pi\colon X\rt B\times [1]$. Since $L$ admits fibrewise right adjoints, this map is an op-Gray fibration as well. Applying the span dualisation from \Cref{thm:dualisation functor} to this op-Gray fibration yields a \lortho{} $\SDualcart(\pi)\colon \SDualcart(X)\rt [1]\times B^{\op}$, which can be straightened over $[1]$ to yield the desired functor $R$.
	
	Let us now describe the behaviour of $R$ on arrows (and hence also on objects). Note that the functor $R$ arises from cartesian transport in $\SDualcart(\pi)\colon \SDualcart(X)\rt [1]\times B^{\op}$ in the direction of $[1]$. Consequently, it suffices to understand those squares in $\SDualcart(\pi)$ whose vertical maps are cartesian lifts of the map in $[1]$ and whose horizontal maps cover an arrow $\beta^{\op}$ in $B^{\op}$. Like in the proof of \Cref{thm:dualfibs}, unraveling the definitions shows that such a square in $\SDualcart(X)$ corresponds to a diagram in the domain of $\pi\colon X\rt B\times [1]$ of the form
	$$
	\begin{tikzcd}
		R_b(x)                                              & \beta^*R_b(x) \arrow[r, "g"] \arrow[l, two heads]                                                                              & R_{b'}(y) \\
		R_b(x) \arrow[u, "\sim"] \arrow[d] & \beta^*R_b(x) \arrow[d] \arrow[l, two heads] \arrow[u, "\sim"] \arrow[r, "g"] & R_{b'}(y) \arrow[d] \arrow[u, "\sim"'] \\
		x                                              & \beta^*x \arrow[l, two heads] \arrow[r, "f"]                                                              & y
	\end{tikzcd}                                                             
	\quad\longmapsto\quad
	\begin{tikzcd}
		(b, 0)                                              & (b', 0) \arrow[r] \arrow[l, "\beta"{swap}]                                                                              & (b', 0) \\
		(b, 0) \arrow[u] \arrow[d] & (b', 0) \arrow[d] \arrow[l, "\beta"{swap}] \arrow[u] \arrow[r] & (b', 0) \arrow[d] \arrow[u] \\
		(b, 1)                                            & (b', 1) \arrow[l, "\beta"{swap}] \arrow[r]                                                              & (b', 1)                                                              
	\end{tikzcd}                                                             
	$$
	where the top right and bottom left squares are (cartesian) ambigressive in $X$. In particular, the left pointing arrows are all $\pi$-cartesian lifts of $\beta$.  Let us explain the rest of the diagram in more detail.
	For the left and right vertical spans to describe cartesian arrows in $\SDualcart(X)$, one needs their upwards pointing leg to be an equivalence (as indicated) and their downwards pointing leg to define a cartesian arrow for the map $\pi_b\colon X_b\rt \{b\}\times [1]$. Because each $\pi_b$ is both a cocartesian and a cartesian fibration classifying the adjoint pair $(L_b, R_b)$, the objects in the middle row are then given by $R_b(x), \beta^*R_b(x)$ and $R_{b'}(y)$. Finally, the right vertical square is entirely contained in the fibre $X_{b'}$. Since $R_{b'}(y)\rt y$ was $\pi_{b'}$-cartesian, the map $g$ is therefore the \emph{unique} one making the square commute. 
	
	It now remains to verify that we can indeed take $g$ to be the adjoint to \eqref{eq:leg}. This follows from an analysis very similar to \cite[Proposition 3.2.7]{PartI}. In fact, when $f$ is the identity, the adjoint to \eqref{eq:leg} is precisely the \emph{mate} of $\lambda_\beta$; In this case, \cite[Proposition 3.2.7]{PartI} precisely asserts that the mate makes the bottom right square commute.
\end{proof}

Finally, we will use \cref{thm:dualisation functor} to identify various ways of straightening orthofibrations. 

\begin{proposition}\label{compatibility}
	The four equivalences given by
	\[\begin{tikzcd}[row sep=0.2pc]
	\ortho(A,B) \arrow[r, "\Dualco"] & \cocart(A^{\op} \times B) \arrow[r, "\Strco"] & \Fun(A^{\op} \times B, \Cat)\\
	\ortho(A,B) \arrow[r, "\Dualcart"] & \cart(A \times B^{\op}) \arrow[r, "\Strcart"] & \Fun(A^{\op} \times B, \Cat)\\
	\ortho(A,B) \arrow[r, "\Strcart"] & \Fun(A^{\op}, \cocart(B)) \arrow[r, "\Strco"] & \Fun(A^{\op} \times B, \Cat)
	\end{tikzcd}\]
	and 
	\[\begin{tikzcd}
	\ortho(A,B) \arrow[r, "\Strco"] & \Fun(B, \cart(A)) \arrow[r, "\Strcart"] & \Fun(A^{\op} \times B, \Cat)
	\end{tikzcd}\]
	are pairwise equivalent, naturally in $A, B\in \Cat$. 
\end{proposition}

\begin{definition}\label{def:unoc}
	We shall refer to any of the functors above as the \emph{orthocartesian (un)straightening equivalence}, in formulae
	\[\begin{tikzcd}
		\Strort \colon \ortho(A,B) \arrow[r, yshift=0.5ex] & \arrow[l, yshift=-0.5ex] \Fun(A^{\op} \times B,\Cat) \cocolon \Unort.\end{tikzcd}\] 	
\end{definition}

For example we learn from \cref{arrowdual}, that 
\[\Strort\big((s,t) \colon \Ar(C) \rightarrow C \times C\big) \simeq \big(\Map_C \colon C^{\op} \times C \rightarrow \Gpd\big),\]
which we will generalise to (op)lax arrow categories of $(\infty,2)$-categories in the next section.

We shall give a fairly direct comparison between the four equivalences above in the present section, but there could be more ways of straightening an orthofibration. For example, in \cite{Stevenson} Stevenson produced an equivalence
\[\bifib(A,B) \simeq \Fun(A^{\op} \times B,\Gpd)\]
by comparing both sides to a model category of correspondences. In order to settle such coherence questions once and for all, we take another cue from \cite{BGN} and show in the appendix that \emph{any} equivalence as in the statement of \Cref{compatibility} that is natural in the input categories, and restricts to the identity for $A=*=B$ agrees with that above (in an essentially unique fashion). We also include similar statements for curved orthofibrations and bifibrations, in particular settling the comparison with Stevenson's construction.

For the direct proof (and also the naturality of the Yoneda embedding) we need:

	\begin{lemma}\label{cor:strintwovar}
		The three functors 
		\[\cocart(A \times B) \xrightarrow{\Strco_{A \times B}} \Fun(A \times B, \Cat)\]
		\[\cocart(A \times B) \xrightarrow{\Strco_A} \Fun(A, \cocart(B)) \xrightarrow{\Strco_B} \Fun(A \times B, \Cat)\]
		\[ \cocart(A \times B) \xrightarrow{\Strco_B} \Fun(B, \cocart(A)) \xrightarrow{\Strco_A} \Fun(A \times B, \Cat)\]
		are pairwise equivalent.
	\end{lemma}
	
\begin{proof}
We again use the formula 
\[\Unco(F) \simeq \colim\Big(\Phi_F\colon \TwR(A \times B) \xrightarrow{(s,t)} (A \times B) \times (A \times B)^\op 
	\xrightarrow{F \times (A \times B)_{-/}} \Cat\Big)\]
from \cite{GHN}. Using $ \TwR(A \times B) = \TwR(A) \times \TwR(B)$ and $(A \times B)_{-/} = A_{-/} \times B_{-/}$ and the fact that colimits over a product can (naturally in the indexing categories) be computed in two steps we find that this colimit agrees with that of
\[\TwR(B) \longrightarrow \Fun(\TwR(A),\Cat) \xrightarrow{\colim} \Cat\]
where the first functor is curried from $\Phi_F$, i.e.\ it takes $f \colon x \rightarrow y \in \TwR(B)$ to 
\[\TwR(A) \xrightarrow{(s,t)} A \times A^{\op} \xrightarrow{F(-,x) \times A_{-/} \times B_{y/}} \Cat.\]
But the colimit of this functor is naturally equivalent to the unstraightening of $F(-,x) \times B_{y/} \colon A \rightarrow \Cat$, so in total $\Unco(F)$ is naturally identified with
\[\colim \left(\TwR(B) \longrightarrow B \times B^\op \xrightarrow{\Unco_A(F) \times B_{-/}} \Cat\right),\]
which is itself unstraightening over $B$; here we regard $\Unco_A$ as the functor
\[\Fun(A \times B,\Cat) \longrightarrow \Fun(B,\cocart(A)) \longrightarrow \Fun(B,\Cat)/\const_A\]
by forgetting the map to $A$. In total this process identifies
\[\Unco_{A \times B} \colon \Fun(A \times B, \Cat) \longrightarrow \cocart(A \times B)\]
 as the composite
\[\Fun(A \times B,\Cat) \xrightarrow{\Unco_A} \Fun(B,\cocart(A)) \xrightarrow{\Unco_B} \cocart(A \times B)\]
as desired.
\end{proof}

\begin{proof}[Proof of \cref{compatibility}]
Recall that $\Dualco$ is defined as the composite 
\[\ortho(A,B)  \xrightarrow{\Strcart} \Fun(A^\op,\cocart(B)) \xrightarrow{\Unco} \Cocart(A^\op \times B).\]
Thus \cref{cor:strintwovar} immediately identifies the functors 
\[\ortho(A,B) \xrightarrow{\Dualco} \cocart(A^{\op} \times B) \xrightarrow{\Strco} \Fun(A^{\op} \times B, \Cat)\]
and
\[\ortho(A,B) \xrightarrow{\Strcart} \Fun(A^\op, \cocart(B)) \xrightarrow{\Strco} \Fun(A^{\op} \times B, \Cat),\]
and the argument for the first and fourth functors is dual. To compare the first composite with
\[\ortho(A,B) \xrightarrow{\Dualcart} \cart(A \times B^{\op}) \xrightarrow{\Strcart} \Fun(A^{\op} \times B, \Cat)\]
we observe that the diagram
\[\begin{tikzcd} \cocart(A^{\op} \times B) \ar[rd,"\SDualcart"'] \ar[rr,"\SDualcart"] &&  \cart(A \times B^{\op}) \\
& \ortho(A,B) \ar[ru,"\SDualcart"'] & \end{tikzcd}\]
commutes essentially by construction (whereas this does not seem clear from the construction for $\Dualcart$). It now follows from 
\cref{corthmb} that $\Dualco \colon \ortho(A,B) \rightarrow\cocart(A^{\op} \times B)$ agrees with 
\[\ortho(A,B) \xrightarrow{\Dualcart} \cart(A \times B^{\op}) \xrightarrow{\Dualco} \cocart(A^{\op} \times B),\]
whence the result follows from \cref{BGNidentification}.
\end{proof}

\begin{remark}%
	\begin{enumerate}
		\item Straightening of bifibrations is also discussed in detail in \cite[Section 5]{HLAS}, \cite[Appendix A]{HMS} and \cite[Section 7.1]{CDH1}, in each case by choosing one of the last two equivalences from \ref{compatibility} as the definition. For example, the stable \icat{} underlying the Poincar\'e \icat{} $\mathrm{Pair}(\mathcal C,\QF)$ from \cite[Section 7.3]{CDH1} is simply the orthocartesian unstraightening of $\Omega^\infty \mathrm{B}_\QF \colon \mathcal C^\op \times \mathcal C^\op \rightarrow \Gpd$ in the language of the present paper.
\item The equivalence
\[\Gray(A,B)\simeq \Fun(A \boxtimes B,\Cat)\] constructed in \cite[Section 5.2]{PartI} by definition restricts to the composite 
\[\cocart(A \times B) \xrightarrow{\Strco} \Fun(A, \cocart(B)) \xrightarrow{\Strco} \Fun(A \times B, \Cat);\]
see \cite[Remark 5.2.10]{PartI}. As another consequence of (the dual of) Corollary \ref{cor:strintwovar} we obtain that our straightening equivalence for Gray fibrations extends the usual one for cocartesian fibrations over $A \times B$. 
\end{enumerate}
\end{remark}

\section{(Op)lax arrow and twisted arrow \icats{}}\label{sec:OplaxAr}
In this section we will discuss an application of the dualisation procedure from \Cref{thm:dualisation functor} which extends Example \ref{arrowdual}: The duality between the arrow and twisted arrow category of an \icat{} $X$ extends to a duality between the \emph{oplax} arrow and twisted arrow category of an $(\infty, 2)$-category $\mathbf{X}$. More precisely, the oplax twisted arrow category of $\mathbf{X}$
\[
(s,t)\colon\TwR(\mathbf{X})\rt X\times X^\op
\]
is introduced in work of Abell{\'a}n Garc{\'\i}a and Stern \cite{GS} as an explicit model for the cartesian fibration classified by the enriched mapping functor of $\mathbf{X}$ (restricted to its underlying $(\infty, 1)$-category $X$)
\[\Map_{\mathbf{X}}\colon X^\op\times X\rt \Cat.\]
We will show that the enriched mapping functor also classifies the orthofibration
\[(s, t)\colon \Ar^{\oplax}(\mathbf{X})\rt X\times X\]
from the $(\infty, 1)$-category underlying the \emph{oplax arrow
  category}, defined using the (a priori unrelated) Gray tensor product of Gagna, Harpaz and Lanari from \cite{GHL-Gray}.

Since the oplax twisted arrow category and the Gray tensor product are both defined using the model for \itcats{} given by scaled simplicial sets \cite{LurieGoo}, we will start with a minimalistic review of these.
\begin{notation}
Recall that a scaled simplicial set is a pair $(\mathbf{X}, S)$
consisting of a simplicial set $\mathbf{X}$ and a subset $S\subseteq
\mathbf{X}_2$ of $2$-simplices that are called \emph{thin}. The
category $\scSet$ of scaled simplicial sets carries a model structure
in which the cofibrations are the monomorphisms, which is related to
the model category $\markCat$ of categories enriched in marked simplicial sets
(with the categorical model structure) by a Quillen equivalence
$\mf{C}^{\scale}\colon \scSet\leftrightarrows \markCat\cocolon
\mm{N}^{\scale}$ \cite[Theorem 4.2.7]{LurieGoo}.
 We define the \icat{} $\Cat_2$ of \itcats{} as the \icat{} associated to any of these two Quillen equivalent model categories (or any of the other standard models, cf.\ \cite[Theorem 0.0.3]{LurieGoo}).
\end{notation}
\begin{notation}
For a simplicial set $X$, we write $X^\sharp$ for the associated scaled simplicial set in which every $2$-simplex is thin. This determines a left Quillen functor $(-)^{\sharp}\colon \sSet\rt \scSet$, where $\sSet$ is equipped with the Joyal model structure. Its right adjoint sends a scaled simplicial set $(\mb{X}, S)$ to the sub-simplicial set of $X$ spanned by the thin $2$-simplices. At the level of \icats{}, this induces the fully faithful inclusion $\Cat\hookrightarrow \Cat_2$ of $(\infty, 1)$-categories into $(\infty, 2)$-categories, together with its right adjoint sending an \itcat{} $\mathbf{X}$ to its underlying $(\infty, 1)$-category $X$.
\end{notation}

To manipulate simple diagrams in \itcats{}, such as lax commuting squares and triangles, let us recall from \cite{BSP} that a \emph{gaunt $2$-category} is a strict $2$-category whose only invertible $1$- and $2$-cells are the identities. For example, a gaunt $1$-category $A$ gives rise to a gaunt $2$-category $[1]_A$ with 
$$
\Map_{[1]_A}(0, 1)=A, \qquad \Map_{[1]_A}(1, 0)=\emptyset, \qquad \text{and}\qquad \Map_{[1]_A}(i, i)=\ast.
$$
Every gaunt $2$-category determines an object in $\markCat$ by taking the nerves of its mapping categories (with degenerate marking) and this determines a fully faithful functor $\cat{Gaunt}_2\hookrightarrow \Cat_2$ from the (ordinary) category of gaunt $2$-categories into the \icat{} of $(\infty, 2)$-categories, with essential image given by the $0$-truncated objects \cite[Corollary 12.3]{BSP}. To decompose diagrams indexed by gaunt $2$-categories into simpler ones, let us record the following:
\begin{proposition}\label{prop:gaunt pushout}
The inclusion $\cat{Gaunt}_2\hookrightarrow \Cat_2$ preserves coproducts and the following pushouts:
\begin{enumerate}
\item\label{it:tadpole} Pushouts of the form $[1]_A \lt \{i\}\xrightarrow{x} \mathbf{X}$ for $i=0, 1$.  

\item\label{it:factorisation} Pushouts of the form $[2]\xleftarrow{\partial^1} [1]\rto{\gamma} \mathbf{X}$, freely adding a factorisation $\gamma=\alpha\beta$.

\item\label{it:2cell} Pushouts of the form $[1]_{[1]} \lt [1]_{\{i\}}\rto{\gamma} \mathbf{X}$ for $i=0, 1$, where $\gamma\colon x\to y$ is an arrow in $\mathbf{X}$ such that $x$ has no nontrivial incoming arrows and $y$ has no nontrivial outgoing arrows. 
\end{enumerate}
\end{proposition}

\begin{proof}
Let us first give explicit descriptions of these three types of strict pushouts in $\markCat$, for an arbitrary marked simplicial category $\mathbf{X}$. These will show that when $\mathbf{X}$ arises from a gaunt $2$-category, then so does the pushout $\mathbf{X}'$.

Case \ref{it:tadpole}: We only treat the pushout along $\{0\}\rt [1]_A$. The pushout $\mathbf{X}'$ contains $\mathbf{X}$ as a full subcategory, together with a new object $z$, so that for any object $u\in \mathbf{X}$
$$
\Map_{\mathbf{X}'}(u, z)= A\times \Map_{\mathbf{X}}(u, x),\qquad \Map_{\mathbf{X}'}(z, u)=\emptyset, \qquad \Map_{\mathbf{X}'}(z, z)=\ast.
$$ 
Composition is defined by acting on the right factor of $A\times \Map_{\mathbf{X}}(u, x)$.

Case \ref{it:factorisation}: The pushout $\mathbf{X}'$ then freely adds a factorisation $\gamma=\alpha\beta$, where $\alpha$ has source $x$, $\beta$ has target $y$, and they a new object $z$ as common target and source respectively. More precisely $\mathbf{X}'$ has one additional object $z$ and mapping categories
\begin{align*}
\Map_{\mathbf{X}'}(u, v)& =\Map_{\mathbf{X}}(u, v), & \Map_{\mathbf{X}'}(u, z)&=\alpha\cdot \Map_{\mathbf{X}}(u, x),\\
\Map_{\mathbf{X}'}(z, v)&=\Map_{\mathbf{X}}(y, v)\cdot \beta, & \Map_{\mathbf{X}'}(z, z) &= \{\mm{id}_z\}\sqcup \alpha\cdot \Map_{\mathbf{X}}(y, x)\cdot \beta
\end{align*}
for $u, v\in \mathbf{X}$. Here the notation $\alpha\cdot \Map_{\mathbf{X}}(u, x)$ indicates that $\alpha\circ(-)\colon \Map_{\mathbf{X}}(u, x)\rt \Map_{\mathbf{X}'}(u, z)$ is an isomorphism for all $u\in \mathbf{X}$. The composition is then the evident one, using the (formal) relation that $\alpha\beta=\gamma$.

Case \ref{it:2cell}: The pushout $\mathbf{X}'$ has the same objects and mapping categories as $\mathbf{X}$, except for $\Map_{\mathbf{X}'}(x, y)$, which is given by a (homotopy) pushout $[1]\lt \{i\}\rt \Map_{\mathbf{X}}(x, y)$. Note that when $\Map_{\mathbf{X}}(x, y)$ is the nerve of a gaunt $1$-category, this pushout of marked simplicial sets is weakly equivalent to the nerve of a gaunt $1$-category (essentially by Case \ref{it:tadpole}).

It remains to show that the above (strict) pushouts also model the pushout in the \icat{} $\Cat_2$. From the above descriptions, one sees that a Dwyer--Kan equivalence of marked simplicial categories  $\mathbf{X}\rt \mathbf{Y}$ induces a Dwyer--Kan equivalence between the resulting pushouts: It clearly induces weak equivalences on mapping objects, and since we add at most one extra object it remains essentially surjective on homotopy categories. Since $\markCat$ is left proper \cite[Proposition A.3.2.4]{HTT}, this means that the above pushouts are all homotopy pushouts, and the result follows.
\end{proof}

	The \icat{} $\Cat_2$ admits a closed (non-symmetric) monoidal structure given by the (oplax) Gray tensor product $\Gtimes$. The Gray tensor product is induced by a Quillen bifunctor $\Gtimes\colon \scSet\times \scSet\rt \scSet$. The scaled simplicial set $(\mathbf{X},S)\Gtimes (\mathbf{Y},T)$ equals $(\mathbf{X}\times \mathbf{Y},S\Gtimes T),$ where a $2$-simplex $(\sigma,\sigma')$ is in $S\Gtimes T$ if the following conditions hold:
	\begin{enumerate}
	\item $(\sigma,\sigma')\in S\times T$,
	\item $\sigma$ factors through $\Delta^{2}\rightarrow \Delta^{\{1,2\}}$ or $\sigma'$ factors through $\Delta^{2}\rightarrow \Delta^{\{0,1\}}$.
	\end{enumerate} See \cite{GHL-Gray} for the details of this construction. Furthermore, by \cite[Proposition 5.1.9]{PartI} its value on the pair $([m], [n])$ is naturally equivalent to the expected (gaunt) \itcat{}
\begin{equation}\label{eq:gray}\begin{tikzcd}[column sep=1.5pc, row sep=1.5pc]
		00\arrow[r]\arrow[d] & 01 \arrow[d]\arrow[r]& 02\arrow[r, dashed]\arrow[d]& 0n\arrow[d]\\
		10\arrow[r]\arrow[d, dashed]\arrow[ur, Rightarrow, start anchor={[xshift=1ex, yshift=1ex]}, end anchor={[xshift=-1.2ex, yshift=-1.2ex]}] & 11 \arrow[r]\arrow[d, dashed]\arrow[ur, Rightarrow, start anchor={[xshift=1ex, yshift=1ex]}, end anchor={[xshift=-1.2ex, yshift=-1.2ex]}] & 12\arrow[r, dashed]\arrow[d, dashed]\arrow[ur, Rightarrow, start anchor={[xshift=1ex, yshift=1ex]}, end anchor={[xshift=-1.2ex, yshift=-1.2ex]}] &  1n\arrow[d, dashed]\\
		m0\arrow[r]\arrow[ur, Rightarrow, start anchor={[xshift=1ex, yshift=1ex]}, end anchor={[xshift=-1.2ex, yshift=-1.2ex]}] & m1 \arrow[r]\arrow[ur, Rightarrow, start anchor={[xshift=1ex, yshift=1ex]}, end anchor={[xshift=-1.2ex, yshift=-1.2ex]}] & m2\arrow[r, dashed]\arrow[ur, Rightarrow, start anchor={[xshift=1ex, yshift=1ex]}, end anchor={[xshift=-1.2ex, yshift=-1.2ex]}] & mn.
\end{tikzcd}\end{equation}
The internal mapping objects induced by the Gray tensor product via
$$
\Map_{\Cat_2}\big(\mb{A}, \mb{Fun}^\lax(\mb{B}, \mb{X})\big)\simeq \Map_{\Cat_2}(\mb{A}\Gtimes \mb{B}, \mb{X})\simeq \Map_{\Cat_2}\big(\mb{B}, \mb{Fun}^\oplax(\mb{A}, \mb{X})\big)
$$
are by definition the \itcats{} of $2$-functors and (op)lax natural transformations between them, see \cite[Definition 3.9]{HaugsengLax}.

Let us now describe the complete Segal \igpd{} models for the fibrations
$$
(s, t)\colon \Ar^\oplax(\mb{X})\rt X\times X, \qquad \qquad\qquad (s, t)\colon \TwR(\mb{X})\rt X\times X^{\op}
$$
associated to an \itcat{} $\mb{X}$. We will start with the oplax arrow \icat{}.
	\begin{definition}
		Let $\mb{X}$ be a \itcat{}. The \emph{oplax arrow \icat{}} $\Ar^\oplax(\mb{X})$ is the \icat{} underlying the oplax functor \icat{} $\mathbf{Ar}^\oplax(\mb{X})=\mb{Fun}^\oplax\big([1], \mb{X}\big)$.
	\end{definition}

Informally, $\Ar^\oplax(\mathbf{X})$ is the \icat{} with objects given by arrows of $\mathbf{X}$, such that morphisms from $f$ to $g$ are given by oplax commuting squares:
\begin{center}
	\begin{tikzcd}[row sep=1.2pc, column sep=1.7pc]
		x\arrow[r] \arrow[d, "f"'] & x'\arrow[d, "g"] \\
		y \arrow[r] \arrow[ur, Rightarrow, start anchor={[xshift=1ex, yshift=1ex]}, end anchor={[xshift=-1.2ex, yshift=-1.2ex]}] & 					   y'
	\end{tikzcd}
\end{center}
More precisely, $\Ar^\oplax(\mathbf{X})$ can be characterised in terms of the Gray tensor product by the natural equivalence
\[\Map_{\Cat}\big(S,\Ar^{\oplax}(\mathbf{X})\big) \simeq \Map_{\Cat_{2}}\big([1] \Gtimes S ,\mathbf{X}\big).\] 
\begin{example}\label{ex:strict oplax arrow}
If $\mathbf{X}$ is a gaunt $2$-category, i.e.\ a $0$-truncated object in $\Cat_2$, then the above natural equivalence becomes a natural bijection of sets. It follows that $\Ar^{\oplax}(\mathbf{X})$ is a gaunt $1$-category, i.e.\ a strict category without nontrivial isomorphisms, and from the description \eqref{eq:gray} of the Gray tensor product one sees that it coincides with the classical oplax arrow category of $\mathbf{X}$.
\end{example}
\begin{remark}\label{rem:lax arrow cat}
Of course, the \emph{lax arrow category} $\Ar^{\lax}(\mathbf{X})$ is defined by a similar universal property:
$$
\Map_{\Cat}\big(S, \Ar^{\lax}(\mathbf{X})\big)\simeq \Map_{\Cat_{2}}\big(S \Gtimes [1] ,\mathbf{X}\big).
$$
The description of the Gray tensor product of simplices \eqref{eq:gray} shows that for any simplex $[m]$, there is a natural equivalence of gaunt $2$-categories $[1]\Gtimes ([m]^{\op})\simeq \big([m]\Gtimes [1]^{\op}\big)^{\op}$.  By adjunction, this determines a natural equivalence
\begin{equation}\label{eq:laxoplaxarrow}
	\Ar^\oplax(\mathbf{X})^{\op} \simeq \Ar^\lax(\mathbf{X}^{1-\op})
\end{equation}
where $\mathbf{X}^{1-\op}$ has only the directions of the $1$-morphisms inverted. On objects, this equivalence sends an arrow in $\mathbf{X}$ to the opposite arrow in $\mathbf{X}^{1-\op}$.
\end{remark}

\begin{proposition}
	Suppose $\mathbf{X}$ is an $(\infty,2)$-category. Then the functor 
	\[(s,t)\colon \Ar^\oplax(\mathbf{X})\rightarrow X\times X\] 
	is an orthofibration where an edge $\sigma$, given by an oplax square
	\begin{center}
		\begin{tikzcd}[row sep=1.2pc, column sep=1.6pc]
			x\arrow[r] \arrow[d, "f"'] & x'\arrow[d, "g"] \\
		y \arrow[r] \arrow[ur, Rightarrow, "\rho", start anchor={[xshift=1ex, yshift=1ex]}, end anchor={[xshift=-1.2ex, yshift=-1.2ex]}]& 					   y',
		\end{tikzcd}
              \end{center}
              over $\iota X\times X$ is cocartesian if and only if
              $\rho$ is invertible, and similarly an edge over $
              X\times \iota X$ is cartesian if and only if $\rho$ is invertible.
              More formally, a morphism $\sigma:[1]\Gtimes [1] \rightarrow \mathbf{X}$ over $\iota X\times X$ (resp. $ X\times \iota X$) is cocartesian (resp. cartesian) if and only if it factors through the 2-functor $[1]\Gtimes [1]\rightarrow [1]\times [1]$.
\end{proposition}
\begin{proof}
	Let us first show that a square $\sigma\colon [1]\times [1]\rt \mb{X}$ (i.e.~a square $\sigma$ for which $\rho$ is an equivalence), corresponding to an arrow in $\Ar^\oplax(\mathbf{X})$ living over $\iota X\times X$, is cocartesian. To see this, consider the following two equivalent unique lifting problems:
	$$\begin{tikzcd}[row sep=1.2pc]
		{\Lambda_0[2]}\arrow[r]\arrow[d] & \Ar^\oplax(\mb{X})\arrow[d, "{(s, t)}"] & & \Lambda_0[2]\Gtimes [1]\coprod\limits_{\Lambda_0[2]\times \{0, 1\}} [2]\times \{0, 1\}\arrow[r]\arrow[d] & \mb{X}\\
		{[2]}\arrow[r]\arrow[ru, dotted] & X\times X & & {[2]\Gtimes [1],}\arrow[ru, dotted]
	\end{tikzcd}$$
	 where $\Lambda_0[2]$ is the {\icat} $1 \leftarrow 0 \rightarrow 2$. Using the explicit description of the Gray tensor product of simplices of \eqref{eq:gray}, the right vertical map is the inclusion of the sub-\itcat{} of $[2]\Gtimes [1]$ given by
	\[\begin{tikzcd}[row sep=1.3pc]
		10 & 00 & 20 \\
		11 & 01 & 21
		\arrow[from=1-2, to=1-1]
		\arrow[from=1-2, to=1-3]
		\arrow[from=1-1, to=2-1]
		\arrow[from=1-2, to=2-2]
		\arrow[from=1-3, to=2-3]
		\arrow[from=2-2, to=2-3]
		\arrow[from=2-2, to=2-1]
		\arrow[curve={height=12pt}, from=2-1, to=2-3]
		\arrow[Rightarrow, start anchor={[xshift=1ex, yshift=1ex]}, end anchor={[xshift=-1.2ex, yshift=-1.2ex]}, from=2-2, to=1-3]
		\arrow[Rightarrow, start anchor={[xshift=-1ex, yshift=1ex]}, end anchor={[xshift=1.2ex, yshift=-1.2ex]}, from=2-2, to=1-1]
		\arrow[curve={height=-12pt}, from=1-1, to=1-3]
	\end{tikzcd}\]
Now suppose that the morphism $0\leq 1$ in $\Lambda_0[2]$ projects to $\iota X\times X$ and the corresponding map $[1]\Gtimes [1]\rt \mb{X}$ factors over $[1]\times [1]$. This means that in the above diagram, the left oplax square commutes and that $00\rt 10$ is an equivalence. The $\infty$-groupoid of such diagrams is given by the $\infty$-groupoid of 2-functors from the \itcat{}
	\[\begin{tikzcd}[row sep=1.3pc]
		& 00 & 20 \\
		11 & 01 & 21
		\arrow[from=1-2, to=1-3]
		\arrow[from=1-2, to=2-2]
		\arrow[from=1-3, to=2-3]
		\arrow[from=2-2, to=2-3]
		\arrow[from=2-2, to=2-1]
		\arrow[curve={height=12pt}, from=2-1, to=2-3]
		\arrow[Rightarrow, start anchor={[xshift=1ex, yshift=1ex]}, end anchor={[xshift=-1.2ex, yshift=-1.2ex]}, from=2-2, to=1-3]
		\arrow[from=1-2, to=2-1]
	\end{tikzcd}\] 
		to $\mathbf{X}$, where the left triangle commutes. For example this follows by decomposing the 2-category before as a pushout in $\Cat_{2}$ using \Cref{prop:gaunt pushout}, and then observing that the 2-category
		\[\begin{tikzcd}
			\bullet \arrow[rd, ""{name=U, left}] \arrow[d, "\sim"'] &         \\
			\bullet \arrow[Rightarrow, from=U, "{\rotatebox[origin=c]{43}{$\sim$}}"{pos=0.3, yshift=2pt,xshift=-1pt}] \arrow[r]                     & \bullet
		\end{tikzcd}\] 
is equivalent to $[1]$. A similar argument shows that the $\infty$-groupoid of diagrams $[2]\Gtimes [1] \rightarrow \mathbf{X}$ such that the morphism $0\leq 1$ in $[2]$ projects to $\iota X\times X$ and the corresponding map $[1]\Gtimes [1]\rt \mb{X}$ factors over $[1]\times [1]$ is also equivalent to the $\infty$-groupoid of 2-functors from the \itcat{} above to $\mathbf{X}$. The functor
\[
 \Hom_{\Cat_2}\big([2]\Gtimes [1],\mathbf{X}\big)\rt \Hom_{\Cat_2}\big(\Lambda_0[2]\Gtimes [1]\coprod\limits_{\Lambda_0[2]\times \{0, 1\}} [2]\times \{0, 1\},\mathbf{X}\big)
\]
lives over the restriction to $[0,1]\Gtimes [1]$. By the previous argument, taking fibres over $\sigma$ we obtain an equivalence, and so we conclude that $\sigma$ is a cocartesian edge.
	
Next note that for any arrow in $\iota X\times X$ and a lift of its domain, there exists an $(s, t)$-cocartesian lift $\sigma\colon [1]\times[1]\rightarrow\mathbf{X}$ and that lifts of this form are closed under equivalence; This shows both that there are enough $(s,t)$-cocartesian edges over $\iota X\times X$, and that, up to equivalence, all $(s, t)$-cocartesian lifts over $\iota X\times X$ are of the form asserted in the proposition.
	
A dual argument proves the existence and form of $(s, t)$-cartesian arrows over $X\times \iota X$. It follows that $(s, t)$ is a \lortho{} in the sense of \Cref{def:curv ortho}.
	
To show that $(s,t)$ is an orthofibration, consider $\alpha\colon a'\rt a$, $\beta\colon b\rt b'$ and let $f\colon a\rt b$ be an object in the fibre $\Ar^\oplax(\mb{X})_{(a, b)}$. We have to show that the interpolating edge associated to this data is an equivalence. Because interpolating edges are unique, we may equivalently exhibit a commutative square in $\Ar^{\oplax}(\mathbf{X})$ living over the square
\[
\begin{tikzcd}
{(a',b)} \arrow[r, "{(\id,\beta)}"] \arrow[d, "{(\alpha,\id)}"{swap}] & {(a',b')} \arrow[d,"{(\alpha,\id)}"] \\
{(a,b)} \arrow[r, "{(\id,\beta)}"] &	{(a,b'),}
\end{tikzcd}	
\]
 whose bottom left corner is $f$ and whose horizontal and vertical morphisms are cocartesian and cartesian, respectively; indeed in this case the interpolating edge is the identity on $\beta f \alpha$ in $\Ar^{\oplax}(\mathbf{X})$. Such a square is given by the following cube in $\mb{X}$
\[
\begin{tikzcd}[column sep=2pc, row sep=1.2pc]
		a'\arrow[rr,equals ]\arrow[rd]\arrow[dd, "f\alpha"{swap}] & & a'\arrow[dd, dotted, "\beta f \alpha"{pos=0.75}]\arrow[rd, "\alpha"]\\
		& a\arrow[rr, equals]\arrow[dd, "f"{pos=0.25, swap}] & & a\arrow[dd, "\beta f"]\\
		b\arrow[rr, dotted]\arrow[rd, equals] & &  b'\arrow[rd, dotted, "\mathrm{id}"] \\
		& b\arrow[rr, "\beta"{swap}] &  & b'
\end{tikzcd}
\]
in which all faces commute.
\end{proof}

Let us now introduce the model for the enhanced twisted arrow \icat{} given in \cite{GS}, which is constructed via the model of scaled simplicial sets.

\begin{notation}
In the following we will denote the image of $i\in \Delta^n$ under the inclusion  $\Delta^n\subset \Delta^n\join\Delta^{n,\op}$ by $i$ and the image of $j\in \Delta^{n,\op}$ under the inclusion of $\Delta^{n,\op} \subset \Delta^n\join\Delta^{n,\op}$ by $\bar{j}$.
\end{notation}

\begin{definition}
	Consider the cosimplicial object
	\[\begin{tikzcd}
	\kappa_{r}\colon \Del\arrow[r] & \scSet; \quad [n]\arrow[r, mapsto] &  \Delta^n\join_{\scale} \Delta^{n,\op}
	\end{tikzcd}\]
sending each $n$ to the join $\Delta^n\join \Delta^{n,\op}$ where a 2-simplex $\sigma\in \Delta^n\join \Delta^{n,\op}$ is thin if one of the following conditions holds:
	\begin{enumerate}[label=(\roman*)]
		\item $\sigma$ is degenerate,
		\item $\sigma$ is contained totally in either $\Delta^n$ or $\Delta^{n,\op}$,
		\item $\sigma = \Delta^{\{i,j,\bar{k}\}}$ for $i< j \leq k$, or
		\item $\sigma = \Delta^{\{k,\bar{j},\bar{i}\}}$ for $i < j \leq k$.
	\end{enumerate}
	The functor $\kappa_{r}$ induces a functor $\TwR=S_{\kappa_{r}}\colon \scSet\rt \sSet$ from scaled simplicial sets to simplicial sets such that $\TwR(\mathbf{X},S)_n \coloneqq \Hom_{\scSet}(\kappa_r([n]),(\mathbf{X},S))$. We enhance this to a functor with values in marked simplicial sets by declaring an edge in $\TwR(\mathbf{X}, S)$ to be marked if it arises from a totally thin diagram $(\Delta^1\join \Delta^{1, \op})^{\sharp}\rt (\mathbf{X}, S)$ and write
	$$\begin{tikzcd}
	\TwR_+(-) \colon \scSet\arrow[r] & \markSet
	\end{tikzcd}$$
	for the resulting functor. Note that the two inclusions $(\Delta^{n})^\sharp \hookrightarrow (\Delta^n) \join_{\scale} \Delta^{n,\op}\hookleftarrow (\Delta^{n, \op})^\sharp$ induce a map \(\TwR(\mathbf{X}, S)\rightarrow X\times X^{\op}\).
\end{definition}

\begin{theorem}[Abell\'an Garc{\'\i}a, Stern]\label{thm:gstern}
	Let $p\colon (\mathbf{Y}, T)\rt (\mathbf{X}, S)$ be a fibration between fibrant scaled simplicial sets. Then the map 
	\[\begin{tikzcd}
	\TwR_+(\mathbf{Y}, T)\arrow[r] & \TwR_+(\mathbf{X}, S)\times_{X\times X^\op} Y\times Y^{\op}
\end{tikzcd}	\]
 is a fibration between fibrant objects in the cartesian model structure on $(\markSet){/Y\times Y^{\op}}$ (where all arrows in $X\times X^{\op}$ and $Y\times Y^{\op}$ are marked). 
 
For each fibrant scaled simplicial set $(\mathbf{X}, S)$, the cartesian fibration $\TwR(\mathbf{X}, S)\rt X\times X^{\op}$ is classified by the functor $\Map_{\mathbf{X}}\colon X^{\op}\times X\rt \Cat$.
\end{theorem}
\begin{proof}
For the first part, it suffices to verify that $\TwR_+(\mathbf{Y}, T)\rt \TwR_+(\mathbf{X}, S)\times_{X\times X^\op} Y\times Y^{\op}$ has the right lifting property with respect to all marked anodyne maps \cite[Proposition B.2.7]{HA}. Let us write $L$ for the left adjoint of $\TwR_+$ and note that for each marked simplicial set $(A, E)$, there is a natural map $A^{\sharp}\amalg A^{\op, \sharp}\rt L(A, E)$, such that mapping into $(\mathbf{X},S)$ we obtain the value of the canonical projection $\TwR(\mathbf{X}, S)\rt X\times X^{\op}$ on $(A,E)$. The desired right lifting property now follows from the fact that for each (generating) marked anodyne map $(A, E)\rt (B, F)$, the map
$$\begin{tikzcd}
L(A, E)\coprod_{A^{\sharp}\amalg A^{\op, \sharp}} \big(B^{\sharp}\amalg B^{\op, \sharp}\big)\arrow[r] & L(B, F)
\end{tikzcd}$$
is a scaled anodyne map \cite[Lemma 2.9 and 2.10]{GS}. 
The second part of the theorem is {\cite[Theorem 3.3]{GS}}.
\end{proof}
\begin{corollary}\label{cor:oplax twist quillen}
The functor $\TwR\colon \scSet\rt \sSet$ is a right Quillen functor, where $\sSet$ is endowed with the Joyal model structure.
\end{corollary}
\begin{proof}
The left adjoint to $\TwR$ preserves cofibrations, which are simply monomorphisms. Furthermore, Theorem \ref{thm:gstern} implies that the map
\[	\TwR_+(\mathbf{Y}, T) \rt \TwR_+(\mathbf{X}, S)\times_{X\times X^\op} Y\times Y^{\op} \rt \TwR_+(\mathbf{X}, S)\] is a composite of fibrations in the cartesian model structure. Since fibrations in the cartesian model structure are in particular categorical fibrations \cite[Proposition 3.1.5.3]{HTT}, we conclude that $\TwR$ sends fibrations between fibrant scaled simplicial sets to categorical fibrations between quasicategories. This suffices by \cite[Proposition E.2.14]{Joyalnotes}.
\end{proof}

\begin{definition}
The right Quillen functor from \Cref{cor:oplax twist quillen} induces a right adjoint functor of \icats{} which we will denote $\TwR\colon \Cat_2\rt \Cat$ and refer to as taking the \emph{oplax twisted arrow \icat{}}.
\end{definition}
To compare with the oplax arrow category, it will be convenient to give a slightly different presentation of $\TwR(\mathbf{X})$ in terms of Segal \igpds{} (which a posteriori does not make use of scaled simplicial sets). To this end, let us make the following observation:
\begin{definition}\label{def:oplax join}
Let us write $[n]\join_{\oplax} [n]^{\op}$ for the gaunt $2$-category informally depicted as
$$\begin{tikzcd}[row sep=1.2pc]
0\arrow[d]\arrow[r] & 1\arrow[r]\arrow[d] & 2\arrow[d]\arrow[r] & \dots\arrow[r] & n\arrow[d]\\
\ol{0} \arrow[r, Rightarrow, shorten=1ex, yshift=1.3pc] & \ol{1}\arrow[l]\arrow[r, Rightarrow, shorten=1ex, yshift=1.3pc] & \ol{2}\arrow[l]\arrow[r, Rightarrow, shorten=1ex, yshift=1.3pc] & \dots\arrow[l]\arrow[r, Rightarrow, shorten=1ex, yshift=1.3pc] & \ol{n}.\arrow[l]
\end{tikzcd}$$
More precisely, the categories $\Map_{[n]\join_{\oplax} [n]^{\op}}(i, j)$ and $\Map_{[n]\join_{\oplax} [n]^{\op}}(\ol{j}, \ol{i})$ are a point if $i\leq j$ and empty otherwise, and $\Map_{[n]\join_{\oplax} [n]^{\op}}(i, \ol{j})$ is the poset of integers $\max(i, j)\leq k\leq n$, ordered by size. We will refer to these integers as \emph{heights}. The composition maps preserve these heights. One readily verifies that this defines a functor $\Del\rt \cat{Gaunt}_2\subseteq \Cat_2$.
\end{definition}
\begin{lemma}\label{lem:oplax join}
For each $n$, there is a natural equivalence between $[n]\join_{\oplax} [n]^{\op}$ and the \itcat{} presented by the scaled simplicial set $\Delta^n\join_{\scale}\Delta^{n, \op}$.
\end{lemma}
\begin{proof}
We will provide a natural weak equivalence of cosimplicial diagrams in $\markCat$
$$\begin{tikzcd}
\phi\colon \mf{C}^{\mathrm{sc}}\big(\Delta^n\join_{\scale}\Delta^{n, \op}\big)\arrow[r] & {[n]\join_{\oplax} [n]^{\op}}
\end{tikzcd}$$
where the target is viewed as a marked simplicially enriched category by taking nerves of the mapping categories. 

Unraveling the definitions, $\mf{C}^{\mathrm{sc}}\big(\Delta^n\join_{\scale}\Delta^{n, \op}\big)$ is given as follows. Its objects are the objects of the poset $0\leq \dots\leq n\leq \ol{n}\leq \dots\leq \ol{0}$. Given two such objects $a, b$, the simplicial set $\Map(a, b)$ is the nerve of the poset of chains $a=x_0\leq x_1\leq \dots\leq x_l=b$, ordered by inclusion. Within (the nerve of) this poset of chains, the following arrows are marked:
\begin{enumerate}[label=(\alph*)]
\item\label{it:half} if $a=i$ and $b=j$ are both contained in the first half or if $a=\ol{i}$ and $b=\ol{j}$ are both contained in the second half, then \emph{all} inclusions of chains are marked.

\item\label{it:height} for $a=i$ and $b=\ol{j}$, let us define the \emph{height} of a chain $\sigma\colon i=x_0\leq x_1\leq \dots\leq x_l=\ol{j}$ to be the largest number $k$ such that either $k$ or $\ol{k}$ is contained in $\sigma$. Then an inclusion of chains is marked if it does not change the height.
\end{enumerate}
The map $\phi$ is then given by the identity on objects and on (nerves of) mapping categories it \ref{it:half} either sends all chains to the point or \ref{it:height} takes the height of a chain. It suffices to verify that these maps between mapping categories are exactly the localisations at the marked arrows. This is evident in case \ref{it:half}, since the poset of chains is either empty or a cube (hence contractible). In case \ref{it:height}, the functor sending a chain to its height has a fully faithful right adjoint, sending each height $k$ to the maximal chain $i\leq i+1\leq \dots \leq k-1\leq k\leq \ol{k}\leq \ol{k-1}\leq \dots \leq \ol{j+1}\leq  \ol{j}$.
\end{proof}
\begin{corollary}
Let $\mathbf{X}$ be an \itcat{}. Then the simplicial $\infty$-groupoid
$$\begin{tikzcd}
{[n]}\arrow[r, mapsto] & \Map_{\Cat_2}\big([n]\join_{\oplax}[n]^{\op}, \mathbf{X}\big)
\end{tikzcd}$$
is a complete Segal $\infty$-groupoid, whose associated \icat{} is naturally equivalent to $\TwR(\mathbf{X})$.
\end{corollary}
\begin{proof}
\Cref{cor:oplax twist quillen} and \Cref{lem:oplax join} identify the
simplicial $\infty$-groupoid $\Map_{\Cat_2}\big([-]\join_{\oplax}[-]^{\op},
\mathbf{X}\big)$ with $\nerve(\TwR(\mathbf{X}))$.
\end{proof}

We will now compute the image of $\Ar^{\oplax}(\mathbf{X})$ under the dualisation functor
\[\SDualcart\colon \ortho(X,X)\rightarrow \Cart(X,X^\op),\] and show it is equivalent to the functor
\[(s,t)\colon\TwR(\mathbf{X})\rightarrow X\times X^\op.\] 

Recall that $\SDualcart$ is given as the composite of $(-)^\op$ and $\ADual$. Since \Cref{rem:lax arrow cat} identifies the opposite of the oplax arrow category of $\mathbf{X}$ with the lax arrow category of $\mathbf{X}^{1-\op}$, we find that the dual of $(s, t)\colon \Ar^{\oplax}(\mathbf{X})\rt X\times X$ is given by
$$\begin{tikzcd}
\SDualcart(s, t)\colon \mathcal{D}\coloneqq \Span^{\perp}(\Ar^\lax(\mathbf{X}^{1-\op}))\arrow[r] & X\times X^{\op}.
\end{tikzcd}$$
Unwinding definitions, we find that $\mathcal{D}$ is the associated \icat{} of the complete Segal \igpd{} given in level $n$ by the sub-\igpd{} of 
\[
\Map_{\Cat_2}\big(\TwR([n])\Gtimes [1], \mathbf{X}^{1-\op}\big)
\]
spanned by the functors \[f\colon \TwR([n])\Gtimes [1] \rt \mathbf{X}^{1-\op}\] which send
\begin{enumerate}[label=(\roman*)]
	\item $\TwR([n])_\ing \Gtimes \{1\}$ to $\core{X^\op}$,
	\item $\TwR([n])^\eg \Gtimes \{0\}$ to $\core{X^\op}$, and 
	\item every 2-morphism $\sigma\subseteq \TwR([n])_\ing \Gtimes [1]$ to an invertible 2-morphism in $\mathbf{X}$.
\end{enumerate} 
This is of course naturally equivalent to the $\infty$-groupoid of functors $f\colon \big(\TwR([n])\Gtimes [1]\big)^{1-\op} \rt \mathbf{X}$ satisfying analogous conditions. To obtain a functor $\mathcal{D}\rightarrow \TwR(\mathbf{X})$ we will write down a natural transformation of cosimplicial objects in $\Cat_2$
\[\begin{tikzcd}
\phi_n\colon \big(\TwR([n])\Gtimes [1]\big)^{1-\op}\arrow[r] &  {[n] \join_{\oplax} [n]^{\op}}.
\end{tikzcd}\] 
We will produce this map at the level of scaled simplicial sets, using
that for a scaled simplicial set $(\mathbf{Y}',T)$ modeling
$\mathbf{Y}$, the opposite $\mathbf{Y}^{1-\op}$ is modelled by
$((\mathbf{Y}')^{\op}, T)$. Let us therefore define 
\[
\phi_n\colon (\TwR([n])\times [1])^{\op}\rt \Delta^n \join \Delta^{n,\op}; \quad ((i \leq j),\epsilon)\longmapsto \begin{cases} \bar{i} & \epsilon = 0 \\
		j & \epsilon = 1 \end{cases}
\]
exactly as in \Cref{arrowdual}. By inspection this defines a map of scaled simplicial sets. For example $\phi_n(\sigma,\sigma') = \Delta^{\{i,j,\bar{k}\}}$ for $i\leq j\leq k$ when $\sigma'$ equals $s_1\colon \Delta^{2}\rightarrow \Delta^{\{0,1\}}$. By definition this is thin in $\Delta\join_{\scale} \Delta^{n,\op}$. Similarly, one can readily show that the image of the other thin simplicies is thin.
We conclude that the $\phi_n$ determine a natural transformation of cosimplicial objects in
$\scSet$.

\begin{example}
Let us describe the induced map $\phi_n\colon \big(\TwR([n])\boxtimes [1]\big)^{\op}\rt [n] \join_{\oplax} [n]^{\op}$ in $\Cat_2$ a bit more precisely in the case
  $n=1$. In this case, the domain is the
  pushout in $\Cat_2$ of two lax squares along a common $1$-morphism, which is a gaunt $2$-category by \Cref{prop:gaunt pushout}. In particular, the domain and codomain of $\phi_1$ are gaunt $2$-categories, and the map $\phi_1$ is then given by the strict $2$-functor 
\begin{equation}\label{eq:DualHomComp}
	\begin{tikzcd}
		00 \arrow[r, "\sim"]  \arrow[rd, Rightarrow, shorten=2ex]               & 01                 & 11 \arrow[l]   \arrow[ld, Rightarrow, shorten=2ex, "\sim"]             &  & \ol{0}                & \bar{1}\arrow[l]                               \\
		\bar{00} \arrow[u] \arrow[r] & \bar{01} \arrow[u] & \bar{11}\arrow[u] \arrow[l, "\sim"]  \arrow[r, mapsto, shift left=4ex, xshift=1.2pc] & {} & 0 \arrow[u]\arrow[r, Rightarrow, shorten=0.8ex, yshift=4ex] \arrow[r] & 1\arrow[u]
	\end{tikzcd}
\end{equation}
obtained by collapsing the $1$- and $2$-morphisms marked by $\sim$. To see that this agrees with the description of $\phi_1$ in terms of scaled simplicial sets, it suffices to note that this strict $2$-functor is uniquely determined by its values on objects and $1$-morphisms.
\end{example}
\begin{remark}
More generally, since $[n]\join_{\oplax} [n]^{\op}$ is a gaunt $2$-category, $\phi_n$ is adjoint to a (strict) functor $\TwR([n])^{\op}\rt \Ar^{\oplax}\big([n]\join_{\oplax} [n]^{\op}\big)$ into its strict oplax arrow category (see \Cref{ex:strict oplax arrow}). This functor sends each $(i\leq j)$ in $\TwR([n])$ to the arrow $j\to \ol{i}$ in $[n]\join_{\oplax} [n]^{\op}$ of minimal height.
\end{remark}
For each $\mathbf{X}$, the maps $\phi_n$ induce a natural transformation of simplicial \igpds{}
\[N\mathcal{D} \rt S_{\kappa_{r}}(\mathbf{X}),\] which after taking associated \icats{} gives a functor $\Phi\colon \mathcal{D}\rightarrow \TwR(\mathbf{X})$ which evidently commutes with the canonical functors to $X\times X^\op$.

\begin{theorem}\label{thm:CartDualAr}
	The functor 
	\[\Phi\colon \mathcal{D} \rt \TwR(\mathbf{X})\] is an equivalence of \icats{}. Consequently, the orthofibration $(s,t)\colon \Ar^\oplax(\mathbf{X}) \rightarrow X \times X$ classifies the mapping \icat{} functor 
	\[\Map_{\mathbf{X}}(-,-)\colon  X^\op\times X\rt \Cat.\]
\end{theorem}

\begin{proof}
Given the first part of the theorem, the second part follows from Corollary \ref{corthmb} and Theorem \ref{thm:gstern}.	Because both sides of the equivalence are the associated \icat{} of Segal \igpds{}, and the functor is induced by a natural transformation of Segal \igpds{}, it suffices to prove that $\Phi_0$ and $\Phi_1$ are equivalences of \igpds{}. 
Note that the map $\phi_0$ is simply the identity on $[1]$, so that $\Phi_0$ is certainly an equivalence.

To see that $\Phi_1$ is an equivalence, note that its domain $\nerve(\mathcal{D})_1$ is the $\infty$-groupoid of $2$-functors from the gaunt $2$-category $(\TwR([1])\times [1])^{\op}$
$$\begin{tikzcd}[column sep=2pc, row sep=1.5pc]
00 \arrow[r, "\sim"]  \arrow[rd, Rightarrow, shorten=2ex]               & 01                 & 11 \arrow[l]   \arrow[ld, Rightarrow, shorten=2ex, "\sim"]                             \\
		\bar{00} \arrow[u] \arrow[r] & \bar{01} \arrow[u] & \bar{11}\arrow[u] \arrow[l, "\sim"]
	\end{tikzcd}$$
sending the marked $1$- and $2$-cells to equivalences in $\mathbf{X}$. It therefore suffices to verify that the functor $\phi_1$ \eqref{eq:DualHomComp} is the universal functor of \itcats{} that collapses these marked morphisms and $2$-morphisms.

To see this, we observe that the above gaunt $2$-category can be obtained by ``cell attachments'', along the lines of the following picture:
$$\begin{tikzcd}
00 \arrow[r, "\sim"]                 & 01                 & 11 \arrow[l] \\
		\bar{00}\arrow[ru, bend right, shorten=0.5ex, ""{name=t1}]\arrow[ru, bend left, shorten=0.5ex, ""{swap, name=s1}] \arrow[Rightarrow, from=s1, to=t1, shorten=-0.3ex]\arrow[u] \arrow[r] & \bar{01} \arrow[u] & \bar{11}\arrow[u] \arrow[l, "\sim"]\arrow[lu, bend left, shorten=0.5ex, ""{name=t2, swap}]\arrow[lu, bend right, shorten=0.5ex, ""{name=s2}]\arrow[Rightarrow, from=s2, to=t2, shorten=-0.3ex, "\sim"]                            
	\end{tikzcd}$$
	Here we start from the middle vertical arrow, attach the bottom horizontal arrows, then the 2-cells and finally the factorisations over $00$ and $11$. Each step is given by a pushout of gaunt $2$-categories as described in \Cref{prop:gaunt pushout} (which is hence also a pushout in $\Cat_2$). 
	
For each pushout where we add a marked cell, collapsing the resulting
cell in $\Cat_2$ simply has the result of not adding a cell at
all. Consequently, the \itcat{} obtained from $(\TwR([1])\times
[1])^{\op}$ by collapsing the marked cells is the gaunt $2$-category
obtained by performing only those cell attachments where no marked
cell was added. The result of this is then
precisely $[1]\join_{\oplax} [1]^{\op}$, and the collapsing map is
precisely $\phi_1$.
\end{proof}

\begin{remark}\label{rmk:universal_cart_fib}
Note that the functor $\Map_{\tCat}(*,-)\colon \Cat \rightarrow \Cat$ is equivalent to the identity on $\Cat$. Therefore we conclude immediately from the previous theorem that the functor
\[
t\colon \Ar^\oplax(\tCat)\times_{\Cat} \{*\}\rt \Cat\]
is a model for the universal cocartesian fibration, i.e.\ the cocartesian fibration which classifies the identity on $\Cat$. The fibre of $s\colon \Ar^\oplax(\tCat)\rightarrow \Cat$ over $*$ is one definition for the oplax slice category $*\oplslice\Cat$, and the previous theorem therefore gives one way to make precise the statement that the forgetful functor $*\oplslice\Cat\rightarrow \Cat$ is the universal cocartesian fibration.

The universal cartesian fibration has a similar description. Indeed, by \Cref{rem:lax arrow cat} the functor
\[
(t^\op,s^\op)\colon \Ar^{\lax}(\mathbf{X})^\op \rt X^\op\times X^\op
\] 
is equivalent to
\[
(s,t)\colon \Ar^\oplax(\mathbf{X}^{1-\op})\rt X^\op\times X^\op.
\]
Again by \Cref{thm:CartDualAr}, this classifies the enhanced mapping functor of $\mathbf{X}^{1-\op}$, which is equivalent to the functor 
\[
X\times X^\op\simeq X^\op\times X\xrightarrow{\Map_{\mathbf{X}}(-,-)} \Cat.
\] 
Applying this in the case $\mathbf{X} = \tCat$, we obtain that the functor 
\[t^\op\colon \Ar^\lax(\tCat)^\op\times_{\Cat^\op} \{*\} \rt X^\op
\]
is a cartesian fibration which classifies the identity on
$\Cat$. Because the fibre of $t\colon \Ar^\lax(\tCat)\rightarrow \Cat$
over $\{*\}$ is one definition for the lax over-category, this makes
precise the statement that $(\ast \laxslice\Cat)^\op \rightarrow
\Cat^\op$ is the universal cartesian fibration.
\end{remark}

\begin{remark}
Let us also briefly discuss the cocartesian fibration encoding the enhanced mapping functor of an \itcat{} $\mathbf{X}$. Just as for \icats{} it is given by a version of the left twisted arrow \icat{}. However, simply taking the opposite of $\TwR(\mathbf{X})$, as one does for \icats{}, would also apply $(-)^{\op}$ to the fibres. Therefore we instead consider the simplicial object $\kappa_{l} = (\kappa_{r})^{\rev,2-\op}$ in $\cat{Cat}_{2}$ obtained from that defining $\TwR$ by reversing the simplicial direction and taking opposites for $2$-morphisms at each level, and
define $\TwL(\mathbf{X})= S_{\kappa_{l}}(\mathbf{X})$. By direct inspection one then finds $\TwL(\mathbf{X}) \simeq \TwR(\mathbf{X}^{2-\op})^{\op}$. In particular, it follows from Theorem \ref{thm:gstern} that
$$(s,t) \colon \TwL(\mathbf{X}) \longrightarrow X^{\op} \times X$$
is a cocartesian fibration, which classifies the composition 
$$X^{\op} \times X \xrightarrow{\Map_{\mathbf X^{2-\op}}} \cat{Cat} \xrightarrow{(-)^{\op}} \cat{Cat}.$$
We claim this composite is the enhanced mapping functor of $\mathbf X$. One way to see this is to note that the functor $(-)^{2-\op}\colon\Cat_2\rightarrow \Cat_2$ is modelled on categories enriched in marked simplicial sets (which we used to define the enhanced mapping functor) by changing the enrichment via  $(-)^{\op}\colon \sSet_{+}\rightarrow\sSet_{+}$. Let us remark that the model of scaled simplicial sets does not seem to admit a simple implementation of reversing $2$-morphisms, so in particular we do not know of a simplicial object in scaled simplicial sets giving rise to $\kappa^{\rev,2-\op}$, and consequently also no explicit model for $\TwL(\mathbf{X})$ as a scaled simplicial set even if $\mathbf{X}$ is given as such. Nevertheless, applying the equivalence of Lemma \ref{lem:oplax join} one can write down a cosimplicial diagram of gaunt $2$-categories equivalent to $\kappa_{l}$.

Note also that $\TwR(\mathbf{X}^{1-\op}) \simeq \TwR(\mathbf{X})$ as cartesian fibrations over $X\times X^\op$, since $(\kappa_{r})^{1-\op}\simeq \kappa_{r}$.
\end{remark}

\begin{remark}
 For the sake of completeness, let us also mention the behaviour of $\Ar^\oplax$ and $\Ar^{\lax}$ under the two kinds of taking  opposites. Recall that we have already seen in \eqref{eq:laxoplaxarrow} that $\Ar^\oplax(\mathbf{X}^{1-\op}) \simeq \Ar^{\lax}(\mathbf{X})^{\op}$. We are left to note \[\Ar^\oplax(\mathbf{X}^{2-\op}) \simeq \Ar^\oplax(\mathbf{X}^{1-\op})^\op \simeq \Ar^{\lax}(\mathbf{X}).\] To see this, observe that both are orthofibrations over $X \times X$ and that the functors $X^{\op} \times X \rightarrow \cat{Cat}$ which both classify are equivalent by Theorem \ref{thm:CartDualAr} and the corresponding statements for the twisted arrow \icats{}.
\end{remark}

\section{Naturality of the Yoneda embedding}\label{sec:Yoneda}

In this short final section we deduce the following result as an application of \Cref{cor:strintwovar}:

\begin{theorem}\label{yonedafunctorial}
	The Yoneda embedding $\cat{A} \rightarrow \mathcal P(\cat{A})$ canonically extends to a natural transformation of functors $\Cat \rightarrow \CAT$  from the inclusion to the composite
	\[\Cat \xrightarrow{\Fun(-^{\op},\Gpd)} (\CAT^\mathrm R)^\op \simeq \CAT^\mathrm{L} \subseteq \CAT.\] 
\end{theorem}

\begin{remark}
As our model of the equivalence $(\CAT^\mathrm R)^\op \simeq \CAT^\mathrm{L}$ we take the equivalence of \cite[Theorem 3.1.11]{PartI}. Intuitively it acts as the identity on objects, and sends a left adjoint $L\colon C\rightarrow D$ to a choice of right adjoint $R\colon D\rightarrow C$.
\end{remark}
\begin{remark}
This question was recently posed to the second author (among others)
by D.~Clausen during a visit to the University of Copenhagen, as it
can be used to simplify a number of arguments in \cite[Section
2]{ClausenJansen} and appears missing from the literature so far. It
is a lucky accident that our methods answer it.

Let us also mention that, building on Theorem \ref{yonedafunctorial}, Ben Moshe and Schlank have recently upgraded the monoidal (and modal) versions of the Yoneda embedding to natural transformations as well, see \cite[Theorem D]{ShayTomer}.
\end{remark}
Before giving a proof of \Cref{yonedafunctorial}, let us briefly comment on related results that already appear in the literature. First, the characterisation of $\mathcal P(\cat{A})$ as the free cocompletion of $\cat{A}$ \cite[Section 5.1.5]{HTT} implies that the assignment $\cat{A} \mapsto \mathcal P(\cat{A})$ extends to a functor equipped with a natural transformation $\cat{A}\rt \mathcal P(\cat{A})$ given pointwise by the Yoneda embedding (essentially by construction). However, it is not a priori clear that this second functor agrees with the one described in Theorem \ref{yonedafunctorial}; let us write $\mathcal P_{\mm{free}}$ and $\mathcal P_{\mm{Kan}}$ to distinguish these two functorialities, the first being the functoriality via the free cocompletion and the second the functoriality of the statement.
	
	Without considering the naturality of the Yoneda embedding, one can produce a natural equivalence $\mathcal P_{\mm{free}}\simeq \mathcal P_{\mm{Kan}}$ as follows: Both of these functors are easily checked to factor as
	\[\Cat \xrightarrow{(-)^\natural} \Cat^\natural \xrightarrow{\mathcal P} \CAT^{\mathrm{ccpt}} \subset \CAT,\]
	where the second term denotes the \icat{} of small idempotent complete \icats{} (and the first functor idempotent completion), and the third term the \icat{} of cocomplete \icats{} admitting a set of completely compact objects that jointly detect equivalences, and functors among them preserving both colimits and completely compact objects. By an argument analogous to \cite[Proposition 5.5.7.8]{HTT}, $\mathcal P_{\mm{free}}$ is an equivalence between the middle two terms. Since $\mathcal P_{\mm{free}}$ and $\mathcal P_{\mm{Kan}}$ are homotopic on individual morphism $\infty$-groupoids \cite[Proposition 5.2.6.3]{HTT}, the same follows for $\mathcal P_{\mm{Kan}}$. But by a minor modification of To\"en's theorem \cite{Toen}, $\Cat^\natural$ has a discrete automorphism $\infty$-groupoid consisting of two objects (the identity and op), whence there is a unique natural equivalence between $\mathcal P_{\mm{free}}$ and $\mathcal P_{\mm{Kan}}$.
	
	To conclude that the Yoneda embedding extends to a natural transformation as in Theorem \ref{yonedafunctorial}, it thus remains to verify that this equivalence is given pointwise by the \emph{identity} of $\mathcal P(\cat{A})$; the trouble is that this latter identification is not clear. 

After this paper appeared, a proof of Theorem \ref{yonedafunctorial} along the lines just sketched was completed by Ramzi, the additional input being a nifty rigidity result for the identity functor on $\Cat$, see \cite[Theorem 2.2]{Ramzi}.

\begin{proof}[Proof of Theorem \ref{yonedafunctorial}]
	Recall that for any \icat{} $\cat{A}$, the Yoneda embedding $\cat{A} \rightarrow \mathcal P(\cat{A})$ is defined as the adjoint of the functor $\Map_{\cat{A}} \colon \cat{A}^\op \times \cat{A} \rightarrow \Gpd$, which in turn is given by the cocartesian unstraightening of $\TwL(\cat{A}) \rightarrow \cat{A}^\op \times \cat{A}$. 
	The definition of the twisted arrow category makes this left fibration functorial in $\cat{A}$. More precisely, we can consider the pullback square
	\[\begin{tikzcd}[column sep=2.2pc]
		E \arrow[r] \arrow[d] & \Ar(\Cat) \arrow[d, "t"] \\
		\Cat \times \Cat \arrow[r, "(-)^\op \times - "] & \Cat
              \end{tikzcd}\]
	and note that $\TwL$ defines a functor from $\Cat$ to the full subcategory of $F\subset E$ spanned by the left fibrations. Note that $F$ is the cartesian unstraightening of the functor $(\cat{A}, \cat{B}) \mapsto \Leftfib(\cat{A}^\op \times \cat{B})$, where the functoriality is in pullback. There are equivalences
	\[\Leftfib(\cat{A}^\op \times \cat{B}) \xrightarrow{\Strco} \Fun(\cat{A}^\op \times \cat{B}, \Gpd) \simeq \Fun(\cat{B},\mathcal P(\cat{A}))\] 
	which are natural in the input \icats{} (where the contravariant functoriality on $\mathcal{P}(-)$ is given by precomposition). Writing $\mathcal F \colon (\Cat^\op)^2 \rightarrow \CAT$ for the functor sending $(\cat{A}, \cat{B})\mapsto \Fun(\cat{B},\mathcal P(\cat{A}))$, we therefore obtain a functor 
	\[\TwL \colon \Cat \longrightarrow F \simeq \Uncart(\mathcal F).\]
	On objects, this takes an \icat{} $\cat{A}$ to its Yoneda embedding and on morphisms it witnesses the lax commutativity of the diagrams
	\[\begin{tikzcd} 
		\cat{A} \arrow[r] \arrow[d,"f"']\arrow[d, Rightarrow, shorten=0.6ex, xshift=1.8pc] & \mathcal P(\cat{A})\\
		\cat{B} \arrow[r] & \mathcal P(\cat{B}) \arrow[u,"f^*"']
	\end{tikzcd}\] given pointwise by the maps $\Map_{\cat{A}}(-,c) \rt \Map_{\cat{B}}\big(f(-),f(c)\big)$.
	It therefore contains all the requisite data for a natural transformation $\cat{A}\rt \mathcal P(\cat{A})$ as in the statement.
	
	To extract the desired natural transformation, curry $\mathcal F$ into a functor
	\[\Cat^\op \rt \Fun(\Cat^\op,\CAT), \quad \cat{B} \longmapsto \big(\cat{A} \mapsto \Fun(\cat{B},\mathcal P(\cat{A}))\big).\]
	Now observe that this functor takes values in the subcategory of $\Fun(\Cat^\op,\CAT)$ spanned by functors to $\CAT^\mathrm{R}$ with left adjointable squares as morphisms. Via cartesian unstraightening this \icat{} is equivalent to $\textsc{Bicart}(\Cat)$, the intersection the two subcategories $\textsc{Cocart}(\Cat)$ and $\textsc{Cart}(\Cat)$ in the over-\icat{} $\CAT/\Cat$. Via cocartesian unstraightening it is then also equivalent to the subcategory of $\Fun(\Cat,\CAT)$ spanned by functors into $\CAT^\mathrm{L}$ and right adjointable squares.
	
	Applying these equivalences to $\mathcal F$ results in a functor $\mathcal G$ in
	\[\Fun(\Cat^\op, \Fun(\Cat,\CAT)) \simeq \Fun(\Cat^{\op} \times \Cat,\CAT),\]
	which still has $\mathcal G(\cat{B}, \cat{A}) \simeq \Fun(\cat{B},\mathcal P(\cat{A}))$, but now the functoriality in the $\cat{A}$-variable is by the left adjoint to restriction, i.e.\ by left Kan extension.
	
	Let $\Unort(\mathcal G)\rt \Cat\times \Cat$ be the orthocartesian unstraightening of $\mathcal G$. We claim that there is a canonical equivalence $\Uncart(\mathcal F) \simeq \Unort(\mathcal G)$. To this end invoke Corollary \ref{cor:strintwovar} to write the cartesian unstraightening functor as 
	\[\Fun(B^\op, \Fun(A^\op,\textsc{Cat})) \xrightarrow{\Uncart} \Fun(B^\op, \textsc{Cart}(A)) \xrightarrow{\Uncart} \textsc{Cart}(A \times B)\]
	and the orthocartesian one as 
	\[\Fun(B^\op, \Fun(A,\textsc{Cat})) \xrightarrow{\Unco} \Fun(B^\op, \textsc{Cocart}(A)) \xrightarrow{\Uncart} \textsc{Ortho}(A \times B).\]
	But now by definition, $\mathcal G$ and $\mathcal F$ have the same image in $\Fun(\Cat^\op,\textsc{Bicart}(\Cat))$, viewed as a subcategory of $\Fun(\Cat^{\op}, \textsc{Cocart}(\Cat))$ and $\Fun(\Cat^{\op}, \textsc{Cart}(\Cat))$ respectively. This gives the claim that $\Uncart(\mathcal F) \simeq \Unort(\mathcal G)$.
	
	The resulting functor 
	\[\TwL \colon \Cat \rightarrow \Unort(\mathcal G)\]
	again takes $\cat{A}$ to its Yoneda embedding, but this time morphisms witness lax commuting squares
	\[\begin{tikzcd} 
		\cat{A} \arrow[r] \arrow[d,"f"'] & \mathcal P(\cat{A}) \arrow[d,"f_!"]\arrow[ld, Rightarrow, shorten=2ex, "\mu"']\\
		\cat{B} \arrow[r] & \mathcal P(\cat{B}).
	\end{tikzcd}\]
	But by \cite[Proposition 5.2.6.3]{HTT} (or rather the second step of its proof), the natural transformation $\mu$ is an equivalence, i.e.\ the diagram actually commutes. 
	
	Armed with this information consider the diagram
	
	\[\begin{tikzcd}
		\Unort(\mathcal G) \arrow[r] \arrow[d] & \Unort(\Fun(-,-)) \arrow[d]\\
		\Cat \times \Cat \arrow[r, "(\mathcal P{,}\mathrm{incl})"] & \CAT \times \CAT
	\end{tikzcd}\]
	whose top right corner is the oplax arrow \icat{} $\Ar^{\oplax}(\CAT)$ of large \icats{} by Theorem \ref{thm:CartDualAr}. Via the inclusion $\Map_{\CAT}(\cat{A},\cat{B}) \subseteq \Fun(\cat{A},\cat{B})$ it contains the actual arrow \icat{} $\Ar(\CAT)$ as a wide subcategory and by the previous observation the composite
	\[\Cat \xrightarrow{\TwL} \Unort(\mathcal G) \longrightarrow \Ar^\oplax(\CAT)\] actually takes values in this subcategory. The resulting functor $\Cat \rightarrow \Ar(\CAT)$ is then the natural transformation we set out to construct.
\end{proof}

\appendix
\section{Uniqueness of straightening and dualisation}\label{sec:straightening}

To resolve coherence questions surrounding the dualisation and straightening equivalences once and for all, we will prove a rigidity result in this appendix, implying that any two ways of straightening or dualising a two-variable fibration are naturally equivalent; namely, we compute the automorphism groups of the functors
\[\begin{tikzcd}[row sep=0.2pc]
	\Funcat\colon \Cat^{\op} \times \Cat^\op \arrow[r] & \CAT, \quad (A,B) \arrow[r, mapsto] & \Fun(A \times B, \Cat), \\
	\Funcat_{\sS}\colon \Cat^{\op} \times \Cat^\op \arrow[r] & \CAT, \quad (A,B) \arrow[r, mapsto] & \Fun(A \times B, \sS)\\
	\Funcat_{\Gray}\colon \Cat^\op \times \Cat^\op\arrow[r] & \CAT; \quad (A, B)\arrow[r, mapsto] & \Fun\big(A\Gtimes B, \tcat{Cat}\big), 
\end{tikzcd}\]
following ideas of \cite{BGN}, who treated the case of a single variable. Here $\Gtimes$ denotes the Gray tensor product (see Section \ref{sec:OplaxAr}) and $\Fun(A\Gtimes B,\tCat)$ denotes the $\infty$-category of functors from the $(\infty,2)$-category $A\Gtimes B$ to the $(\infty,2)$-category of \icats{}. By \cite[5.2.9]{PartI} the functor $\Funcat_{\Gray}$ is equivalent to the functor $\Gray(-,-)$ which sends a pair $(A,B)$ to the $\infty$-category of Gray fibrations over $A\times B$ (see Remark \ref{rem:gray}).

Before we describe the result of our calculations, note that the functor $\Funcat$ carries a natural action of the group $\Aut(\Cat)$, acting pointwise on $\Fun(A\times B, \Cat)$ by postcomposition. When evaluated at $([0], [0])\in \Cat^{\times 2}$, this simply gives the canonical action of $\Aut(\Cat)$ on $\Cat$. By a theorem of To\"en,  $\Aut(\Cat)$ is discrete with two path components, corresponding to the identity and $(-)^\op \colon \Cat \rightarrow \Cat$ \cite{Toen}. Our goal will be to prove:

\begin{theorem}\label{discrete autos}
	Acting by postcomposition and evaluation at $([0],[0]) \in \Cat^{\times 2}$ determines inverse equivalences
	\[\begin{tikzcd}
		\Aut(\Cat) \arrow[r, yshift=0.5ex] & \arrow[l, yshift=-0.5ex] \Aut(\Funcat) .
	\end{tikzcd}\]
	In particular, $\Aut(\Funcat)$ is discrete with $\pi_0\mm{Aut}(\Funcat) = \mathbb Z/2$ having its nontrivial element induced by postcomposition with $(-)^\op \colon \Cat \rightarrow \Cat$. Furthermore, $\Aut(\Funcat_{\sS})\simeq \ast$ and $\Aut(\Funcat_{\Gray})\simeq \ast$.
\end{theorem}

\begin{remark}
Orthocartesian (un)straightening restricts to a natural (un)straightening equivalence $\Strbi\colon \bifib \rt \Funcat_{\sS}$ for bifibrations by \cite[2.3.15]{PartI}. Theorem \ref{discrete autos} implies that the equivalence $\Strbi$ in fact agrees with the equivalence from \cite{Stevenson}, provided the latter is natural. This can be verified using the techniques of \cite[Appendix A]{GHN}; we refrain from working out further details, as we do not need the statement.	
\end{remark}

 \begin{remark}\label{remarksaboutbarwickII}
As mentioned the analogue of \cref{discrete autos} for the functor $A \mapsto \Fun(A,\Cat)$ is one step in the proof of \cite[Theorem 1.4]{BGN} and our proof below follows their strategy in the large. There is, however, one crucial difference: Barwick, Glasman and Nardin use Yoneda's lemma to deduce that the automorphisms of the functor \[\Cat^\op \rightarrow \GPD, \quad A \mapsto \Map_{\CAT}(A,\Cat)\] are given by $\mm{Aut}(\Cat)$. This suffices to establish the analogue of \Cref{discrete autos} after taking \igpd{} cores, but it is unclear to us how to obtain the actual statement from this information. 

By contrast, our more elaborate proof of \Cref{discrete autos} below also makes use of intermediate steps of To\"en's results (and the strategy applies equally well in the situation of \cite{BGN}).
\end{remark}

In the proof we will make use of two tangential results. For the first, recall that $f\colon A\to B$ is called a \emph{faithful} functor if each $\Map_A(a, b)\to \Map_B(f(a), f(b))$ is an inclusion of path components, or equivalently if $\ho f$ is faithful and the commutative square
\[
  \begin{tikzcd}
    A \arrow{r}{f} \arrow{d}& B \arrow{d} \\
    \ho A \arrow{r}{\ho f} & \ho B
  \end{tikzcd}
\]
is cartesian.
\begin{lemma}\label{inclusionfull}
	If $f \colon \cat{A} \rt \cat{B}$ is faithful then the diagram
	\[\begin{tikzcd}
		\Fun(\cat{C},\cat{A}) \arrow[r] \arrow[d] & \Fun(\cat{C},\cat{B}) \arrow[d] \\
		\Fun(\ho\cat{C},\ho\cat{A}) \arrow[r] & \Fun(\ho\cat{C},\ho\cat{B})
	\end{tikzcd}
	\]
	is cartesian and, in particular, $f_* \colon \Fun(\cat{C},\cat{A}) \longrightarrow \Fun(\cat{C},\cat{B})$ is again faithful for any $\cat{C} \in \Cat$. If, furthermore, the restriction $\core(\cat{A}) \rightarrow \core(\cat{B})$ is an inclusion of path components, then so is \[\core(\Fun(\cat{C},\cat{A})) \longrightarrow \core(\Fun(\cat{C},\cat{B})).\]
\end{lemma}
\begin{remark}
	Note that the functor
	\[\ho\Fun(\cat{C},\cat{A}) \longrightarrow \Fun(\ho\cat{C},\ho\cat{A})\]
	is not usually an equivalence, so the square of the lemma does not agree with the square
	\[\begin{tikzcd}
	\Fun(\cat{C},\cat{A}) \arrow[r] \arrow[d] & \Fun(\cat{C},\cat{B}) \arrow[d] \\
	\ho\Fun(\cat{C},\cat{A}) \arrow[r] & \ho\Fun(\cat{C},\cat{B})
	\end{tikzcd}\]
	which is cartesian if and only if $f_*$ is faithful. 
\end{remark}
\begin{proof}[Proof of \cref{inclusionfull}]
	The first assertion is immediate from $\Fun(\ho \cat{C}, \ho\cat{B}) \simeq \Fun(\cat{C},\ho\cat{B})$ and the analogous assertion for $\cat{A}$ in place of $\cat{B}$. The second assertion then follows, since the lower horizontal functor is clearly faithful and faithful functors are closed under pullback: This is most easily seen from the characterisation that all induced maps on morphism complexes have empty or contractible fibres, which is evidently stable under pullback. The third statement similarly follows from the analogue for ordinary categories by applying cores to the diagram of the lemma.
\end{proof}
In the following, recall that a full subcategory $I\hookrightarrow C$ is said to be \emph{dense} if the restricted Yoneda embedding $C\rt \mathcal P(I)$ is fully faithful, or equivalently, if for each $c\in C$ the canonical map $\colim_{i\in I/c} i\to c$ in $C$ is an equivalence.
\begin{lemma}\label{lem:free functors}
Let $\mathbf{X}$ be an \itcat{} and $\Fun(\mathbf{X}, \mathbf{Cat})$ the \icat{} of $2$-functors $\mathbf{X}\rt \mathbf{Cat}$. For each $x\in \mathbf{X}$ and $A\in \Cat$, there exists an object $h_{x, A}\in \Fun(\mathbf{X}, \mathbf{Cat})$  with the universal property
\begin{equation}\label{eq:yoneda universal}
\Map_{\Fun(\mathbf{X}, \mathbf{Cat})}\big(h_{x, A}, F\big)\simeq \Map_{\Cat}(A, F(x)).
\end{equation}
The full subcategory $\Del_{\mathbf{X}}\subseteq \Fun(\mathbf{X}, \mathbf{Cat})$ spanned by all $h_{x, [k]}$ with $x\in \mathbf{X}$ and $k\geq 0$ then forms a dense subcategory.
\end{lemma}
\begin{observation}\label{obs:free functors LKE}
For any functor $f\colon \mathbf{X}\rt \mathbf{Y}$, let $f_!\colon \Fun(\mathbf{X}, \mathbf{Cat})\rt \Fun(\mathbf{Y}, \mathbf{Cat})$ be the left adjoint to the restriction functor. By the universal property \eqref{eq:yoneda universal}, this sends $h_{x, A}$ to $h_{f(x), A}$ and hence restricts to $f_!\colon \Del_{\mathbf{X}}\rt \Del_{\mathbf{Y}}$.
\end{observation}
\begin{proof}
Let us start by observing that by the enriched Yoneda lemma  \cite[Proposition 6.2.7]{Hin}, these universal functors $h_{x, A}$ indeed exist and are given by
\begin{equation}\label{eq:yoneda}\begin{tikzcd}
h_{x, A}=\Map_{\mathbf{X}}(x, -)\times A\colon \mathbf{X}\arrow[r] & \mathbf{Cat}.
\end{tikzcd}\end{equation}
Moreover, Hinich's results show that the \icat{} $\Fun(\mathbf{X}, \mathbf{Cat})$ can be described as that of $\mathbf{X}$-modules in $\Cat$ (in the sense of algebras for a many-object module $\infty$-operad whose underlying many-object algebra is the $\Cat$-enriched \icat{} $\mathbf{X}$). From this description it follows that $\Fun(\mathbf{X}, \mathbf{Cat})$ admits all colimits, and so the inclusion $j\colon \Del_{\mathbf{X}}\hookrightarrow \Fun(\mathbf{X}, \mathbf{Cat})$ induces an adjoint pair
$j_!\colon \mathcal{P}(\Del_{\mathbf{X}}^{\op})\leftrightarrows \Fun(\mathbf{X}, \mathbf{Cat})\cocolon j^*$
and we have to verify that the counit map 
\begin{equation}\label{eq:can colimit}\begin{tikzcd}
j_!j^*F\simeq \colim\limits_{h_{x, [k]}\in \Del_{\mathbf{X}}/F} h_{x, [k]}\arrow[r] & F
\end{tikzcd}\end{equation}
is an equivalence for all $F$. We will prove this in increasing levels of generality.

First, suppose that $F=h_{y, A}$. In this case, we claim that there is a cofinal functor $h_y\colon \Del/A\rt \Del_{\mathbf{X}}/h_{y, A}$ sending each $[k]\to A$ to $h_{y, [k]}\rt h_{y, A}$. The result then follows from the fact that $A\longmapsto h_{y, A}$ preserves colimits and that $A\simeq \colim_{[k]\in \Del/A} [k]$. To see that $h_y$ is indeed cofinal, it suffices to verify that the natural map
$$\begin{tikzcd}
\colim_{[k]\in \Del/A} \Map_{\Fun(\mathbf{X}, \mathbf{Cat})}\big(h_{x, [n]}, h_{y, [k]}\big)\arrow[r] & \Map_{\Fun(\mathbf{X}, \mathbf{Cat})}\big(h_{x, [n]}, h_{y, A}\big)
\end{tikzcd}$$
is an equivalence. From the universal property \eqref{eq:yoneda universal} and the explicit formula \eqref{eq:yoneda} for $h_{y, A}$, one sees that this map is equivalent to
$$\begin{tikzcd}
\colim\limits_{[k]\in \Del/A} \Map_{\Cat}\big([n], \Map_{\mathbf{X}}(y, x)\times [k]\big)\arrow[r] & \Map_{\Cat}\big([n], \Map_{\mathbf{X}}(y, x)\times A\big).
\end{tikzcd}$$
Taking the constant factor $\Map_{\Cat}\big([n], \Map_{\mathbf{X}}(y, x)\big)$ out of the colimit, one sees that this is an equivalence since $\Del\hookrightarrow \Cat$ is a dense subcategory \cite{Rezk}.

Next, note that the class of $F$ for which \eqref{eq:can colimit} is an equivalence is closed under all colimits that are preserved by $j^*$. For example, one easily sees from the universal property \eqref{eq:yoneda universal} that $j^*$ preserves all coproducts.

To conclude, we now write each functor $F$ as the geometric realisation of a bar construction. More precisely, let us fix a set $S$ and an essentially surjective functor $f\colon S\rt \mathbf{X}$. The description of $\Fun(\mathbf{X}, \Cat)$ as modules then gives a free/forgetful adjunction
\begin{equation}\label{diag:free functor adj}\begin{tikzcd}
f_!\colon \prod_S \Cat \arrow[r, yshift=0.5ex] & \Fun\big(\mathbf{X}, \mathbf{Cat}\big)\cocolon f^*\arrow[l, yshift=-0.5ex]
\end{tikzcd}\end{equation}
where $f^*$ evaluates at each object in $S$. Note that $f^*$ preserves all limits and colimits and detects equivalences, so that $(f_!, f^*)$ is a monadic adjunction. The left adjoint $f_!$ sends each $(A_s)_{s\in S}\in \Cat$ to $\coprod_{s\in S} h_{s, A_s}$.

Any $F\colon \mathbf{X}\rt \mathbf{Cat}$ can now be written as the geometric realisation $F\simeq |B_\bullet|$ where $B_\bullet=\mm{Bar}(f_!, f^*f_!, f^*F)$ is the usual bar resolution of $F$ for the monadic adjunction \eqref{diag:free functor adj}. The functor $j^*$ preserves these geometric realisations, i.e.\ for any $x\in \mathbf{X}$ and $k\geq 0$, the map $\big|\Map_{\Cat}([k], B_\bullet(x)\big)\big|\rt \Map_{\Cat}\big([k], F(x)\big)$ is an equivalence. Indeed, up to equivalence we may take $x=f(s)$, in which case this follows from the augmented simplicial object $f^*B_\bullet\to f^*F$ having extra degeneracies \cite[Example 4.7.2.7]{HA}. Since each term in the bar construction is a coproduct of functors of the form $h_{f(s), A}$, for which we have already verified that the counit map \eqref{eq:can colimit} is an equivalence, it follows that \eqref{eq:can colimit} is an equivalence for $F$ as well. 
\end{proof}

\begin{proof}[Proof of Theorem \ref{discrete autos}] 
	We will spell out the case of $\Funcat$ and discuss the modifications required for $\Funcat_{\sS}$ and $\Funcat_{\Gray}$ at the end of the proof.
	
	\subsubsection*{Step 1: reducing to simplices} 
	Recall that precomposition with the localisation $\asscat\colon \sGpd \rightarrow \Cat$ gives a fully faithful embedding
	\[\Fun(\Cat^\op \times \Cat^\op,\cat{X}) \longrightarrow \Fun(\sGpd^\op \times \sGpd^\op,\cat{X})\]
	for any \icat{} $\cat{X}$. Since $\asscat$ preserves colimits, this inclusion furthermore preserves the property of commuting with small limits in each variable separately, and on the right the full subcategory spanned by such functors is equivalent to $\Fun(\Del^\op \times \Del^\op,\cat{X})$. Since $\Funcat \colon \Cat^\op \times \Cat^\op \rightarrow \CAT$ does indeed preserve small limits in each variable, we find that $\Aut(\Funcat)$ agrees with the automorphisms of the (large) bisimplicial \icat{} $\Funcat_{|\Del^2}$ given by $(n,m) \mapsto \Fun([n] \times [m],\Cat)$, with functoriality arising from restriction.

\subsubsection*{Step 2: restricting to generators}
Note that $\Funcat_{|\Del^2}$ takes values in the subcategory $\cat{Pr^R}$ of $\CAT$ consisting of presentable \icats{} and right adjoint functors. Although the inclusion $\cat{Pr^R} \subseteq \CAT$ is not fully faithful, it does induce a fully faithful map on \igpd{} cores, and by Lemma \ref{inclusionfull} this feature persists to simplicial objects. We may therefore compute the automorphisms of $\Funcat_{|\Del^2}$ as a bisimplicial object in $\cat{Pr^R}$ instead. 
	
	From the equivalence $\cat{Pr^R} \simeq (\cat{Pr^L})^\op$, given by taking adjoints, we find
	\[\Fun(\Del^\op \times \Del^\op,\cat{Pr^R}) \simeq \Fun(\Del \times \Del,\cat{Pr^L})^\op.\]
	By another application of Lemma \ref{inclusionfull}, it follows that $\Aut(\Funcat_{|\Del^2})\simeq \Aut(\mathcal G)$ where $\mathcal G$ denotes the (large) bicosimplicial \icat{} \[(n,m) \longmapsto \Fun([n] \times [m],\Cat)\]
	with functoriality arising by left Kan extension. 
	
Now consider the dense subcategories $\mathcal{H}(n, m)=\Del_{[n]\times [m]}\subseteq \mathcal{G}(n, m)$ from \Cref{lem:free functors}. By \Cref{obs:free functors LKE}
every induced map on homotopy categories $\ho\mathcal G(n,m) \rightarrow \ho\mathcal G(n',m')$ carries $\ho \mathcal H(n,m)$ into $\ho \mathcal H(n',m')$. In fact, 
$\mathcal H(n, m)$ coincides with the full subcategory spanned by the images of $\Del\subseteq \Cat=\mathcal{G}(0, 0)$ under all structure maps $\mathcal{G}(0, 0)\to \mathcal G(n, m)$. Since there are pullbacks
	\[\begin{tikzcd} \mathcal H(n,m) \arrow[r]\arrow[d] & \mathcal G(n,m) \arrow[d] \\
		\ho\mathcal H(n,m) \arrow[r] & \ho\mathcal G(n,m)
	\end{tikzcd}\]
	it follows from the functoriality of pullbacks that the $\mathcal H(n,m)$ assemble into a functor $\cat{\Del}^2 \rightarrow \Cat$ equipped with a natural transformation $\mathcal H \rightarrow \mathcal G$. 
	
	Let us now consider the induced maps
	\[\Map_{\cat{c^2Pr^L}}(\mathcal G, \mathcal G) \longrightarrow \Map_{\cat{c^2CAT}}(\mathcal H, \mathcal G) \longleftarrow \Map_{\cat{c^2Cat}}(\mathcal H, \mathcal H).\]
	We claim that both are inclusions of path components. To see this, note that the universal property of $\mathcal P(\mathcal{H}(n, m))$ as the free cocompletion of $\mathcal{H}(n, m)$ produces a natural transformation $\mathcal P(\mathcal{H})\rt \mathcal G$. By \Cref{lem:free functors}, this exhibits each $\mathcal G(n,m)$ as a (left Bousfield) localisation of $\mathcal P(\mathcal H(n,m))$ for each $(n,m) \in \Del^2$. Consequently, writing out the terms as $\infty$-groupoids of natural transformations and applying \cite[Proposition 5.1]{GHN}, we find 
	\begin{align*}
		\Map_{\cat{c^2Pr^L}}(\mathcal G, \mathcal G) &\simeq \lim \Map_{\cat{Pr^L}}(\mathcal G(n,m),\mathcal G(n',m')) \\
		&\subseteq \lim\Map_{\cat{Pr^L}}(\mathcal P(\mathcal H(n,m)),\mathcal G(n',m'))\\
		&\simeq \lim\Map_{\CAT}(\mathcal H(n,m),\mathcal G(n',m')) \\
		& \simeq \Map_{\cat{c^2CAT}}(\mathcal H, \mathcal G)
	\end{align*}
	where the limits run over $[f \colon (n,m) \rightarrow (n',m')] \in \TwR(\Del \times \Del)$. The second term is a set of path components in the third; indeed, this is so before taking limits, so that the fibre over a point in the target is the limit of a diagram only taking values $\emptyset$ and $\ast$, and thus also either empty or contractible itself. Similarly, we find
	\begin{align*}
		\Map_{\cat{c^2Cat}}(\mathcal H, \mathcal H) &\simeq \lim \Map_{\Cat}(\mathcal H(n,m),\mathcal H(n',m')) \\
		&\subseteq \lim \Map_{\CAT}(\mathcal H(n,m),\mathcal G(n',m' )) \\
		&\simeq \Map_{\cat{c^2CAT}}(\mathcal H, \mathcal G).
	\end{align*}
	Now notice that any automorphism $\varphi$ of $\mathcal G$ preserves the sub-diagram $\mathcal H$. Indeed, the induced automorphism on $\mathcal G(0,0) = \Cat$ preserves the full subcategory $\Del \subset \Cat$ by \cite[Corollary 4.4.11 \& Proposition 4.4.13]{LurieGoo} and naturality with respect to left Kan extension then implies that $\varphi$ preserves the full subcategory $\mathcal H(n,m) \subset \mathcal G(n,m)$ as well. The inclusions of path components above therefore refine to inclusions of path components
	\[\mathrm{Aut}(\mathcal G, \mathcal G) \subseteq \mathrm{Aut}(\mathcal H, \mathcal H) \subseteq \Map_{\cat{c^2CAT}}(\mathcal H, \mathcal G).\]
	In particular, the claim that $\mathrm{Aut}(\mathcal G, \mathcal G)$ is discrete with two components will follow from the analogous statement for $\mathcal H$. 
	
	\subsubsection*{Step 3: automorphism group of the generators}
	Recall that $\mathcal H(n, m)$ is spanned by the functors $h_{(i, j), [k]}$ for $(i, j)\in [n]\times [m]$ and $k\geq 0$, defined by the universal property \eqref{eq:yoneda universal} or by the explicit formula \eqref{eq:yoneda}. The latter formula reduces to
	\[
	h_{(i, j), [k]}\colon [n] \times [m] \longrightarrow \Cat, \quad (a,b) \longmapsto \begin{cases} [k]       & a \geq i, b \geq j \\ \emptyset & \text{otherwise}.
	\end{cases}\]
	Combined with \eqref{eq:yoneda universal}, one sees that the mapping $\infty$-groupoids between such functors are discrete and that there are equivalences 
	$$\begin{tikzcd}
	{[n]^{\op}\times [m]^{\op}\times \Del}\arrow[r, "\sim"] & \mathcal{H}(n, m); \quad (i, j, [k])\arrow[r, mapsto] & h_{(i, j), [k]}.
	\end{tikzcd}$$
	In particular, the $\infty$-categories $\mathcal{H}(n, m)$ are 0-truncated, i.e.\ equivalent to ordinary categories with discrete core; one can use this to conclude that the above equivalence is natural in $(n, m)$ (which is now a property, rather than a structure), for the obvious functoriality on the left leaving $\Del$ fixed.
	
	Since all $\mathcal{H}(n, m)$ are $0$-truncated, $\mathcal{H}$ is a $0$-truncated object in the \icat{} of bicosimplicial \icats{} and $\Aut(\mathcal{H})$ is discrete. Furthermore any automorphism of $\mathcal H$ induces one on $\mathcal H(0,0) = \Del$, and this restriction determines the entire transformation: The composite 
	\[\begin{tikzcd} {[n]^\op \times [m]^\op \times \Del} \arrow[r, "\varphi_{n,m}"] & {[n]^\op \times [m]^\op \times \Del} \arrow[r] & \Del\end{tikzcd}\]
	is determined by naturality for the codegeneracy map $(n,m) \rightarrow (0,0)$ and the composite
	\[\begin{tikzcd} {[n]^\op \times [m]^\op \times \Del} \arrow[r, "\varphi_{n,m}"] & {[n]^\op \times [m]^\op \times \Del} \arrow[r] & {[n]^\op \times [m]^\op}\end{tikzcd}\]
	by naturality with respect to the boundary maps $(0,0) \rightarrow (n,m)$. Thus 
	\[\mathrm{Aut}(\mathcal H) \simeq \mathrm{Aut}(\Del) = \mathbb Z/2\]
	as desired.
	
	\subsubsection*{The other two cases}
	To conclude, we briefly describe the modifications to be made to the above argument to prove that $\Funcat_{\sS}$ and $\Funcat_{\Gray}$ have trivial automorphism groups. 
	
	For $\Funcat_{\sS}$, we repeat the entire argument, but in Step 2 we use the natural subcategory $\mathcal{H}_{\sS}(n, m)\subseteq \Fun([n]\times [m], \sS)$ of diagrams of $\infty$-groupoids given by left Kan extensions along $(i, j)\colon [0]\rt [n]\times [m]$ of the constant diagram on the point. One then identifies the diagram $\mathcal{H}_{\sS}$, with functoriality by left Kan extension, with the obvious diagram sending $(n, m)\mapsto [n]^{\op}\times [m]^\op$. This has trivial automorphisms.
	
	For $\Funcat_{\Gray}$, one uses the dense subcategories
	\(\mathcal{H}_{\Gray}(n, m)=\Del_{[n]\Gtimes [m]}\subseteq \Fun([n]\Gtimes [m], \tcat{Cat}).\)
	Again, these full subcategories can also be characterised as those spanned by the essential images of $\Del$ under all structure maps $\mathcal G_{\Gray}(0, 0)\to \mathcal G_{\Gray}(n, m)$. Formula \eqref{eq:yoneda} shows that the object $h_{(i, j), [k]}$ in $\mathcal{H}_{\Gray}(n, m)$ is given by the $2$-functor
	\[
		h_{(i, j), [k]}\colon [n] \Gtimes [m] \longrightarrow \tcat{Cat}, \quad (a,b) \longmapsto  \Map_{[n]\Gtimes [m]}\big((i, j), (a, b)\big)\times [k].
	\]
	The explicit description of the Gray tensor product \eqref{eq:gray} (see \cite[Proposition 5.1.9]{PartI}) identifies the mapping \icat{} $\Map_{[n]\Gtimes [m]}((i, j), (a, b))$ with the poset of maximal chains from $(i, j)$ to $(a, b)$ in the grid $[n]\times [m]$, ordered using that for each square, ``right-after-down" is smaller than ``down-after-right". Using this and the equivalence \eqref{eq:yoneda universal} to compute mapping $\infty$-groupoids, one readily sees that each $\mathcal{H}_{\Gray}(n, m)$ is a $0$-truncated \icat{}. This already implies that $\Aut(\mathcal{H}_{\Gray})$ is discrete, so it remains to verify that it has only one component.
	
	To see this, consider the fully faithful functor (which we only need for gaunt $2$-categories)
	$$\begin{tikzcd}
		\Phi \colon \Cat_2\arrow[r] & \cart(\Del)
	\end{tikzcd}$$
	defined as follows: for each \itcat{} $\tcat{C}$ and $[k]\in \Del$, let $\Phi(\tcat{C})_k$ be the \icat{} obtained by applying the monoidal functor $\Map_{\Cat}([k], -)\colon \Cat\rt \sS$ to each mapping object. By naturality in $[k]$, this defines a simplicial diagram of \icats{}, whose underlying simplicial diagram of $\infty$-groupoids of objects is constant; we define $\Phi(\tcat{C})\rt \Del$ to be the cartesian unstraightening of this simplicial diagram of \icats{}. 
	Unraveling the definitions then shows that the codegeneracy map $\mathcal{H}_{\Gray}(n, m)\rt \mathcal{H}_{\Gray}(0, 0)\cong \Del$ can be identified with the cartesian fibration between $0$-truncated \icats{} $\Phi([n]\Gtimes [m])\rt \Del$ (naturally in $[n]$ and $[m]$).
	
	Let us now consider the map of sets $\Aut(\mathcal{H}_{\Gray})\rt \Aut(\mathcal{H}_{\Gray}(0, 0))$ restricting an automorphism $\varphi$ of $\mathcal{H}_{\Gray}$ to the component $\varphi_{0, 0}$. For an automorphism $\psi$ of $\mathcal{H}_{\Gray}(0, 0)$, the fibre $\Aut(\mathcal{H}_{\Gray})_{\psi}$ can be identified with the set of isomorphisms
	$$\begin{tikzcd}
		\mathcal{H}_{\Gray}\arrow[r] & \psi^*\mathcal{H}_{\Gray}
	\end{tikzcd}$$
	of bicosimplicial diagrams in the over-\icat{} $\Cat/\mathcal{H}_\Gray(0, 0)$. By the above discussion, $\mathcal{H}$ actually determines a cosimplicial diagram in the subcategory $\cart(\mathcal{H}_{\Gray}(0, 0))\subseteq \Cat/\mathcal{H}_\Gray(0, 0)$. By Lemma \ref{inclusionfull}, it then suffices to compute the sets of isomorphisms $\mathcal{H}_{\Gray}\rt \psi^*\mathcal{H}_{\Gray}$ within $\cart(\mathcal{H}_{\Gray}(0, 0))$.
	
	Identifying $\mathcal{H}(0, 0)\cong \Del$, it now suffices to compute the set of isomorphisms of bicosimplicial diagrams in $\cart(\Del)$
	$$\begin{tikzcd}
		\Phi\big([-]\Gtimes [-]\big)\arrow[r] & \psi^*\Phi\big([-]\Gtimes [-]\big)
	\end{tikzcd}
	$$
	for all $\psi\in \Aut(\Del)\cong \{\mm{id}, \op\}$. Since $\Phi\colon \Cat_2\rt \cart(\Del)$ is fully faithful (at least on gaunt $2$-categories), this comes down to computing the sets of natural isomorphisms
	$$
	\varphi\colon [n]\Gtimes [m] \rt [n]\Gtimes [m] \qquad\qquad \text{and} \qquad\qquad\varphi'\colon [n]\Gtimes [m] \rt \big([n]\Gtimes [m]\big)^{2\mm{-op}}.
	$$
	Since such isomorphisms are determined by their behaviour on objects, compatibility with the vertex inclusions shows that there is a unique natural $\varphi$ (the identity) and no natural $\varphi'$. We conclude that $\Aut(F_{\Gray})\subseteq \Aut(\mathcal{H}_{\Gray})\simeq \ast$.
\end{proof}

\bibliographystyle{amsalpha}

\end{document}